\documentclass[envcountsect]{svjour3}

\smartqed 
\usepackage{graphicx}
\usepackage[T1]{fontenc}
\usepackage[utf8]{inputenc}
\usepackage[english]{babel}

\usepackage{soul}
\setstcolor{red}
\usepackage{amssymb}
\usepackage{amsmath}

\usepackage[colorlinks=true]{hyperref} 
\usepackage[table,xcdraw]{xcolor}
\usepackage[nameinlink,capitalize]{cleveref}

\usepackage{enumitem}
\usepackage{tikz-cd}
\usepackage{float}

\usepackage{tabularx}

\hypersetup{
    colorlinks,
    linkcolor={blue!80!black},
    citecolor={blue!80!black},
    urlcolor={blue!80!black},
    pdftitle={Conic nonholonomic constraints on surfaces and control systems},
    pdfauthor={Timothée Schmoderer and Witold Respondek},
    pdfsubject={This paper addresses the equivalence problem of conic submanifolds in the tangent bundle of a smooth 2-dimensional manifold. Those are treated as nonholonomic constraints whose admissible curves are trajectories of the corresponding control systems. We deal with this problem under the prism of feedback equivalence of control systems, both control-affine and fully nonlinear. The first main result of this work is a complete description of regular conic submanifolds. We also give equivalence results for a special class of conic submanifolds via the study of the Lie algebra of infinitesimal symmetries of the corresponding control systems. Then, we consider the classification problem of conic submanifolds, which is achieved via feedback classification of nonlinear control system. Our results describe and completely characterise quadratic systems, and include several normal and canonical forms.},
    pdfkeywords={nonlinear control system, feedback equivalence, conic submanifolds, nonholonomic constraint, normal forms, pseudo-Riemannian geometry},
}

\crefformat{equation}{(#2#1#3)}
\crefformat{enumi}{#2#1#3}
\crefmultiformat{equation}{(#2#1#3)}%
{ and~(#2#1#3)}{, (#2#1#3)}{ and~(#2#1#3)}
\usepackage[width=150mm,top=25mm,bottom=25mm,headheight=15pt]{geometry}
\graphicspath{{img/}}

\newcommand{\RR}{\ensuremath{\mathbb{R}}}
\renewcommand{\SS}{\ensuremath{\mathbb{S}}}

\newcommand{\AAA}{\ensuremath{\mathcal{A}}}
\newcommand{\BBB}{\ensuremath{\mathcal{B}}}

\newcommand{\GGG}{\ensuremath{\mathcal{G}}}

\newcommand{\MMM}{\ensuremath{\mathcal{M}}}

\newcommand{\SSS}{\ensuremath{\mathcal{S}}}

\newcommand{\WWW}{\ensuremath{\mathcal{W}}}
\newcommand{\XXX}{\ensuremath{\mathcal{X}}}

\crefformat{equation}{(#2#1#3)}
\crefformat{enumi}{#2#1#3}
\crefmultiformat{enumi}{#2#1#3}{ and ~#2#1#3}{, #2#1#3}{ and~#2#1#3}
\crefname{problem}{problem}{problems}
\Crefname{problem}{Problem}{Problems}

\crefformat{figure}{#2figure~#1#3}
\Crefformat{figure}{#2Figure~#1#3}

\newcommand{\lccref}[1]{\hyperref[{#1}]{\lcnamecref{#1}~\labelcref{#1}}} 

\newcommand\lb[2]{\left[#1,#2\right]}              
\renewcommand\vec[1]{\frac{\partial}{\partial #1}} 
\newcommand\distrib[1]{\mathrm{span}\left\{#1\right\}} 
\newcommand\ad[2]{\mathrm{ad}_{#1}#2}                
\newcommand\adk[3]{\mathrm{ad}_{#1}^{#2}#3}          
\newcommand\dL[2]{\mathrm{L}_{#1}\left(#2\right)}  
\newcommand\dLk[3]{\mathrm{L}_{#1}^{#2}\left(#3\right)\,} 
       
\newcommand\sgn[1]{\mathrm{sgn}\left(#1\right)}      
\newcommand\ann[1]{\mathrm{ann}\left(#1\right)}      

\newcommand\rk[1]{\mathrm{rk}\,#1}        

\newcommand\diff{\mathrm{d}} 

\newcommand\structfunctA{\ensuremath{\rho}}
\newcommand\structfunctB{\ensuremath{\tau}}

\newcommand\SM{\ensuremath{\SSS}}
\newcommand\Sm{\ensuremath{S}}
\newcommand\sm{\ensuremath{s}}

\makeatletter
\newcommand*\owedge{\mathpalette\@owedge\relax}
\newcommand*\@owedge[1]{%
  \mathbin{%
    \ooalign{%
      $#1\m@th\bigcirc$\cr
      \hidewidth$#1\m@th\wedge$\hidewidth\cr
    }%
  }%
}
\makeatother

\DeclareMathOperator\arcsinh{arcsinh}

\newcommand{\ttiny}[1]{\text{\tiny{$#1$}}}

\newcommand\antigamma{\ensuremath{\rotatebox[origin=c]{180}{$\gamma$}}}

\renewcommand{\st}[1]{}

\journalname{myjournal}

\begin{document}

\title{Conic nonholonomic constraints on surfaces and control systems}
\titlerunning{Quadratic control systems}

\author{Timoth\'ee Schmoderer \and Witold Respondek.}
\authorrunning{T. Schmoderer \& W. Respondek}

\institute{%
LMI - UR 3226 - FR CNRS 3335, INSA Rouen Normandie, Avenue de l'Université 76800 St Etienne du Rouvray \at %
\email{timothee.schmoderer@insa-rouen.fr, witold.respondek@insa-rouen.fr} %
}

\date{Received: date / Accepted: date}

\maketitle

\begin{abstract}
This paper addresses the equivalence problem of conic submanifolds in the tangent bundle of a smooth 2-dimensional manifold. Those {are given by a quadratic relation between the velocities and} are treated as nonholonomic constraints whose admissible curves are trajectories of the corresponding control systems{, called quadratic systems}. We deal with \st{this problem}{the problem of characterising and classifying conic submanifolds} under the prism of feedback equivalence of control systems, both control-affine and fully nonlinear. The first main result of this work is a complete description of \st{regular}{non-degenerate} conic submanifolds {via a characterisation under feedback transformations of the novel class of quadratic control-affine systems. This characterisation can explicitly be tested on structure functions defined for any control-affine system and gives a normal form of quadratizable systems and of conic submanifolds}. \st{We also give equivalence results for a special class of conic submanifolds via the study of the Lie algebra of infinitesimal symmetries of the corresponding control systems.} Then, we consider the classification problem of {regular} conic submanifolds {(ellipses, hyperbolas, and parabolas)}, which is \st{achieved}{treated} via feedback classification of {quadratic control-}nonlinear \st{control} systems. {Our classification includes several normal forms of quadratic systems (in particular, normal forms not containing functional parameters as well as those containing neither functional nor real parameters), and, as a consequence, gives a classification of regular conic submanifolds.} \st{Our results completely describe and characterise quadratic systems, and include several normal and canonical forms.}
\end{abstract}

%
%
\setcounter{tocdepth}{2}
\tableofcontents
\section{Introduction}\label{sec:introduction}
Let $\XXX$ be a smooth connected manifold of dimension $n= 2$ (a surface), equipped with local coordinates $x$\st{$=(z,y)$; we choose the order $(z,y)$ to be consistent with some normal forms existing in the literature}. In the tangent bundle $T\XXX$ of $\XXX$, we consider a smooth 3-dimensional submanifold $\SM${, a hypersurface,} given by
\begin{align*}
\SM=\left\{(x,\dot{x})\in T\XXX,\ \Sm(x,\dot{x})=0\right\},
\end{align*}
\noindent
where $\Sm\,:\,T\XXX\rightarrow\RR$ is a smooth {scalar} function satisfying $\rk{\frac{\partial \Sm}{\partial \dot{x}}}(x,\dot{x})=1$ for all $(x,\dot{x})\in\SM$. Two submanifolds $\SM=\{\Sm(x,\dot{x})=0\}\subset T\XXX$ and $\tilde{\SM}=\{\tilde{\Sm}(\tilde{x},\dot{\tilde{x}})=0\}\subset T\tilde{\XXX}$ are {said to be} equivalent if there exists a diffeomorphism $\phi\,:\,\XXX\rightarrow\tilde{\XXX}$ and a smooth nonvanishing function $\delta\,:\,T\XXX\rightarrow\RR$ such that 
\begin{align*}
    \tilde{\Sm}\left(\phi(x),D\phi(x)\dot{x}\right)&=\delta(x,\dot{x})\,S\left(x,\dot{x}\right),
\end{align*}
\noindent
where $D\phi$ is the derivative (tangent map) of $\phi$. This definition implies that for all $(x,\dot{x})\in\SM$ we have $(\tilde{x},\dot{\tilde{x}}) =(\phi(x),D\phi(x)\dot{x})\in\tilde{\SM}$, hence the \st{diffeomorphism}{application} $(\phi,D\phi)$ maps the graph of $\Sm^{-1}(0)$ into that of $\tilde{\Sm}^{-1}(0)$. Equivalence of submanifolds $\SM$ and $\tilde{\SM}$ means simply that implicit underdetermined ordinary differential equations $\Sm(x,\dot{x})=0$ and $\tilde{\Sm}(\tilde{x},\dot{\tilde{x}})=0$ are equivalent; see e.g. \cite[Definition 2]{bogaevsky2014Implicitordinarydifferential}.

It is natural to ask how to characterise and classify submanifolds $\SM$ of certain particular classes, for instance the class of linear submanifolds given by $\Sm_{lin}(x,\dot{x})=\omega(x)\dot{x}=0$ or the class of affine submanifolds given by $\Sm_{aff}(x,\dot{x})=\omega(x)\dot{x}+h(x)=0$, where $\omega$ is a smooth differential 1-form on $\XXX$ and $h$ is a smooth function on $\XXX$. Those questions have been widely studied under the prism of Pfaffian equations (linear and affine) and go back to Pfaff, Darboux, Cartan \cite{pfaff1815Methodusgeneralisaequationes,cartan1931theoriesystemesinvolution,darboux1882problemePfaff}. Although the problem of classification of Pfaffian equations is still open in its full generality, many important results have been obtained for various classes of linear Pfaffian equations (contact and quasi-contact case, Martinet case, singularities, \cite{zhitomirskii1992TypicalSingularitiesDifferential,zhitomirskii1995Singularitiesnormalforms,martinet1970singularitesformesdifferentielles,jakubczyk2001Localreductiontheorems,elkin2012ReductionNonlinearControl}) and of affine Pfaffian equations (dimension two \cite{jakubczyk1990EquivalenceInvariantsNonlinear}, three \cite{respondek1995Feedbackclassificationnonlinear,respondek1998Feedbackclassificationnonlinear}, and arbitrary \cite{zhitomirskii1998Simplegermscorank,elkin2012ReductionNonlinearControl}). 

We will call $\SM_q$ a quadratic, or a conic, submanifold if it is given by 
\begin{align*}
    \Sm_q(x,\dot{x})=\dot{x}^t\textsf{g}(x)\dot{x}+2\omega(x)\dot{x}+h(x)=0, 
\end{align*}
\noindent
with all involved objects being smooth, i.e. for each point $x\in\XXX$, the set $\SM_q(x)$ forms a conic curve in $T_x\XXX$. The image of a quadratic submanifold $\SM_q=\left\{\Sm_q=0\right\}$ under a diffeomorphism $\tilde{x}=\phi(x)$ is a quadratic submanifold $\tilde{\SM}_q=\left\{\tilde{\Sm}_q=0\right\}$, where $\tilde{\Sm}_q(\tilde{x},\dot{\tilde{x}})=\dot{\tilde{x}}^t\tilde{\textsf{g}}\dot{\tilde{x}}+\tilde{\omega}\dot{\tilde{x}}+\tilde{h}$, with $\textsf{g}=\phi^*\tilde{\textsf{g}}$, $\omega=\phi^*\tilde{\omega}$, and $h=\phi^*\tilde{h}$. On the other hand, a submanifold $\tilde{\SM}$ equivalent to a quadratic $\SM_q$ is given by $\phi^*\tilde{\Sm}=\delta\cdot\Sm_q=0$ and $\tilde{\SM}$ need not be quadratic because of the presence of the function $\delta$. The first purpose of this work is to provide a local characterisation of submanifolds $\SM$ that are equivalent to quadratic submanifolds $\SM_q$.
In particular, we will identify submanifolds equivalent to elliptic, hyperbolic, and parabolic conics, given {respectively} by 
\begin{multline*}
    \SM_E=\left\{a^2(\dot{z}-c_0)^2+b^2(\dot{y}-c_1)^2=1\right\},\qquad\SM_H=\left\{a^2(\dot{z}-c_0)^2-b^2(\dot{y}-c_1)^2=1\right\},\\
    \textrm{and}\quad \SM_P=\left\{a\dot{y}^2-\dot{z}+b\dot{y}+c = 0\right\},
\end{multline*}
\noindent
\st{respectively, }where \st{$r$}{$a$, $b$}, $c_0$, $c_1$, \st{$a$, $b$,}and $c$, are smooth functions of $x=(z,y)$, we choose the order $(z,y)$ to be consistent with some normal forms existing in the literature, satisfying \st{$r(\cdot)\neq0$ and $a(\cdot)\neq0$}{$a(\cdot)\neq0$ and $b(\cdot)\neq0$ in the elliptic and hyperbolic cases and $a(\cdot)\neq0$ in the parabolic case}. We will also discuss the case of passing smoothly from one type to another. {To state and discuss general facts on elliptic, hyperbolic, and parabolic submanifolds we set $\SM_Q=\{\SM_E,\SM_H,\SM_P\}$.}

The second goal of this work is to provide a classification of elliptic {$\SM_E$}, hyperbolic {$\SM_H$}, and parabolic {$\SM_P$} conics. {In the elliptic and hyperbolic cases, we will first describe the submanifolds with the functions $a(x)=b(x)$ and call them \emph{conformally-flat}; second, if, additionally, $a=b=1$, then we call them \emph{flat} elliptic or hyperbolic submanifold. Finally, we will describe the forms for which we also have $c_0,c_1\in\RR$ (in the same coordinate system as the one for which we have $a=b=1$), called a \emph{constant-form} elliptic and hyperbolic submanifolds; the particular case of $c_0=c_1=0$ is called a \emph{null-form} to emphasis the absence of any (functional, continuous, or discrete) parameter. In the parabolic case, we will describe the submanifolds with $a=1$, called \emph{weakly-flat}, and characterise those with, additionally, $b=0$ (called \emph{strongly-flat}) and, moreover, with $c\in\RR$ (called \emph{constant-form}); in the particular case where $c=0$ we call the parabolic submanifold a \emph{null-form}}. Our classification is summarised by \cref{tab:classification_nomenclature} in \cref{sec:classification_conic}, where a characterisation of all above subclasses is given. \st{The several normalisations of submanifolds that we propose are reflected in normal and canonical forms of control-nonlinear systems.} \st{We will give several normalisations and, in particular, we will characterise and propose canonical forms of submanifolds $\Sm_q$ with constant coefficients (called strongly flat in the paper).}


Our analysis {of the equivalence problem of submanifolds} will be based on attaching to a submanifold $\SM=\{\Sm(x,\dot{x})=0\}\subset T\XXX$ two control systems. First, 
\begin{align*}
    \Xi_{\SM}\,:\,\dot{x}=F(x,w),\quad x\in\XXX,\quad w\in\WWW\subset\RR,
\end{align*}
\noindent
where $\dot{x}-F(x,w)=0$ is a {regular} parametric representation of $\SM$, {that is for all values of the parameter $w\in\WWW$ (interpreted as a scalar control) we have $\Sm(x,F(x,w))=0$ and $\rk\frac{\partial F}{\partial w}(x,w)=1$,} and second,
\begin{align*}
    \Sigma_{\SM}\,:\,\left\{\begin{array}{rl}
        \dot{x} &= F(x,w)  \\
        \dot{w} &= u 
    \end{array}\right.,\quad (x,w)\in\XXX\times\WWW, \quad u\in\RR, 
\end{align*}
called, respectively, a first and second {extension}\st{prolongation} of $\SM$. Notice that for $\Xi_{\SM}$ the control $w$ enters in a nonlinear way, whereas for $\Sigma_{\SM}$ the control $u$ enters in an affine way, but \st{this costs the}{for the price of} augmentation of the dimension of the state space. Observe that, since \st{$\SM$}{$\Sm(x,\dot{x})=0$, defining $\SM$,} relates the positions $x$ with the velocities $\dot{x}$, it describes a nonholomic constraint. We say that a smooth curve $x(t)\in\XXX$ satisfies the nonholonomic constraint given by $\SM$ if we have $(x(t),\dot{x}(t))\in\SM$. Clearly, $x(t)$ satisfies the nonholonomic constraint described by $\SM$ (equivalently, satisfies the {implicit} differential equation $\Sm(x(t),\dot{x}(t))=0$) if and only if $x(t)$ is a trajectory of $\Xi_{\SM}$ for a certain smooth control $w(t)$ or, equivalently, $(x(t),w(t))$ is a trajectory of $\Sigma_{\SM}$ for a {smooth} control $u(t)$. \st{An}{A crucial} observation that links studying submanifolds $\SM\subset T\XXX$ and their extensions\st{prolongations} $\Xi_{\SM}$ and $\Sigma_{\SM}$ is that the equivalence of submanifolds corresponds to the equivalence of control systems $\Xi_{\SM}$ and $\Sigma_{\SM}$ via feedback transformations, general for $\Xi_{\SM}$, and control-affine for $\Sigma_{\SM}$, as assured by \cref{prop:equiv_of_equiv}.

\begin{figure}[!h]
    \centering
    \begin{tikzpicture}[%
rectnode/.style={above, rectangle, draw=black, thick, inner sep=3pt},
dblarrow/.style={{Implies}-{Implies},double},
]
%
\node[rectnode] (main) at (0,7) {Pbm characterisation $\SM_q$ and classification $\SM_Q$};
\node[rectnode] (mainCh) at (-4.5,5) {Pbm characterisation $\Sigma_q$};
\node[rectnode] (mainCh2) at (-4.5,3.5) {Characterisation $\Sigma_q$};
\node[rectnode] (mainCh3) at (-5.5,2) {Characterisation $\Sigma_Q$};
\node[rectnode] (mainCh4) at (-3.1,0.8) {Normal form of $\Sigma_q$ and $\SSS_q$};
\node[rectnode] (mainCl) at (3.5,5) {Pbm classification $\Xi_Q$};
\node[above left] at (0.9,4.8) {\Cref{subsec:classification_eh}};
\node[rectnode, text width=4cm] (ClNF1) at (1.4,1.9) {\scriptsize Classification $\Xi_E$ and $\Xi_H$:\newline \cref{prop:m1:equivalence_prenormal_elliptic_hyperbolic}: conformally-flat\\ \cref{prop:m1:EH:existence_commutative_frame}: flat\newline \cref{prop:m1:charact_class_eh_flat_systems}: constant-form \& null-form};
\node[rectnode] (mainCl3) at (1.4,0.6) {Classification of $\SM_E$ and $\SM_H$};
\node[above right] at (5.45,4) {\Cref{subsec:classification_p}};
\node[rectnode, text width=3cm] (ClNF2) at (5.5,1.9) {\tiny Classification $\Xi_P$:\newline \cref{cor:m1_prenormal_form}: weakly-flat\\ \cref{thm:m1_normalisation_quadratic_system_gamma_1_0}: strongly-flat, constant-form \& null-form};
\node[rectnode] (mainCl4) at (5.5,0.6) {Classification of $\SM_P$};
%
%
\draw[dblarrow] (main) -- (mainCh);
\draw[dblarrow] (main) -- (mainCl);
\draw[-latex] (mainCh) -- (mainCh2);
\draw[-latex] (mainCh2) -- (mainCh3);
\draw[-latex] (mainCh2) -- (mainCh4);
\draw[-latex] (mainCl) -- (ClNF1);
\draw[-latex] (mainCl) -- (ClNF2);
\draw[-latex] (ClNF1) -- (mainCl3);
\draw[-latex] (ClNF2) -- (mainCl4);
\draw (-7.3,0.5) rectangle (-1.1,6.5); 
\node[above right] at (-7.4,0) {\Cref{sec:feedback-m-1}: Characterisation results};
\draw (-0.9,0.5) rectangle (7.3,6.5); 
\node[above right] at (3,0) {\Cref{sec:classification_conic}: Classification results};
%
\node[fill=white] at (-3,6) {\cref{prop:equiv_of_equiv}};
\node[fill=white,above] at (-4.5,4.4) {\cref{thm:feedback_quadratization_m_1}};
\node[fill=white] at (-5.5,3) {\cref{cor:thm_feedback_quadratization}};
\node[fill=white,above] at (-2.9,1.6) {\cref{thm:m1_normal_form_quadratizable_systems}}; 
\node[fill=white] at (0.9,1.55) {\cref{lem:classification_ehp_system_to_submanifolds}};
\node[fill=white] at (6,1.55) {\cref{lem:classification_ehp_system_to_submanifolds}};
\node[fill=white] at (3,6) {\cref{prop:equiv_of_equiv}};
%
%
\end{tikzpicture}
    \caption{{Walk-through the paper with separation between characterisation and classification results}}
    \label{fig:paper_orga}
\end{figure}
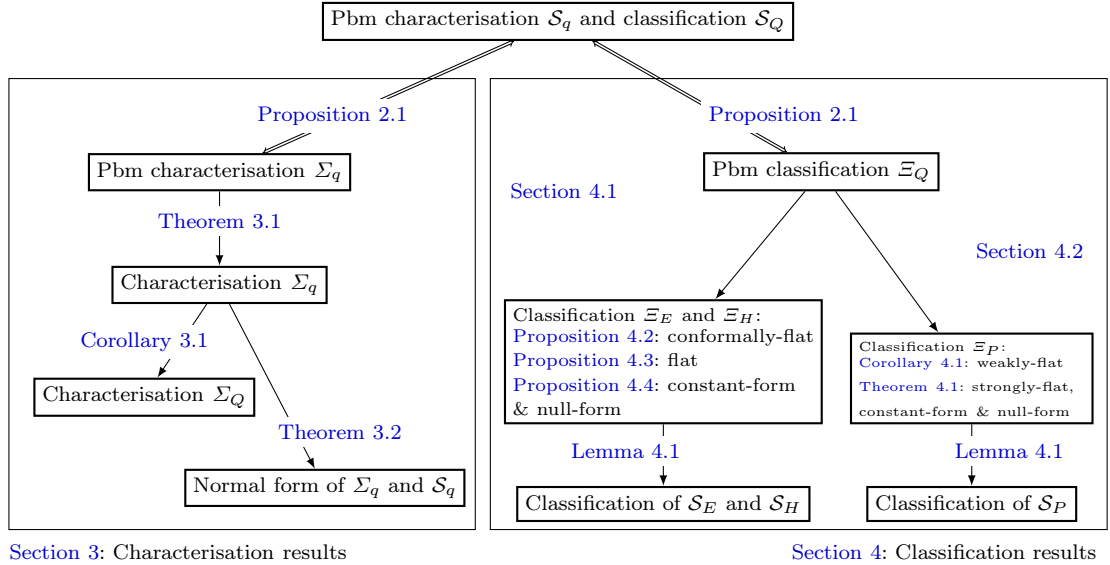

\paragraph{Organisation of the paper.} \st{The paper is organised as follows.}In the next section, we will {recall some definitions of control theory and we will} show that the problem of characterising and classifying submanifolds of $T\XXX$ can be replaced by that of characterising and classifying their first and second extensions\st{prolongations} $\Xi_{\SM}$ and $\Sigma_{\SM}$ {under feedback transformations} (see \cref{prop:equiv_of_equiv}). Moreover, we will give a first rough classification of non-degenerate conic submanifolds, introducing elliptic $\SM_E$, hyperbolic $\SM_H$, and  parabolic $\SM_P$ {sub}classes (see \cref{lem:m1_classification_conic_submanifolds}). In \cref{sec:feedback-m-1}, we will define a general second extension\st{prolongation} of a conic submanifold $\SM_q${, called a quadratic system $\Sigma_q$. {Just as we do for conic submanifolds, we} \st{We} will identify elliptic $\Sigma_E$, hyperbolic $\Sigma_H$, and parabolic $\Sigma_P$ systems {as particular cases of $\Sigma_q$},} see \cref{m1:def:quadratizable_systems} and \cref{prop:three_normal_forms_quadratic_systems}. \st{and, in}{In} \cref{thm:feedback_quadratization_m_1}, we will fully characterise the class of quadratic systems by means of a checkable relation between well-defined structure functions attached to any control-affine system. The conditions obtained in that theorem allow to give a normal form for all quadratizable control-affine systems {(i.e. systems $\Sigma$ that are feedback equivalent to $\Sigma_q$)}, see \cref{thm:m1_normal_form_quadratizable_systems}, which in turn leads to a normal form for all {non-degenerate} conic submanifolds $\SM_q$. We will also show how our characterisation and the normal form apply to the classes of elliptic, hyperbolic, and parabolic control-affine systems, which gives us a deeper insight into our conditions and in our normal form (see \cref{cor:thm_feedback_quadratization} and \cref{cor:thm:m1_normal_form_quadratizable_systems}). 
%
%
Finally, in \cref{sec:classification_conic}, we will be interested in the classification of elliptic, hyperbolic, and parabolic submanifolds {as presented in \cref{tab:classification_nomenclature} given there}. {This problem is dealt with using the classification of their first extensions (treated as control-nonlinear systems) under feedback transformations. We first show that the classification of elliptic, hyperbolic, and parabolic submanifolds presented in \cref{tab:classification_nomenclature} in \cref{sec:classification_conic} is reflected in properties of a triple of vector fields attached to their parametrisations; see \cref{lem:classification_ehp_system_to_submanifolds} statements \cref{lem:classification_ehp_system_to_submanifolds:1} to \cref{lem:classification_ehp_system_to_submanifolds:4} for the elliptic and hyperbolic cases and statements \cref{lem:classification_ehp_system_to_submanifolds:5} to \cref{lem:classification_ehp_system_to_submanifolds:8} for the parabolic case.} To every quadratic control-nonlinear system $\Xi_E$, $\Xi_H$, or $\Xi_P$ {(first extensions of elliptic, hyperbolic, and parabolic submanifolds, respectively)}, we will attach a frame of the tangent bundle (see the paragraph before \cref{prop:reparam_m1_quadratic_systems}) and we give conditions for that frame to be commutative: it turns out that in the elliptic and hyperbolic cases this requires that a certain pseudo-Riemanian metric is flat (see \cref{prop:m1:EH:existence_commutative_frame}), whereas in the parabolic case this problem can be solved without any extra assumptions (see \cref{prop:m1_exits_commutative_frame}). Then we show how we can additionally normalize the systems while preserving the commutativity of that frame. Our classifications include several normal and canonical forms, given by \cref{prop:m1:charact_class_eh_flat_systems} for elliptic and hyperbolic systems and by \cref{thm:m1_normalisation_quadratic_system_gamma_1_0} for parabolic systems. 
{We summarise the structure of the paper in \cref{fig:paper_orga}.} 
%
%
\paragraph{Related works.} A classification of quadratic control systems was initiated by Bonnard in \cite{bonnard1991Quadraticcontrolsystems}. His work differs from our as he considered homogeneous systems of degree $2$ with respect to all state variables. Hence, his class of quadratic control systems is a subclass of our parabolic systems (where we require that only one variable enters quadratically) but he considers the general dimension $n$ while our results concern $3$-dimensional systems only. In \cite{kang1992Extendedquadraticcontroller}, Krener and Kang studied the problem of equivalence, via feedback, to polynomial systems of degree $2$ modulo higher order terms. This work was continued in \cite{kang1996Extendedcontrollerform} and \cite{tall2002FeedbackClassificationNonlinear} for any degree but all those results are given for formal classification only. 
Examples of control systems subject to conic nonholonomic constraints appear in various domains of physics and engineering applications. In the next section, we will discuss Dubin's car \cite{dubins1957CurvesMinimalLength} which is a simple model of a vehicle, as well as its hyperbolic counterpart \cite{monroy-perez1998NonEuclideanDubinsProblem}. We also mention \cite{zhu2015Planartiltingmaneuver}, where the planar tilting manoeuvre problem is considered under small angle assumption, the studied control system is elliptic with respect to the states.
\section{Preliminaries}\label{sec:preliminary}
\textbf{\textit{{Main notations.}}}
\begin{table}[H]
\centering
\begin{tabularx}{\textwidth}{lX}
{$\XXX$, $T\XXX$, $x=(z,y)$} &  {A smooth 2d manifold, its tangent bundle, and its local coordinates.} \\
{$\MMM$, $\xi=(x,w)=(z,y,w)$} & {A smooth 3d manifold and its local coordinates.} \\
{$\phi$, $\phi_*$}  & {A diffeomorphism and its tangent map.}  \\
{$\SM$} & {Smooth submanifold of $T\XXX$ given by an equation of the form $\Sm(x,\dot{x})=0$, with $\rk\frac{\partial \Sm}{\partial\dot{x}}(x,\dot{x})=1$.} \\
{$\SM_q$} & {Quadratic, conic, submanifold described by $\Sm_q(x,\dot{x})=\dot{x}^t\textsf{g}(x)\dot{x}+2\omega(x)\dot{x}+h(x)=0$.}  \\
{$\SM_Q=\{\SM_E, \SM_H, \SM_P\}$} & {Set of elliptic, hyperbolic, and parabolic submanifolds.} \\
{$\Xi_{\SM}$}  & {Regular parametrisation of a submanifold $\SM$, called a first extension, seen as a control-nonlinear system.} \\
{$\Sigma_{\SM}$} & {Extension of $\Xi_{\SM}$, called a second extension of a submanifold $\SM$, seen as a control-affine system on a three dimensional manifold $\MMM$.}   \\
%
\end{tabularx}%
\end{table}%
In this section, we introduce all tools and concepts that we will need in this paper. {First, we recall some notions on the equivalence of submanifolds (of the tangent bundle) and of control systems.} {Second}, most importantly, we show that the problem of equivalence of submanifolds can be replaced by that of equivalence of their first and second {extensions}\st{prolongations}, see \cref{prop:equiv_of_equiv}. {Finally, we present quadratic submanifolds and we introduce the subclasses of elliptic, hyperbolic, and parabolic submanifolds.} In the paper the word \emph{smooth} will always mean $C^{\infty}$-smooth and throughout all systems, functions, manifolds and submanifolds are assumed to be smooth. 
\paragraph{Equivalence notions.} Recall that we have defined two submanifolds $\SM=\{\Sm(x,\dot{x})=0\}\subset T\XXX$ and $\tilde{\SM}=\{\tilde{\Sm}(\tilde{x},\dot{\tilde{x}})=0\}\subset T\tilde{\XXX}$ to be equivalent if there exists a diffeomorphism $\phi\,:\,\XXX\rightarrow\tilde{\XXX}$ and a smooth nonvanishing function $\delta\,:\,T\XXX\rightarrow\RR$ such that 
\begin{align*}
    \tilde{\Sm}(\phi(x),D\phi(x)\dot{x})&=\delta(x,\dot{x})\,S(x,\dot{x}).
\end{align*}
\noindent
If $\phi$ is defined in a neighbourhood $\XXX_0$ of $x_0$ and $\delta$ is defined in a neighbourhood of $(x_0,\dot{x}_0)\in\SM$ only, then we say that $\SM$ and $\tilde{\SM}$ are locally equivalent {at $(x_0,\dot{x}_0)$ and $(\tilde{x}_0,\dot{\tilde{x}}_0)=(\phi(x_0),D\phi(x_0)\dot{x}_0)$, respectively}.
\paragraph{\textbf{Example 1}.} The submanifolds $\SM$ and $\tilde{\SM}$ given by $\Sm(x,\dot{x})=\dot{z}-\left(-1+\sqrt{1+\dot{y}}\right)^2=0$, around $(x_0,\dot{x}_0)= (0,0)$, and $\tilde{\Sm}(\tilde{x},\dot{\tilde{x}})=\dot{\tilde{z}}-\left(\frac{\dot{\tilde{y}}}{2}\right)^2=0$, around $(\tilde{x}_0,\dot{\tilde{x}}_0)=(0,0)$, respectively, are equivalent via
\begin{align*}
    (\tilde{z},\tilde{y})=\phi(x)=(z,y-z)\quad\textrm{and}\quad \delta(x,\dot{x})=-\frac{1}{4}\left(\dot{z}-\dot{y}-2-2\sqrt{1+\dot{y}}\right).
\end{align*}
\noindent
\paragraph{\textbf{Example 2}.} Since on a $2$-dimensional manifold all metrics are locally conformally flat {(see \cite[Addendum~1 of Chapter 9]{spivak1999ComprehensiveIntroductionDifferentiala} for the Riemannian case and \cite[Theorem 7.2]{schottenloher2008MathematicalIntroductionConformal} for the Lorentzian one)}, it follows that the submanifold given by  $\Sm(x,\dot{x})=a(x)\dot{z}^2+2b(x)\dot{z}\dot{y}+c(x)\dot{y}^2-1=0$ (where $ac-b^2\neq0$) is locally equivalent to that given by $\tilde{\Sm}(\tilde{x},\dot{\tilde{x}})=\left(\dot{\tilde{z}}^2\pm\dot{\tilde{y}}^2\right)-{r(\tilde{x})^2}=0$.
\noindent
\paragraph{\textbf{Example 3}.} Since $\rk{\frac{\partial \Sm}{\partial\dot{x}}}(x,\dot{x})=1$, assume that $\frac{\partial \Sm}{\partial \dot{z}}(x_0,\dot{x}_0)\neq 0${; recall that $x=(z,y)$}. By {the implicit function theorem}, 
we can {locally} write $\Sm(x,\dot{x})=\delta(x,\dot{x})\left(\dot{z}-\sm(x,\dot{y})\right)$, where $\delta(x_0,\dot{x}_0)\neq 0$. So{, a submanifold} $\SM$ is {always locally} equivalent to {the one given by} $\dot{z}-\sm(x,\dot{y})=0$.\\

Two general control systems $\Xi\,:\,\dot{x}=F(x,w)$ and $\tilde{\Xi}\,:\,\dot{\tilde{x}}=\tilde{F}(\tilde{x},\tilde{w})$, not necessarily coming from parametrisations of submanifolds, are called \emph{feedback equivalent} if there exists a diffeomorphism $\Phi\,:\,\XXX\times\WWW\rightarrow\tilde{\XXX}\times\tilde{\WWW}$ of the form 
\begin{align*}
(\tilde{x},\tilde{w})=\Phi(x,w)=(\phi(x),\psi(x,w)),
\end{align*}
\noindent
which transforms the first system into the second, i.e.
\begin{align*}
D\phi(x)F(x,w)=\tilde{F}(\phi(x),\psi(x,w)).
\end{align*}
\noindent
The map $\phi(x)$ plays the role of a coordinates change in the state space $\XXX$, and the map $\psi$ is called feedback as it changes the parametrisation by control $w$ in a way that depends on the state $x$. If the diffeomorphism $\Phi$ is defined in a neighbourhood of $(x_0,w_0)$ only, and $\Phi(x_0,w_0)=(\tilde{x}_0,\tilde{w}_0)$, then we say that the two systems are \emph{locally feedback equivalent} at $(x_0,w_0)$ and $(\tilde{x}_0,\tilde{w}_0)$, respectively. If the diffeomorphism $\Phi$ is defined on the product $\XXX_0\times\WWW$, that is, it is global with respect to $w$, where $\XXX_0$ is a neighbourhood of a state $x_0$, and $\phi(x_0)=\tilde{x}_0$, then we say that $\Xi$ and $\tilde{\Xi}$ are locally feedback equivalent at $x_0$ and $\tilde{x}_0$, respectively. The latter local feedback equivalence will be especially useful for the class of control-affine systems
\begin{align*}
\Sigma\,:\, \dot{\xi}=f(\xi)+g(\xi)u,\quad u\in\RR,
\end{align*}
\noindent
where $f$ and $g$ are smooth vector fields on $\MMM$ (equipped with local coordinates $\xi$); below the state space $\MMM$ will be of dimension $3$ so $\MMM$ is not $\XXX$, which is of dimension $2$. For control-affine systems we will restrict the feedback transformations to the control-affine transformations 
\begin{align*}
\psi(\xi,u)=\alpha(\xi)+\beta(\xi)u,
\end{align*} 
\noindent
where $\alpha(\xi)$ and $\beta(\xi)$ are smooth functions satisfying $\beta(\cdot)\neq0$. In that case, we denote the feedback transformation by the triple $(\phi,\alpha,\beta)$ and if $\phi=\textrm{id}$, then this action is called a \emph{pure feedback} transformation and is denoted $(\alpha,\beta)$. Observe that a general control system $\Xi:\dot{x}=F(x,w)$ can be {extended}\st{prolonged} to a control-affine system  $\Xi^{{e}}$ by augmenting the state space with the control $w$ and introducing the new control $u=\dot{w}$, which gives
\begin{align*}
\Xi^{{e}}\,:\,\left\{\begin{array}{cl}
\dot{x}&=F(x,w)\\
\dot{w}&= u
\end{array}\right.,\quad u\in\RR.
\end{align*}
\noindent
Notice that $\Xi^{{e}}$ lives on the manifold $\MMM=\XXX\times \WWW$ of dimension $n=3$. \\

In introduction, we already attached to a submanifold $\SM\subset T\XXX$ a {regular} parametrisation given by a control-nonlinear system $\Xi_{\SM}\,:\,\dot{x}=F(x,w)$ satisfying $\Sm(x,F(x,w))=0$ {and $\rk\frac{\partial F}{\partial w}(x,w)=1$}. Since constructing $\Xi_{\SM}$ requires to introduce an extra variable $w$ (a control), we will also call it \st{the}{a} \emph{first {extension}\st{prolongation}} of $\SM$. Then we can attach to $\Xi_{\SM}$ its {extension}\st{prolongation} $\Xi_{\SM}^{{e}}$, denoted $\Sigma_{\SM}$ and given by 
\begin{align*}
\Sigma_{\SM}\,:\,\left\{\begin{array}{cl}
\dot{x}&= F(x,w)\\
\dot{w}&= u
\end{array}\right.,\quad u\in\RR,
\end{align*}
\noindent
and called \st{the}{a} \emph{second {extension}\st{prolongation}} of $\SM$. To distinguish different control systems attached to the submanifolds $\SM$ and $\tilde{\SM}$, we will denote $\Xi_{\tilde{\SM}}$ (resp. $\Sigma_{\tilde{\SM}}$) by $\tilde{\Xi}_{\tilde{\SM}}$ (resp. by $\tilde{\Sigma}_{\tilde{\SM}}$).
The following proposition shows that the problem of equivalence of submanifolds \st{is reflected by}{corresponds to} the equivalence {under feedback transformations} of their corresponding first and second {extensions}\st{prolongations}.
\begin{proposition}[Equivalence of equivalence relations
]\label{prop:equiv_of_equiv} Consider submanifolds $\SM=\left\{\Sm(x,\dot{x})=0\right\}$ of $T\XXX$ and $\tilde{\SM}=\left\{\tilde{\Sm}(\tilde{x},\dot{\tilde{x}})=0\right\}$ of $T\tilde{\XXX}$. The following statements are locally equivalent:
\begin{enumerate}[label=(\roman*),ref=\textit{(\roman*)}]
    \item \label{prop:equiv_of_equiv:1} The submanifolds $\SM$ and $\tilde{\SM}$ are equivalent via $\phi(x)$ and $\delta(x,\dot{x})$.
    \item \label{prop:equiv_of_equiv:2} Their regular parametrisations (first {extensions}\st{prolongations}) $\Xi_{\SM}$ and $\tilde{\Xi}_{\tilde{\SM}}$ are feedback equivalent via $\phi(x)$ and $\psi(x,w)$.
     \item\label{prop:equiv_of_equiv:3} Their second {extensions}\st{prolongations} $\Sigma_{\SM}$ and $\tilde{\Sigma}_{\tilde{\SM}}$ are feedback equivalent via $\varphi(x,w)=(\phi(x),\psi(x,w))$ and $(\alpha,\beta)$.
\end{enumerate}
That is, the following diagram is commutative:
\begin{figure}[H]
    \centering
    \begin{tikzcd}
    \mathcal{S} \arrow[d, leftrightarrow, "{(\phi,\delta)}"'] \arrow[rrr, "\textrm{parametrisation}"] &  &  & \Xi_{\mathcal{S}} \arrow[d, leftrightarrow, "{(\phi,\psi)}"'] \arrow[rr, "\textrm{{extension}}"] &  & \Sigma_{\mathcal{S}} \arrow[d, leftrightarrow, "{(\varphi,\alpha,\beta)}"] \\
    \tilde{\mathcal{S}} \arrow[rrr, "\textrm{parametrisation}"]                       &  &  & \tilde{\Xi}_{\tilde{\mathcal{S}}} \arrow[rr, "\textrm{{extension}}"]                     &  & \tilde{\Sigma}_{\tilde{\mathcal{S}}}                           
    \end{tikzcd}
\end{figure}
\end{proposition}
\begin{proof}
    \cref{prop:equiv_of_equiv:1}$\Rightarrow$\cref{prop:equiv_of_equiv:2}. {To show that the equivalence of submanifolds implies the equivalence under feedback of their first extensions, we will first use a diffeomorphism (of the state space) to analyse both first extensions in the same coordinate system. Then, we will show that the parameters (controls) of those two first extensions are related by a pure feedback transformation. To obtain invertibility of the feedback transformation, our argument will strongly rely on the regularity of the submanifold, $\rk\frac{\partial \Sm}{\partial \dot{x}}=1$, and of their first extensions, $\rk\frac{\partial F}{\partial w}=1$.} Assume that two submanifolds $\SM$ and $\tilde{\SM}$ given, respectively, by $\Sm(x,\dot{x})=0$ and $\tilde{\Sm}(\tilde{x},\dot{\tilde{x}})=0$ are equivalent via a diffeomorphism $\tilde{x}=\phi(x)$ and a nonvanishing function $\delta(x,\dot{x})$, that is $\tilde{\Sm}(\phi(x),D\phi(x)\dot{x})=\delta(x,\dot{x}),S(x,\dot{x})$. Consider $\Xi_{\SM}\,:\,\dot{x}=F(x,w)$ and $\tilde{\Xi}_{\tilde{\SM}}\,:\,\dot{\tilde{x}}=\tilde{F}(\tilde{x},\tilde{w})$, two regular parametrisations of those submanifolds. {Then, using $0=\tilde{\Sm}(\tilde{x},\dot{\tilde{x}})=\tilde{\Sm}(\tilde{x},\tilde{F}(\tilde{x},\tilde{w}))$, we have 
    \begin{align*}
        \tilde{\Sm}(\phi(x),\tilde{F}(\phi(x),\tilde{w}))=\delta(x,\left(D\phi(x)\right)^{-1}\tilde{F}(\phi(x),\tilde{w}))\Sm(x,\left(D\phi(x)\right)^{-1}\tilde{F}(\phi(x),\tilde{w})),\quad\forall\,\tilde{w}\in\tilde{\WWW},
    \end{align*}
    \noindent
    implying $\Sm(x,\hat{F}(x,\tilde{w})=0$, where $\hat{F}(x,\tilde{w})=\left(D\phi(x)\right)^{-1}\tilde{F}(\phi(x),\tilde{w})$. Therefore, $\dot{x}=F(x,w)$ and $\dot{x}=\hat{F}(x,\tilde{w})$ are two regular parametrisations of the same submanifold ${\SM}$. We will prove that ${F}({x},{w})$ and $\hat{F}({x},{\tilde{w}})$ are related by an invertible  (pure) feedback transformation of the form $\tilde{w}=\psi(x,w)$.
    }
    
    Since $\rk{\frac{\partial \Sm}{\partial \dot{x}}}(x_0,\dot{x}_0)=1$, we may assume that $\frac{\partial \Sm}{\partial \dot{z}}(x_0,\dot{x}_0)\neq0$, where $x=(z,y)$ and, by {the} implicit function theorem , we have $\Sm(x,\dot{x})=\delta(x,\dot{x})(\dot{z}-\sm(x,\dot{y}))=0$. Let $\dot{z}=F_1(x,w),\quad\dot{y}=F_2(x,w)$ be a regular parametrisation of $\SM$, it follows that $\frac{\partial F_2}{\partial w}(x_0,w_0)\neq0$. Indeed, if $\frac{\partial F_2}{\partial w}(x_0,w_0)=0$, then $\dot{z}-\sm(x,\dot{y})=0$ (recall that $\delta\neq0$), hence $F_1-\sm(x,F_2)=0$ and thus $\frac{\partial F_1}{\partial w}-\frac{\partial \sm}{\partial \dot{y}}\frac{\partial F_2}{\partial w}=0$ implying that $\frac{\partial F_1}{\partial w}(x_0,w_0)=0$; contradiction. 
    
    {Hence, the two regular parametrisations of $\SM=\{\Sm(x,\dot{x})=0\}$ given by
    \begin{align*}
        {\Xi}_{{\SM}}\,:\,\left\{\begin{array}{rl}
            \dot{{z}} &= {F}_1({x},{w}) \\
            \dot{{y}} &= {F}_2({x},{w})
        \end{array}\right.\quad\textrm{and}\quad\hat{\Xi}_{{\SM}}\,:\,\left\{\begin{array}{rl}
            \dot{{z}} &= \hat{F}_1({x},\tilde{w}) \\
            \dot{{y}} &= \hat{F}_2({x},\tilde{w})
        \end{array}\right.
    \end{align*}
    \noindent
    satisfy $\frac{\partial {F}_2}{\partial {w}}({x}_0,{w}_0)\neq0$ and $\frac{\partial \hat{F}_2}{\partial \tilde{w}}({x}_0, \tilde{w}_0)\neq0$. Therefore, $\tilde{w}=\hat{F}_2^{-1}({x},\dot{{y}})=\hat{F}_2^{-1}({x},{F}_2({x},w))${, where $\hat{F}_2^{-1}$ is the inverse with respect to the second argument}. And using $\dot{{z}}=\sm({x},\dot{{y}})$ we obtain $\dot{{z}}={F}_1({x},{w})=\sm({x},\dot{{y}})=\hat{F}_1({x},\tilde{w})$. Thus, ${\Xi}_{{\SM}}$ and $\hat{\Xi}_{{\SM}}$ are feedback equivalent via $\tilde{w}=\psi(x,w)$ with $\psi(x,w)=\hat{F}_2^{-1}(x, {F}_2(x,w))$ and the systems ${\Xi}_{{\SM}}$ and $\tilde{\Xi}_{\tilde{\SM}}$ are feedback equivalent since $\tilde{\Xi}_{\tilde{\SM}}$ \st{is the image of $\Xi_{\SM}$ under}{is the system $\hat{\Xi}_{{\SM}}$ mapped via} the diffeomorphism $\tilde{x}=\phi(x)$.\\
    }
    \indent\cref{prop:equiv_of_equiv:2}$\Rightarrow$\cref{prop:equiv_of_equiv:1}. Assume that the two regular parametrisations $\Xi_{\SM}\,:\,\dot{x}=F(x,w)$ and $\tilde{\Xi}_{\tilde{\SM}}\,:\,\dot{\tilde{x}}=\tilde{F}(\tilde{x},\tilde{w})$  of $\SM$ and $\tilde{\SM}$, respectively, are feedback equivalent via $\tilde{x}=\phi(x)$ and $\tilde{w}=\psi(x,w)$. Denote $\Phi=\phi^{-1}$, and apply the diffeomorphism ${x}=\Phi(\tilde{x})$ to $\tilde{\Xi}_{\tilde{\SM}}$ to obtain a new vector field (parametrized by $\tilde{w}$) $\hat{F}({x},\tilde{w})=D\Phi(\phi({x}))\tilde{F}(\phi({x}),\tilde{w})$ related with ${F}({x},{w})$ by the pure feedback transformation $\tilde{w}={\psi}({x},w)$. Denote ${F}=({F}_1,{F}_2)^t$ and $\hat{F}=(\hat{F}_1,\hat{F}_2)^t$ in the ${x}=({z},{y})$ coordinates. Without loss of generality we can assume that $\frac{\partial {F}_2}{\partial {w}}({x}_0,{w}_0)\neq0$ and $\frac{\partial \hat{F}_2}{\partial \tilde{w}}({x}_0, \tilde{w}_0)\neq0$.
    
    Now, apply to $\tilde{\SM}$ the same diffeomorphism ${x}=\Phi(\tilde{x})${, whose inverse we denoted by $\tilde{x}=\phi(x)$,} and \st{denote}{set} $\hat{\Sm}({x},\dot{{x}})=\tilde{\Sm}(\phi({x}), D\phi({x})\dot{{x}})$. Since by definition of $\tilde{\Xi}_{\tilde{\SM}}$ we have $\tilde{\Sm}(\tilde{x},\tilde{F}(\tilde{x},\tilde{w}))=0$, we conclude that $\hat{\Sm}({x},\hat{F}({x},\tilde{w}))=0$. Then we claim that $\frac{\partial {\Sm}}{\partial \dot{{z}}}({x}_0,\dot{{x}}_0)\neq0$; indeed, if $\frac{\partial {\Sm}}{\partial \dot{{z}}}({x}_0,\dot{{x}}_0)=0$, then $\frac{\partial }{\partial w}{\Sm}({x},{F}(x,w))=0$ yields $\frac{\partial {\Sm}}{\partial \dot{{z}}}\frac{\partial {F}_1}{\partial w}+\frac{\partial {\Sm}}{\partial \dot{{y}}}\frac{\partial {F}_2}{\partial w}=0$ and thus $\frac{\partial {\Sm}}{\partial \dot{{y}}}({x}_0,\dot{{x}}_0)=0$ giving a contradiction. The same observation holds for $\frac{\partial \hat{\Sm}}{\partial \dot{{z}}}({x}_0,\dot{{x}}_0)$. Thus by the implicit function theorem we have ${\Sm}({x},\dot{{x}}) = {\delta}({x},\dot{{x}})(\dot{{z}}-{\sm}({x},\dot{{y}}))$ and $\hat{\Sm}({x},\dot{{x}})=\hat{\delta}({x},\dot{{x}})(\dot{{z}}-\hat{\sm}({x},\dot{{y}}))$, with ${\delta}\neq0$ and $\hat{\delta}\neq0$.
    
    Using $\dot{{z}}={\sm}({x},\dot{{y}})$ and $\dot{{z}}=\hat{\sm}({x},\dot{{y}})$, we obtain ${\sm}({x},\dot{{y}})={F}_1({x},w)=\hat{F}_1({x},\tilde{w})=\hat{\sm}({x},\dot{{y}})$, hence we have 
    \begin{align*}
        \Sm({x},\dot{{x}})&= \delta(x,\dot{x})(\dot{z}-\sm(x,\dot{y}))=  \delta(x,\dot{x})(\dot{z}-\hat{\sm}({x},\dot{{y}})) = \frac{\delta}{\hat{\delta}}\tilde{\Sm}(\phi(x),D\phi(x)\dot{x})
    \end{align*}
    \noindent
    establishing the equivalence between $\SM$ and $\tilde{\SM}$.\\
    \indent\cref{prop:equiv_of_equiv:2}$\Rightarrow$\cref{prop:equiv_of_equiv:3}. If $\Xi_{\SM}$ and $\tilde{\Xi}_{\tilde{\SM}}$ are feedback equivalent , then 
    \begin{align*}
        D\phi(x)F(x,w)=\tilde{F}(\phi(x),\psi(x,w)).
    \end{align*}
    \noindent
    Thus the diffeomorphism $(\phi(x),\psi(x,w))$, of the augmented state space $(x,w)$, together with the control-affine feedback 
    \begin{align*}
        \tilde{u}=\frac{\partial \psi}{\partial x}F(x,w)+\frac{\partial \psi}{\partial w}u
    \end{align*}
    \noindent
    transform $\Sigma_{\SM}$ into $\tilde{\Sigma}_{\tilde{\SM}}$.\\
    \indent\cref{prop:equiv_of_equiv:3}$\Rightarrow$\cref{prop:equiv_of_equiv:2}. Assume that $\Sigma_{\SM}$ into $\tilde{\Sigma}_{\tilde{\SM}}$ are {feedback} equivalent via $(\tilde{x},\tilde{w})=\Phi(x,w)$ and $\tilde{u}=\alpha(x,w)+\beta(x,w)u$. Since the distribution $\distrib{\vec{w}}$ is sent by $\phi_{*}$ into $\distrib{\vec{\tilde{w}}}$, it follows that $\Phi$ has the triangular form $(\phi(x),\psi(x,w))$. Therefore, feedback equivalence of the systems $\Xi_{\SM}$ and $\tilde{\Xi}_{\tilde{\SM}}$ is established via the diffeomorphism $\tilde{x}=\phi(x)$ and the reparametrisation $\tilde{w}=\psi(x,w)$.
\qed\end{proof}
\begin{remark}
    {The use of extensions of control-nonlinear systems was introduced in \cite{vanderschaft1982ObservabilityControllabilitySmooth} and used to study controlability and observability of nonlinear systems, and then to analyse linearizability and decoupling \cite{vanderschaft1984Linearizationinputoutputdecoupling}.} Moreover, n\st{N}otice that the same proof as that of \cref{prop:equiv_of_equiv:2}$\Leftrightarrow$\cref{prop:equiv_of_equiv:3} shows that any two {control-}nonlinear systems $\Xi$ and $\tilde{\Xi}$ (which need not be regular parametrisations of submanifolds) are feedback equivalent if and only if their {extensions}\st{prolongations} $\Xi^{{e}}$ and $\tilde{\Xi}^{{e}}$ are {equivalent via } {control-affine} feedback\st{ equivalent}; see \cite{jakubczyk1990EquivalenceInvariantsNonlinear}.
\end{remark}
\begin{remark}
     The equivalence \cref{prop:equiv_of_equiv:1}$\Leftrightarrow$\cref{prop:equiv_of_equiv:2} does not hold, in general, if \st{a}{the} parametrisation {$\Xi_{\SM}:\dot{x}=F(x,w)$} does not satisfy the regularity condition $\frac{\partial F}{\partial w}(x,w)\neq0$ that we assume. To see that, consider the submanifold $\SM$ given by $\dot{z}-\dot{y}^2=0$. The parametrisations of $\SM$
    \begin{align*}
        \Xi_{\SM}\,:\,\left\{\begin{array}{rl}
            \dot{z} &= w^2  \\
            \dot{y} &= w 
        \end{array}\right.\quad\textrm{and}\quad \tilde{\Xi}_{\SM}\,:\,\left\{\begin{array}{rl}
            \dot{z} &= \tilde{w}^6  \\
            \dot{y} &= \tilde{w}^3
        \end{array}\right.
    \end{align*}
    \noindent
    are not feedback equivalent around $w_0=0$, and the reason is that $\tilde{\Xi}_{\SM}$ fails to satisfy $\frac{\partial \tilde{F}}{\partial \tilde{w}}(\tilde{w}_0)\neq0$ at $\tilde{w}_0=0$.
\end{remark}
According to the previous proposition, in order to deal with the problem of equivalence of submanifolds of the tangent bundle $T\XXX$ it is interesting, first, to study the problem of equivalence of general control-nonlinear systems $\Xi$ and, second, the problem of equivalence of general control-affine systems $\Sigma$. 
\paragraph{Conic submanifolds.} Recall that a conic, or quadratic, submanifold is given (in a{, thus in any, }\st{ suitable }local coordinate system) by $\Sm_q(x,\dot{x})=0$, where $\Sm_q$ is a smooth polynomial of degree two in $\dot{x}${, i.e.}
\begin{align*}
     \Sm_q(x,\dot{x})=\dot{x}^t\textsf{g}(x)\dot{x}+2\omega(x)\dot{x}+h(x).
\end{align*}
\noindent
The map $\Sm_q$ can be represented by the triple $\Sm_q=(\textsf{g},\omega,h)$, where $\textsf{g}$ is a symmetric $(0,2)$-tensor (possibly degenerated), $\omega$ is a one-form and $h$ is a function. Clearly, two conics $\SM_q$ of $T\XXX$ and $\tilde{\SM}_q$ of $T\tilde{\XXX}$, given by $(\textsf{g},\omega,h)$ and $(\tilde{\textsf{g}},\tilde{\omega},\tilde{h})$, respectively, are equivalent if and only if there exists a diffeomorphism $\tilde{x}=\phi(x)$ and a non-vanishing function $\delta=\delta(x)$ on $\XXX$ such that $\delta\textsf{g}=\phi^*\tilde{\textsf{g}}$, $\delta {\omega}=\phi^*\tilde{\omega}$, and $\delta h=\phi^*\tilde{{h}}$. In particular, observe that the tensors $\textsf{g}$ and $\tilde{\textsf{g}}$, which are (pseudo-)Riemanian metrics (possibly degenerated), are conformally equivalent.

It is well-known in affine geometry that such conic equations can be classified by the signature of the matrix $M_{q}(x)=\begin{pmatrix} \textsf{g}(x) & {\omega(x)^t}\\ {\omega(x)}  & h(x)\end{pmatrix}$ and that of $\textsf{g}(x)$. We will use the following two determinants
\begin{align}
    \Delta_1(x)=\det(M_q(x))
    \quad \textrm{and}\quad \Delta_2(x)= \det(\textsf{g}(x)).\label{def:determinant_conic_manifold} 
\end{align}
\noindent
Of course, $\Delta_1$ and $\Delta_2$ depend on the choice of coordinates{, however, since a diffeomorphism $\tilde{x}=\phi(x)$ transforms them as $\Delta_i=\theta^2\phi^*(\tilde{\Delta}_i)$, for $i=1,2$, where $\theta(x)=\det D\phi(x)$, their zero-level sets $\left\{\Delta_i(x)=0\right\}$ are invariantly related to the submanifold $\SM_q$.} \st{but the  ideals generated by them, in the ring of $C^{\infty}(\XXX)$-functions and thus their zero level set are invariantly related to the submanifold $\SM_q$.} In this work we will characterise non-degenerate conics, that is, non-empty and satisfying $\Delta_1\neq0${; notice that non-empty conics at points of degeneration $\Delta_1(x)=0$ form in fact (pairs of) linear subspaces of $T_x\XXX$}\st{notice that degenerate conics are in fact linear submanifolds of $T_x\XXX$}. Excluding empty $\SM_q$ is needed when considering elliptic submanifolds (see lemma below) and it implies that $M_q$ is indefinite. The non-degeneracy assumption $\Delta_1(x)\neq0$ implies $\frac{\partial \Sm_q}{\partial\dot{x}}(x,\dot{x})\neq0$ (but the converse does not hold in general) and that, locally around $x_0$, $\rk{\textsf{g}(x)}\geq1$. 

If $\Delta_2(x_0)\neq0$ or if $\Delta_2\equiv 0$ in a neighbourhood of $x_0$, then we can describe three particular types of conics given by the classification lemma below. Notice, however, that this lemma does not describe conics for which we pass smoothly from one type to another (see remark below the proof of \cref{thm:m1_normal_form_quadratizable_systems} for that case).  
\begin{lemma}[Classification of non-degenerate conics]\label{lem:m1_classification_conic_submanifolds}
Consider a nonempty conic $\SM_q$, given by $(\textsf{g},\omega,h)$, and assume $\Delta_1(x_0)\neq 0$. Then, locally around $x_0$, we have
\begin{enumerate}[label=(\roman*),ref=\textit{(\roman*)}]
    \item \label{lem:m1_classification_conic_submanifolds:1} If $\Delta_2(x_0) > 0$, then $\SM_q$ is equivalent to $\SM_E=\left\{a^2(\dot{z}-c_0)^2+b^2(\dot{y}-c_1)^2=1\right\}$, 
    \item \label{lem:m1_classification_conic_submanifolds:2} If $\Delta_2(x_0) < 0$, then $\SM_q$ is equivalent to $\SM_H=\left\{a^2(\dot{z}-c_0)^2-b^2(\dot{y}-c_1)^2=1\right\}$, 
    \item \label{lem:m1_classification_conic_submanifolds:3} If $\Delta_2 \equiv 0$, then $\SM_q$ is equivalent to $\SM_P=\left\{a\dot{y}^2-\dot{z}+b\dot{y}+c = 0\right\}$, 
\end{enumerate}
where {$a$, $b$}\st{$r$}, $c_0$, $c_1$,\st{ $a$, $b$,} and $c$ are smooth functions satisfying {$a\neq0$ and $b\neq0$ in the elliptic and hyperbolic cases, and $a\neq0$ in the parabolic case}\st{$r(\cdot)\neq0$ and $a(\cdot)\neq0$}.
\end{lemma}
We call $\SM_E$, resp. $\SM_H$, resp. $\SM_P$, an elliptic, resp. a hyperbolic, resp. a parabolic, submanifold and we will use the notation $\SM_Q$ to \st{design}{denote} the set $\{\SM_E,\SM_H,\SM_P\}$ of those three particular forms. Observe that for the parabolic form $\SM_P$, the non-degeneracy assumption $\Delta_1\neq0$ implies $\rk\textsf{g}\equiv1$ and the existence of a nonvanishing one-form $\omega=-dz+b\,dy$ satisfying $\omega\notin\ann{\ker\textsf{g}}$, whereas in the elliptic and hyperbolic cases this one-form, given by $\omega=-(c_0dz\pm c_1dy)$, can vanish at some points. Those three classes of submanifolds are related to the signature of the metric $\textsf{g}$; indeed if $\sgn{\textsf{g}}$ is constant in a neighbourhood of $x_0$, then $\SM_E$, resp. $\SM_H$, resp. $\SM_P$, corresponds to $\sgn{\textsf{g}}=(+,+)$, resp. $\sgn{\textsf{g}}=(+,-)$, resp. $\sgn{\textsf{g}}=(+,0)$; notice that we can always assume that there is at least one positive eigenvalue, otherwise take the equivalent submanifold given by $\tilde{\Sm}=-\Sm$.
\begin{proof}\leavevmode
Consider a submanifold $\SM_q$ given in local coordinates by $\Sm_q(x,\dot{x})=0$, with
\begin{align*}
    \Sm_q(x,\dot{x})=\dot{x}^t\begin{pmatrix}\textsf{g}_{11}&\textsf{g}_{12}\\\textsf{g}_{12}&\textsf{g}_{22}\end{pmatrix}\dot{x} + 2\left(\omega_1,\;\omega_2\right)\dot{x} + h,
\end{align*}
\noindent
where all functions $\textsf{g}_{ij}$, $\omega_1$, $\omega_2$, and $h$ depend smoothly on $x\in\XXX$. 
\begin{enumerate}[label=\textit{(\roman*)}]
    \item[$(i)$-$(ii)$] We deal with $\Delta_2(x_0)\neq0$, that is, the elliptic and hyperbolic cases together. In those cases, $\textsf{g}$ is a non-degenerate symmetric $(0,2)$-tensor, therefore it can be interpreted as a pseudo-Riemanian metric. Since on $2$-dimensional manifolds all metrics are {smoothly diagonalisable} ({actually, all metrics are }conformally flat{; see \cite[Addendum 1 of Chapter 9]{spivak1999ComprehensiveIntroductionDifferentiala} for the Riemannian case and \cite[Theorem 7.2]{schottenloher2008MathematicalIntroductionConformal} for the Lorentzian one)}, introduce coordinates $\tilde{x}=\phi(x)=(z,y)$ such that in those coordinates, $\SM_q$ can be written (we drop the tildes for more readability), 
    \begin{align*}
        \Sm_q=\lambda_1(x)^2\dot{z}^2\pm\lambda_2(x)^2\dot{y}^2+2\omega(x)\dot{x}+h(x)
    \end{align*}
    \noindent
    or, equivalently (since $\lambda_1(\cdot)\neq0$ and $\lambda_2(\cdot)\neq0$), as 
    \begin{align*}
        \Sm_q&= \lambda_1^2\left(\dot{z}+\frac{\omega_1}{\lambda_1^2}\right)^2\pm\lambda_2^2\left(\dot{y}\pm\frac{\omega_2}{\lambda_2^2}\right)^2 +h-\lambda_1^2\left(\frac{\omega_1}{\lambda_1^2}\right)^2\mp\lambda_2^2\left(\frac{\omega_2}{\lambda_2^2}\right)^2
    \end{align*}
    \noindent
    Notice that for this form we have $\Delta_1= \pm \lambda_1^2\lambda_2^2\left(h-\frac{(\omega_1)^2}{\lambda_1^2} \mp \frac{(\omega_2)^2}{\lambda_2^2}\right)$ which, by {our} assumption, does not vanish. Denote $c_0=-\frac{\omega_1}{\lambda_1^2}$, $c_1=\mp\frac{\omega_2}{\lambda_2^2}$, and divide by $\tilde{h}=-h+\lambda_1^2\left(\frac{\omega_1}{\lambda_1^2}\right)^2\pm \lambda_2^2\left(\frac{\omega_2}{\lambda_2^2}\right)^2$, observe that $\tilde{h}\neq0$ as $\tilde{h}=\mp\frac{1}{\lambda_1^2\lambda_2^2}\Delta_1$, to obtain 
    \begin{align*}
        \Sm_q&= \frac{\lambda_1^4\lambda_2^2}{\mp\Delta_1}\left(\dot{z}-c_0\right)^2\pm \frac{\lambda_1^2\lambda_2^4}{\mp\Delta_1}\left(\dot{y}- c_1\right)^2-1,
    \end{align*}
    \noindent
    {where the upper, resp. lower, sign corresponds to the elliptic, resp. hyperbolic, case.} In the elliptic case, if $\Delta_1>0$ the conic is empty which is excluded by assumption, therefore $\Delta_1<0$ and set {$a^2=\frac{\lambda_1^4\lambda_2^2}{-\Delta_1}$} and {$b^2=\frac{\lambda_1^2\lambda_2^4}{-\Delta_1}$} to obtain $\SM_E$. In the hyperbolic case, if $\Delta_1<0$, then permute the variables $(z,y)$ and we can thus always obtain a conic defined by $\Sm_q$ in the above form with $\Delta_1>0$. Then, set {$a^2=\frac{\lambda_1^4\lambda_2^2}{\Delta_1}$} and {$b^2=\frac{\lambda_1^2\lambda_2^4}{\Delta_1}$} to obtain $\SM_H$.
    \setcounter{enumi}{2}
    \item Assume $\Delta_2\equiv 0$. Since $\Delta_1\neq0$, we have $\rk{\textsf{g}(x)}=1$ in a neighbourhood {of $x_0$}, which implies that $\textsf{g}_{11}(x_0)\textsf{g}_{22}(x_0)\neq 0$ and thus, without loss of generality, we can assume that $\textsf{g}_{22}(x_0)\neq0$. Then, by $\rk{\textsf{g}(x)}=1$,\st{ we have $\textsf{g}_{11}=\frac{(\textsf{g}_{12})^2}{\textsf{g}_{22}}$, and} the distribution $\ker{\textsf{g}}=\distrib{\textsf{g}_{22}(x)\vec{z}-\textsf{g}_{12}(x)\vec{y}}$ is locally of constant rank and thus we introduce coordinates $\tilde{x}=(\tilde{z},\tilde{y})$ \st{satisfying}{such that}{} \st{$d\tilde{y}\perp\ker{\textsf{g}}$}$\ker{\textsf{g}}=\distrib{\vec{\tilde{z}}}$ in which we have
    \begin{align*}
       \tilde{\Sm}_q=\tilde{a}(\tilde{x})\dot{\tilde{y}}^2+2\tilde{\omega}(\tilde{x})\dot{\tilde{x}}+\tilde{h}(\tilde{x}), 
    \end{align*}
    \noindent
    whose determinant $\Delta_1=-{\tilde{a}}(\tilde{\omega}_1)^2\neq0$ implies that $\tilde{\omega}_1\neq0$. Dividing $\tilde{\Sm}_q$ by $-2\tilde{\omega}_1$ we obtain the desired form $\SM_P$ with $a=\frac{\tilde{a}}{-2\tilde{\omega}_1}$, $b=\frac{\tilde{\omega}_2}{-\tilde{\omega}_1}$, and $c=\frac{\tilde{h}}{-2\tilde{\omega}_1}$.
\end{enumerate}
\qed\end{proof}
\paragraph{\textbf{Example 4}.} There is a well known example of a {control-}nonlinear system subject to an elliptic constraint $\SM_E$, namely Dubin's car \cite{dubins1957CurvesMinimalLength}. The state of the system is\st{ $(z,y,w)$ where $(z,y){\in\RR^2}$ is} the centre of mass of the vehicle $(z,y)\in\RR^2$, and {we control}\st{$w{\in\SS^1}$ is} the orientation of the vehicle (with respect to the $z$-axis) via $w\in\SS^1$.\st{ Assume that the vehicle has {a} constant velocity and that we control the angular velocity of the orientation.} Then, the dynamics of Dubin's car reads
\begin{align*}
    \left\{\begin{array}{cl}
        \dot{z} &= r\cos(w)  \\
        \dot{y} &= r\sin(w)  
    \end{array}\right., \quad r\in\RR^*,
\end{align*}
\noindent
which clearly is a {first}\st{ second} {extension}\st{prolongation}{, regular parametrisation,} of the elliptic submanifold $\dot{z}^2+\dot{y}^2=r^2$.
\newline

We have established all notions necessary for our characterisation and classification of conic submanifolds $\SM_q$ by \st{means of }studying the feedback equivalence of their first and second {extensions}\st{prolongations} $\Xi_{\SM_q}$ and $\Sigma_{\SM_q}$.
\section{Quadratizable control-affine systems}\label{sec:feedback-m-1}
\textbf{\textit{{Notations.}}}
\begin{table}[H]
\centering
\begin{tabularx}{\textwidth}{lX}
%
%
{$\Sigma=(f,g)$} & {A control-affine system on a 3-dimensional manifold and with scalar control.}\\
{$(\phi,\alpha,\beta)$}  & {Feedback transformations acting on control-affine systems.} \\
{$\lb{g}{f}$, $\adk{g}{k}{f}$, $\dL{g}{}$}  & {Lie bracket, iterated Lie bracket, Lie derivative.}\\
{$\GGG$} & {Distribution spanned by the vector field $g$.} \\
{$\Sigma_q=(f_q,g_q)$} & {Quadratic control-affine system given, in coordinates $(x,w)$, by the vector fields $f_q=\textsf{f}_q\vec{x}$ and $g_q=\vec{w}$, which satisfy \cref{m1:def:quadratizable_systems}.} \\
{$(\rho,\tau)$} & {Structure functions attached to a control-affine system $\Sigma=(f,g)$
, see condition \cref{con:thm_nf_m1_c2} of \cref{thm:feedback_quadratization_m_1}.}\\
{$\Sigma_Q=\{\Sigma_E, \Sigma_H, \Sigma_P\}$} & {Set of elliptic, hyperbolic, parabolic subclasses of quadratic systems $\Sigma_q$.} \\
{$\Sigma_h$} & {Second extension of a submanifold $\SM=\{\dot{z}=h(x,\dot{y})\}$.}\\
%
\end{tabularx}%
\end{table}%
In this section, we {introduce the novel class of quadratic control-affine systems $\Sigma_q$ that describes second extensions of quadratic submanifolds {$\SM_q$ given by 
\begin{align*}
    \Sm_q(x,\dot{x})=\dot{x}^t\textsf{g}(x)\dot{x}+2\omega(x)\dot{x}+h(x)=0
\end{align*}
\noindent
and satisfying $\Delta_1\neq0$}.} Next, we address the equivalence problem of a control-affine system $\Sigma$ to a quadratic control-affine system $\Sigma_q${, and in that way, due to \cref{prop:equiv_of_equiv}, we provide a characterisation of quadratic submanifolds $\SM_q$}. {As a corollary, we will give a characterisation of the elliptic, hyperbolic, and parabolic submanifolds via a characterisation of corresponding subclasses of quadratic control-affine systems. Moreover, by studying our conditions, we will give a normal form of control-affine systems that are feedback equivalent to a quadratic one, and, as a consequence, this will give us a normal form of quadratic submanifolds that smoothly passes from the elliptic to the hyperbolic classes.} 

On a $3$-dimensional manifold $\MMM$, equipped with local coordinates $\xi$, we consider the control-affine system
\begin{align*}
    \Sigma\,:\, \dot{\xi}=f(\xi)+g(\xi) u,
\end{align*}
\noindent
with a scalar control $u\in\RR$ and smooth vector fields $f$ and $g$. {A control-affine system $\Sigma$ is denoted by the pair $\Sigma=(f,g)$, and we set $\GGG=\distrib{g}$ the distribution spanned by the vector field $g$. Moreover, we}\st{For this system, we denote $\GGG=\distrib{g}$, the distribution spanned by $g$, and we} will use the following notations: given two vector fields $g$ and $f$ on $\MMM$, by $\lb{g}{f}$ we denote the Lie bracket of $g$ and $f$, in coordinates we have $\lb{g}{f}=\frac{\partial f}{\partial \xi}g-\frac{\partial g}{\partial \xi}f$, and $\adk{g}{k}{f}=\lb{g}{\adk{g}{k-1}{f}}$ stands for the iterated Lie bracket, with {the convention} $\adk{g}{0}{f}=f$\st{, and $\phi_*$ denotes the tangent map of a diffeomorphism $\phi$}.
\begin{definition}[Quadratic and quadratizable systems]\label{m1:def:quadratizable_systems}
We say that a control-affine system $\Sigma=(f,g)$ is \emph{quadratizable} if it is feedback equivalent to
{\begin{align*}
    \Sigma_q\,:\,\left\{\begin{array}{cl}
        \dot{x} &= \textsf{f}_q(x,w) \\
        \dot{w} &= u
    \end{array} \right. ,
\end{align*}
where $x\in\XXX$, a $2$-dimensional manifold, $w\in\WWW\subset\RR$, and $\textsf{f}_q(x,w)$ satisfies $\left(\frac{\partial^2 \textsf{f}_q}{\partial w^2}\wedge\frac{\partial \textsf{f}_q}{\partial w}\right)(x,w)\neq0$
, and 
\begin{align*}
\frac{\partial^3 \textsf{f}_q}{\partial w^3}=\tau_q(x)\, \frac{\partial \textsf{f}_q}{\partial w},
\end{align*}
\noindent
with $\tau_q$ a smooth function of the indicated variable $x\in\XXX$. A system $\Sigma_q$ of the above form is called quadratic.}
\end{definition}
{For a quadratic system $\Sigma_q=(f_q,g_q)$, where $f_q=\textsf{f}_q\vec{x}$ and $g_q=\vec{w}$,}\st{Set$x=(z,y)$ and denote $\textsf{f}_q=f^1\vec{z}+f^2\vec{y}$;} it is shown in \cref{appendix:m1_solution_equation_quadratizable} that any smooth $\textsf{f}_q$ can be written, locally around $(x_0,w_0)$, as
\begin{align}\label{eq:inifite_series_form_sigma_Q}
    \textsf{f}_q(x,w) &= A(x)\sum_{k=0}^{+\infty}\frac{(w-w_0)^{2k+2}}{(2k+2)!}\tau_q(x)^k+B(x)\sum_{k=0}^{+\infty}\frac{(w-w_0)^{2k+1}}{(2k+1)!}\tau_q(x)^k+C(x),
\end{align}
\noindent
where $A,B,C$ are smooth. {Since $\GGG=\distrib{g_q}=\distrib{\vec{w}}$ is involutive and of constant rank one, it is integrable and its integral curves can be identified with the points $x$ of $\XXX$ and $A$, $B$, and $C$ can be seen as smooth vector fields on $\XXX$, for which we have $A\wedge B\neq0$.}\st{where $A$, $B$, and $C$ can be seen as smooth vector fields on  $\XXX\cong\MMM/\GGG$ for which we have $A\wedge B\neq0$.} The following proposition shows that $\Sigma_q$ is a second {extension}\st{prolongation} of a {quadratic}\st{conic} submanifold {$\SM_q$}, thus justifies to call $\Sigma_q$ a quadratic system, and {identifies the elliptic, hyperbolic, and parabolic subclasses by describing}\st{describes} three normal forms of $f_q=\textsf{f}_q\vec{x}$ given for $\tau_q \neq0$ and $\tau_q\equiv0$.
\begin{proposition}\leavevmode\label{prop:three_normal_forms_quadratic_systems}Locally around $\xi_0\in\MMM$ we have the followings
    \begin{enumerate}[label=(\roman*),ref=\textit{(\roman*)}]
        \item \label{prop:three_normal_forms_quadratic_systems:1} $\Sigma_q$ is a second {extension}\st{prolongation} of a conic submanifold $\SM_q$ {and, conversely, any second extension $\Sigma_{\SM_q}$ of a conic submanifold $\SM_q$ is feedback equivalent to a system of the form $\Sigma_q$}.
        \item \label{prop:three_normal_forms_quadratic_systems:2} If $\tau_q(\xi_0) <0$ in a neighbourhood, resp. $\tau_q(\xi_0) >0$, resp. $\tau_q \equiv 0$, then $\Sigma_q$ is locally feedback equivalent to $\Sigma_E$, resp. $\Sigma_H$, resp. $\Sigma_P$, given by, respectively, $\textsf{f}_q$ of the form
        \begin{multline*}
            \textsf{f}_E = {A}(x)\cos(\tilde{w})+{B}(x)\sin(\tilde{w})+{C}(x),\qquad
            \textsf{f}_H = {A}(x)\cosh(\tilde{w})+{B}(x)\sinh(\tilde{w})+{C}(x),\\
            \textsf{f}_P = A(x)w^2+B(x)w+C(x),
        \end{multline*}
        \noindent
        {where in all three cases $A\wedge B\neq0$.}
        \item \label{prop:three_normal_forms_quadratic_systems:3} $\Sigma_E$, resp. $\Sigma_H$, resp. $\Sigma_P$, is a second {extension}\st{prolongation} of a conic submanifold $\SM_q$ satisfying $\Delta_2>0$, resp. $\Delta_2<0$, resp. $\Delta_2\equiv0$.
    \end{enumerate}
\end{proposition}
The conic submanifold $\SM_q$ of item \cref{prop:three_normal_forms_quadratic_systems:3} is, by \cref{lem:m1_classification_conic_submanifolds}, equivalent to $\SM_E$ (if $\Delta_2>0$), resp. $\SM_H$ (if $\Delta_2<0$), resp. $\SM_P$  (if $\Delta_2\equiv0$). So it is natural to call $\Sigma_E$ an elliptic system, $\Sigma_H$ a hyperbolic system, and $\Sigma_P$ a parabolic system. We will denote by $Q$ the set $\{E,H,P\}$, and, consequently, $f_Q=\{f_E,f_H,f_P\}$, {$g_Q=\{g_E,g_H,g_P\}$,  $\textsf{f}_Q=\{\textsf{f}_E,\textsf{f}_H,\textsf{f}_P\}$,} and $\Sigma_Q=\{\Sigma_E,\Sigma_H,\Sigma_P\}$. 
\begin{proof}\leavevmode
    \begin{enumerate}[label=\textit{(\roman*)},ref=\textit{(\roman*)}]
        \item Consider $\textsf{f}_q$ given by \cref{eq:inifite_series_form_sigma_Q}, for simplicity {of the notations} we assume $w_0=0$. By \cref{appendix:m1:rectification_2_distrib_rk_1}, we choose coordinates $(z,y)=\phi(x)$ that rectify simultaneously the distributions spanned by $A$ and $B$ so we may assume that $A=a\vec{z}$ and $B=b\vec{y}$, where $a$ and $b$ are smooth functions satisfying $a(x_0)b(x_0)\neq0$. {We set $C=c_0\vec{z}+c_1\vec{y}$, where $c_0$ and $c_1$ are smooth functions of $x$.} Using Cauchy products we compute
        \begin{align*}
            \left(\frac{\dot{z}-c_0}{a}\right)^2 &= \left(\sum_{k=0}^{+\infty}\frac{w^{2k+2}}{(2k+2)!}\tau_q^k\right)^2=\sum_{k=0}^{+\infty}\frac{(8\cdot4^k-2)w^{2k+4}\,\tau_q^k}{(2k+4)!}, \\ 
            \left(\frac{\dot{y}-c_1}{b}\right)^2 &= \left(\sum_{k=0}^{+\infty}\frac{w^{2k+1}}{(2k+1)!}\tau_q^k\right)^2=\sum_{k=0}^{+\infty}\frac{2\cdot4^kw^{2k+2}\,\tau_q^k}{(2k+2)!}, \\ 
            &= w^2+\sum_{k=1}^{+\infty}\frac{2\cdot4^kw^{2k+2}\,\tau_q^k}{(2k+2)!} 
            = w^2+\sum_{k=0}^{+\infty}\frac{2\cdot4^{k+1}w^{2k+4}\,\tau_q^{k+1}}{(2k+4)!},\\ 
            &=w^2+\tau_q \sum_{k=0}^{+\infty}\frac{8\cdot4^{k}w^{2k+4}\,\tau_q^{k}}{(2k+4)!} 
            = w^2+\tau_q \left(\frac{\dot{z}-c_0}{a}\right)^2 + \tau_q \sum_{k=0}^{+\infty}\frac{2 w^{2k+4}\,\tau_q^k}{(2k+4)!}, \\
            \left(\frac{\dot{y}-c_1}{b}\right)^2-\tau_q \left(\frac{\dot{z}-c_0}{a}\right)^2&= w^2+\sum_{k=0}^{+\infty}\frac{2 w^{2k+4}\,\tau_q^{k+1}}{(2k+4)!}=\sum_{k=-1}^{+\infty}\frac{2 w^{2k+4}\,\tau_q^{k+1}}{(2k+4)!},\\
            &= \sum_{k=0}^{+\infty}\frac{2 w^{2k+2}\,\tau_q^k}{(2k+2)!} = 2\left(\frac{\dot{z}-c_0}{a}\right).
        \end{align*}
        \noindent
        Hence, \st{the} $\Sigma_q$ is a second {extension}\st{prolongation} of the submanifold $\SM_q$ given by
        \begin{align}\label{prf:normal_sm_q}
            \left(\frac{\dot{y}-c_1}{b}\right)^2-\tau_q \left(\frac{\dot{z}-c_0}{a}\right)^2-2\left(\frac{\dot{z}-c_0}{a}\right)=0,
        \end{align}
        \noindent
        and for which we have $\Delta_1=-\frac{1}{a^2b^2}$ and $\Delta_2=-\frac{\tau_q}{a^2b^2}$. 
        
        {To prove that any second extension of a conic submanifold $\SM_q$ is feedback equivalent to a system of the form $\Sigma_q$, we show that any conic submanifold is equivalent via a diffeomorphism to a conic submanifold of the form given by \cref{prf:normal_sm_q}. To this end, consider $\Sm_q=\dot{x}^t\textsf{g}\dot{x}+2\omega\dot{x}+h=a\dot{z}^2+2b\dot{z}\dot{y}+c\dot{y}^2+2d\dot{z}+2e\dot{y}+h=0$, where all functions $a,b,c,d,e,h$ depend smoothly on $x=(z,y)$. By $\Delta_1(x_0)\neq0$, we have $\rk \textsf{g}(x_0)\geq1$, where $\textsf{g}=\left(\begin{smallmatrix}a&b\\b&c\end{smallmatrix}\right)$. If $\rk\textsf{g}(x_0)=2$, then we can apply the results of \cref{lem:m1_classification_conic_submanifolds} to get $\SM_q$ in the form of \cref{prf:normal_sm_q}. If $\rk\textsf{g}(x_0)=1$, then we may assume that $c(x_0)\neq0$. Indeed, if $a(x_0)=c(x_0)=0$, then $b(x_0)\neq0$ implying that $\rk\textsf{g}(x_0)=2$ so either $a(x_0)\neq0$ or $c(x_0)\neq0$ and we can always suppose $c(x_0)\neq0$ by permuting $y$ and $z$, if necessary. Dividing $\Sm_q=0$ by $c$ we obtain (we keep the same names for all remaining functions): 
        \begin{align*}
            \Sm_q&=a\dot{z}^2+2b\dot{z}\dot{y}+\dot{y}^2+2d\dot{z}+2e\dot{y}+h=(a-b^2)\dot{z}^2+(b\dot{z}+\dot{y})^2+2d\dot{z}+2e\dot{y}+h=0.
        \end{align*}
        \noindent
        Choose local coordinates $(\tilde{z},\tilde{y})=(z,\psi(z,y))$, hence $\frac{\partial \psi}{\partial y}\neq0$, that rectify the line-distribution $\distrib{\vec{z}-b\vec{y}}$, i.e. $\frac{\partial\psi}{\partial z}-b\frac{\partial\psi}{\partial y}=0$. Therefore, $\dot{\tilde{z}}=\dot{z}$ and $\dot{\tilde{y}}=\frac{\partial \psi}{\partial z}\dot{z}+\frac{\partial \psi}{\partial y}\dot{y}=\frac{\partial \psi}{\partial y}\left(b\dot{z}+\dot{y}\right)$. Thus $\tilde{\SM}_q$, which is $\SM_q$ in $(\tilde{z},\tilde{y})$-coordinates, reads 
        \begin{align*}
            \tilde{\Sm}_q=\frac{1}{\tilde{b}^2}\dot{\tilde{y}}^2+\tilde{a}\dot{\tilde{z}}^2+2\tilde{d}\dot{\tilde{z}}+2\tilde{e}\dot{\tilde{y}}+\tilde{h}=0,
        \end{align*}
        \noindent
        where $\tilde{b}=\frac{\partial \psi}{\partial y}$ and $\tilde{a}(\tilde{x}_0)=0$ because $\rk\tilde{\textsf{g}}(\tilde{x}_0)=1$. By $\tilde{\Delta}_1(\tilde{x}_0)\neq0$ and $\tilde{a}(\tilde{x}_0)=0$, we conclude that $\tilde{d}(\tilde{x}_0)\neq0$ and thus dividing by $\tilde{d}$ (and removing the "tildes" from the coordinates) we get 
        \begin{align*}
            \bar{\Sm}_q=\frac{\tilde{\Sm}_q}{\tilde{d}}=\frac{1}{\bar{b}^2}\dot{{y}}^2+\bar{A}\dot{{z}}^2+2\dot{{z}}+2\bar{e}\dot{{y}}+\bar{h}=\left(\frac{\dot{{y}}-\bar{c}_1}{\bar{b}}\right)^2+\bar{A}\dot{{z}}^2+2\dot{{z}}+\bar{H}=0,
        \end{align*}
        \noindent
        with $\bar{b}=\tilde{b}\sqrt{\tilde{d}}$ (if $\tilde{d}<0$, then we take ${z}=-\tilde{z}$ to get $\tilde{d}>0$), $\bar{A}=\frac{\bar{a}}{\tilde{d}}$ (satisfying $\bar{A}(x_0)=0$), $\bar{c}_1=-\bar{e}\bar{b}^2$, and $\bar{H}=\bar{h}-\bar{e}^2\bar{b}^2$ (see also \cite{427307} for another proof for the smooth diagonalisation of a symmetric $(0,2)$-tensor in the case of rank deficiency). To get the form \cref{prf:normal_sm_q}, we will now prove that there exists, locally around ${x}_0$, functions $\bar{\tau}_q$, $\bar{a}\neq0$, and $\bar{c}_0$ such that $\bar{A}\dot{{z}}^2+2\dot{{z}}+\bar{H}=-\bar{\tau}_q\left(\frac{\dot{{z}}-\bar{c}_0}{\bar{a}}\right)^2-2\left(\frac{\dot{{z}}-\bar{c}_0}{\bar{a}}\right)$. We obtain 
        \begin{align*}
            \frac{-\bar{\tau}_q}{\bar{a}^2}=\bar{A},\quad \frac{\bar{\tau}_q\bar{c}_0}{\bar{a}^2}-\frac{1}{\bar{a}}=1,\quad\textrm{and}\quad-\frac{\bar{\tau}_q(\bar{c}_0)^2}{\bar{a}^2}+2\frac{\bar{c}_0}{\bar{a}}=\bar{H}.
        \end{align*}
        \noindent
        Hence, from the second equation, $\frac{\bar{c}_0}{\bar{a}}=-\bar{A}(\bar{c}_0)^2-\bar{c}_0$ implying that $\bar{A}(\bar{c}_0)^2+2\bar{c}_0+\bar{H}=0$. The latter equation possesses smooth solutions $\bar{c}_0(z,y)$ because the discriminant (evaluated at $x_0$) of that degree two polynomial is $\bar{\Delta}(x_0)=4>0$ (recall that $\bar{A}(x_0)=0$). So 
        \begin{align*}
            \bar{\Sm}_q=\left(\frac{\dot{{y}}-\bar{c}_1}{\bar{b}}\right)^2-\bar{\tau}_q\left(\frac{\dot{{z}}-\bar{c}_0}{\bar{a}}\right)^2-2\left(\frac{\dot{{z}}-\bar{c}_0}{\bar{a}}\right)=0,
        \end{align*}
        \noindent
        where $\bar{a}=\frac{-1}{1+\bar{A}\bar{c}_0}$, which is well defined since $\bar{A}(x_0)=0$, and $\bar{\tau}_q=-\bar{A}\bar{a}^2$, proving that $\bar{\Sm}_q$ is of the desired form \cref{prf:normal_sm_q}.
        }
        \item If $\tau_q < 0$, then $\frac{\partial^3\textsf{f}_q}{\partial w^3}=\tau_q \frac{\partial \textsf{f}_q}{\partial w}$ implies that $\textsf{f}_q=A\cos(\sqrt{-\tau_q}w)+B\sin(\sqrt{-\tau_q}w)+C${, where $A$, $B$, and $C$ depend on $x$ and, clearly, $\frac{\partial^2\textsf{f}_q}{\partial w^2}\wedge\frac{\partial \textsf{f}_q}{\partial w}\neq0$ implies $A\wedge B\neq0$ (for this case, as well as for the next two cases).} Via the change of coordinate $\tilde{w}=\sqrt{-\tau_q}w$, we get $\textsf{f}_q=\textsf{f}_E$. If $\tau_q >0$, then $\frac{\partial^3\textsf{f}_q}{\partial w^3}=\tau_q \frac{\partial \textsf{f}_q}{\partial w}$ implies that $\textsf{f}_q=A\cosh(\sqrt{\tau_q}w)+B\sinh(\sqrt{\tau_q}w)+C$ which, via the change of coordinate $\tilde{w}=\sqrt{\tau_q}w$, gives $\textsf{f}_q=\textsf{f}_H$. Finally, if $\tau_q \equiv0$, then $\frac{\partial^3\textsf{f}_q}{\partial w^3}\equiv 0$ implies that $\textsf{f}_q=\textsf{f}_P=Aw^2+Bw+C$.
        \item {By item \cref{prop:three_normal_forms_quadratic_systems:1}, $\Sigma_q$ is a second extension of a conic submanifold $\Sm_q$ for which we have $\Delta_2=-\frac{\tau_q}{a^2b^2}$. Thus, $\Delta_2>0$, resp. $\Delta_2<0$, resp. $\Delta_2\equiv0$ if and only if $\tau_q< 0$, resp. $\tau_q>0$, resp. $\tau_q\equiv0$, which, by item \cref{prop:three_normal_forms_quadratic_systems:2}, correspond to $\Sigma_E$, resp. $\Sigma_H$, resp. $\Sigma_P$.} 
    \end{enumerate}
\qed\end{proof}

Notice that $\tau_q$ plays for $\Sigma_q$ \st{a similar role as the one that $\Delta_2$ plays for $\SM_q$ and, indeed, we see in the proof of the first item that in a suitable coordinate system we have $\Delta_2=-\tau_q$.}an analogous role to that played by $\Delta_2$ for $\SM_q$; indeed, we have $\sgn{\Delta_2(x)}=-\sgn{\tau_q(x)}$. In particular, the sign of $\tau_q$ identifies the subclasses of elliptic, hyperbolic, and parabolic control-affine systems. {Moreover, statement \cref{prop:three_normal_forms_quadratic_systems:1} shows that every second extension of a quadratic submanifold $\SM_q$ is feedback equivalent to a quadratic system $\Sigma_q$ and the other way around, therefore to obtain a characterisation of quadratic submanifolds, it is crucial to characterise the class of quadratizable systems.}\\

The remaining part of this section is organised as follows. First, we will state our main theorem giving necessary and sufficient conditions characterising the class of quadratic control-affine systems $\Sigma_q$. Second, by carefully studying the conditions of that theorem, we will give a normal form of the quadratizable systems. 

\subsection{Characterisation of quadratizable control-affine systems} \label{sub-sec:quadratization_of_submanifolds}
We now focus on the {feedback} equivalence of a general control-affine system $\Sigma\,:\,\dot{\xi}=f(\xi)+g(\xi)u$ with {a} quadratic control-affine system {of the form} $\Sigma_q$. The theorem below gives checkable necessary and sufficient conditions in terms of the vector fields $f$ and $g$ of $\Sigma$ for the existence of a smooth feedback transformation $(\phi,\alpha,\beta)$ that locally brings $\Sigma$ into a quadratic system $\Sigma_q$. Equivalence to particular subcases of $\Sigma_q$, namely elliptic $\Sigma_E$, hyperbolic $\Sigma_H$, and parabolic $\Sigma_P$, is provided by \cref{cor:thm_feedback_quadratization} below. 

\begin{theorem}[Feedback quadratisation]\label{thm:feedback_quadratization_m_1}
    Let $\Sigma=(f,g)$ be a control-affine system on a 3-dimensional smooth manifold with a scalar control. {The system }$\Sigma$ is, locally around $\xi_0\in\MMM$, feedback equivalent to a quadratic system $\Sigma_q$ if and only if
    \begin{enumerate}[label=\textrm{(C\arabic*)}, ref=\text{(C\arabic*)}, font=\normalfont]
        \item \label{con:thm_nf_m1_c1} $g\wedge\ad{g}{f}\wedge\adk{g}{2}{f}\; (\xi_0) \neq 0$, 
        \item \label{con:thm_nf_m1_c2} The structure functions $\rho$ and $\tau$ in the decomposition $\adk{g}{3}{f} = \structfunctA\,\adk{g}{2}{f} + \structfunctB\,\ad{g}{f}\mod \GGG$ satisfy, locally around $\xi_0$,
        \begin{align}\label{eq:thm_nf_m1_c2}
            \dL{g}{\chi}-\frac{2}{3}\structfunctA\chi=0,
        \end{align}
        \noindent
        where $\chi = 3\dL{g}{\structfunctA}-2\structfunctA^2-9\structfunctB$.
    \end{enumerate}
\end{theorem}
\noindent
Condition \cref{con:thm_nf_m1_c1} is a regularity condition, it ensures that the vector fields $g,\ad{g}{f},$ and $\adk{g}{2}{f}$ are locally linearly independent and thus that they form a local frame, hence the structure functions $(\structfunctA,\structfunctB)$ of \cref{con:thm_nf_m1_c2} are well defined. The main idea behind this theorem is to observe that if for $\Sigma$ we have  $\adk{g}{3}{f}=\tau(x)\,\ad{g}{f}${, modulo $\GGG=\distrib{g}$}, i.e. the third Lie derivative of $f$ along $g$ is proportional to the first Lie derivative of $f$ along $g${, modulo $\GGG$}, then with the help of a diffeomorphism we can obtain the form $\Sigma_q$, see the sufficiency part of the proof for details. Thus condition \cref{con:thm_nf_m1_c2} shows how that relation changes when we allow for feedback transformations $(\alpha,\beta)$. 

\begin{proof}\emph{Necessity.} 
    Consider the affine control system $\Sigma$ given by two smooth vector fields $f$ and $g$ and recall that $\GGG$ is the distribution $\GGG=\distrib{g}$. Let $(\alpha,\beta)$ form a control-affine feedback and let $\phi$ be a diffeomorphism such that $\Sigma$ is, locally, \st{equivalent to}{transformed into} $\Sigma_q$ {via $\phi$ and $(\alpha,\beta)$}. In coordinates $\tilde{\xi}=\phi(\xi)=(\tilde{x},\tilde{w})=(\tilde{z},\tilde{y},\tilde{w})$, we denote $\tilde{f}_q{=\tilde{\textsf{f}}_q\vec{\tilde{x}}}$ and $\tilde{g}_q{=\vec{\tilde{w}}}$ the vector fields of $\Sigma_q$, $\tilde{\GGG}$ the distribution spanned by $\tilde{g}_q$, and $(\tilde{\structfunctA},\tilde{\structfunctB})$ the structure functions of $\Sigma_q${, defined as in \cref{con:thm_nf_m1_c2}}. By definition of feedback equivalence the following relations between $(f,g)$ and $(\tilde{f}_q,\tilde{g}_q)$ hold: $\tilde{f}_q=\phi_*\left(f+\alpha g\right)$ and $\tilde{g}_q=\phi_*\left(g\beta\right)$.\\ 
    The system $\Sigma_q$ is quadratic, so by \cref{m1:def:quadratizable_systems}, we have $\frac{\partial^2 \tilde{\textsf{f}}_q}{\partial \tilde{w}^2}\wedge\frac{\partial \tilde{\textsf{f}}_q}{\partial \tilde{w}}(\tilde{x}_0,\tilde{w}_0)\neq0$, which implies that \cref{con:thm_nf_m1_c1} holds for $\Sigma_q$, and we also have $\frac{\partial^3\tilde{\textsf{f}}_q}{\partial \tilde{w}^3}=\tilde{\structfunctB}_q\frac{\partial \tilde{\textsf{f}}_q}{\partial \tilde{w}}$, thus we get $\tilde{\structfunctA}=0$ and $\tilde{\structfunctB}=\tilde{\structfunctB}_q(\tilde{z},\tilde{y})$. Therefore for $\Sigma_q$ we have $\tilde{\chi}{(\tilde{z},\tilde{y})}=-9\tilde{\structfunctB}(\tilde{z},\tilde{y})$ implying $\dL{\tilde{g}_q}{\tilde{\chi}}-\frac{2}{3}\tilde{\rho}\tilde{\chi}=\frac{\partial \tilde{\chi}}{\partial \tilde{w}}=0$. Hence $\Sigma_q$ satisfies \cref{con:thm_nf_m1_c1} and \cref{con:thm_nf_m1_c2} and we will now prove that those conditions are invariant under diffeomorphisms $\phi$ and feedback transformations $(\alpha,\beta)$.\\ 
    Clearly, \cref{con:thm_nf_m1_c1} is invariant under diffeomorphisms (as $\lb{\phi_{*}{f}}{\phi_{*}{g}}=\phi_{*}\lb{{f}}{{g}}$) and under feedback $(\alpha,\beta)$ since $\beta\neq 0$. We have checked that \cref{con:thm_nf_m1_c2} holds for $\Sigma_q=(\tilde{f}_q,\tilde{g}_q)$ an, clearly, \cref{con:thm_nf_m1_c2} is invariant under diffeomorphisms since they conjugate structure functions. Moreover \cref{con:thm_nf_m1_c2} is invariant under the transformation $\tilde{f}_q\mapsto \tilde{f}_q+\alpha \tilde{g}_q$, since the expression of $\adk{\tilde{g}_q}{3}{\tilde{f}_q}$ is considered modulo the distribution $\tilde{\GGG}$. Finally, under the action of $\beta$ the brackets, with $\tilde{g}_q=g\beta$, are transformed by
    \begin{align*}
        \ad{\tilde{g}_q}{f_q} &= \beta\ad{g}{f} \mod \tilde{\GGG}, \\ 
        \adk{\tilde{g}_q}{2}{f_q} &= \beta^2\adk{g}{2}{f} + \beta\dL{g}{\beta}\ad{g}{f} \mod \tilde{\GGG}, \\ 
        \adk{\tilde{g}_q}{3}{f_q} &= (\beta^3\structfunctA + 3\beta^2\dL{g}{\beta})\adk{g}{2}{f} + (\beta^3\structfunctB + \beta\dL{g}{\beta\dL{g}{\beta}}) \ad{g}{f}  \mod \tilde{\GGG}, \\
        &= (\structfunctA \beta+3\dL{g}{\beta})\adk{\tilde{g}_q}{2}{f} + \left(\structfunctB\beta^2 + \dL{g}{\beta\dL{g}{\beta}}-\structfunctA \beta\dL{g}{\beta}-3\left(\dL{g}{\beta}\right)^2 \right)\ad{\tilde{g}_q}{f} \mod \tilde{\GGG}
    \end{align*}
    \noindent
    This implies that the structure functions $\tilde{\structfunctA}$ and $\tilde{\structfunctB}$ of $\Sigma_q$ defined by $\adk{\tilde{g}_q}{3}{f_q} = \tilde{\structfunctA}\,\adk{\tilde{g}_q}{2}{f_q} + \tilde{\structfunctB}\ad{\tilde{g}_q}{f_q}\mod \tilde{\GGG}$ are given in terms of the feedback transformation $\beta$ and the structure functions $\structfunctA$ and $\structfunctB$ of $\Sigma$ by
    \begin{align}\label{eq:structure-functions-m1-transformation}
        \tilde{\structfunctA }= \structfunctA \beta+3\dL{g}{\beta},\quad \tilde{\structfunctB}=\structfunctB\beta^2+ \dL{g}{\beta\dL{g}{\beta}}-\structfunctA \beta\dL{g}{\beta}-3\left(\dL{g}{\beta}\right)^2.
    \end{align}
    \noindent
    Since for $\Sigma_q$ the structure function $\tilde{\structfunctA}=0$\st{vanishes}, \st{$\tilde{\structfunctA }=0$, }we have the relation $\dL{g}{\beta}=-\frac{\beta\structfunctA}{3}$ and thus $\tilde{\chi}=-9\tilde{\structfunctB}$, which is equal to 
    \begin{align*}
        \tilde{\chi}&=-9\left(\structfunctB\beta^2+\dL{g}{\beta\dL{g}{\beta}}\right)=-9\left(\structfunctB\beta^2+\dL{g}{\frac{-\beta^2\structfunctA}{3}}\right),\\ 
        &=-9\left(\structfunctB\beta^2-\frac{1}{3}\left(\structfunctA\dL{g}{\beta^2}+\beta^2\dL{g}{\structfunctA}\right)\right)= -9\beta^2\left(\tau-\frac{1}{3}\dL{g}{\structfunctA}+\frac{2}{9}\structfunctA^2\right)= \beta^2\chi.
    \end{align*}
    \noindent
    And finally, 
    \begin{align*}
        \dL{\tilde{g}_q}{\tilde{\chi}}&=\beta\dL{g}{\beta^2\chi}=\beta^3\dL{g}{\chi}+2\beta^2\chi\dL{g}{\beta}= \beta^3\dL{g}{\chi}-\frac{2}{3}\beta^3\chi\structfunctA =0,
    \end{align*}
    showing the necessity of relation \cref{eq:thm_nf_m1_c2} and concludes the necessity part of the proof. \\ 
    \paragraph{Sufficiency.} There are two steps in the sufficiency part. The first one consists in building a vector field $g$ such that $\adk{g}{3}{f}=\structfunctB\,\ad{g}{f}\mod\GGG$ with $\dL{g}{\tau}=0$. Then we will construct a diffeomorphism $\phi$ that brings $\Sigma$ into the form $\Sigma_q$. \\ 
    Consider the system $\Sigma\,:\,\dot{\xi}=f+g u$, for which we assume $g\wedge\ad{g}{f}\wedge\adk{g}{2}{f}\,(\xi_0)\neq 0$ and suppose that relation \cref{eq:thm_nf_m1_c2} holds for the structure functions $\structfunctA$ and $\structfunctB$ of $\Sigma$. Choose a function $\beta\neq0$ satisfying $\dL{g}{\beta}=\frac{-\beta \structfunctA}{3}${,which exists locally since $g\neq0$ by condition \cref{con:thm_nf_m1_c1}; to guarantee that $\beta\neq 0$, we actually may solve the equation $\dL{g}{\ln(\beta)}=-\frac{\structfunctA}{3}$} \st{(to guarantee that $\beta\neq 0$, we actually may solve the equation $\dL{g}{\ln(\beta)}=-\frac{\structfunctA}{3}$)}. Define the system $\tilde{\Sigma}\,:\, \dot{\xi}=\tilde{f}+\tilde{g}\tilde{u}$, where $\tilde{g}=g\beta$ and $\tilde{f}=f$, then by \cref{eq:structure-functions-m1-transformation} the structure function $\tilde{\structfunctA}$ of $\tilde{\Sigma}$ vanishes. Therefore, we have $\tilde{\chi}=-9\tilde{\structfunctB}$ and thus relation \cref{eq:thm_nf_m1_c2} implies that $\dL{\tilde{g}}{\tilde{\structfunctB}}=0$.
    
    Since $\tilde{g}\neq0$, we apply a diffeomorphism $(z,y,w)=\phi(\xi)$ such that $\phi_*\tilde{g}=g_q=\vec{w}$ and denote $f_q=\phi_*\tilde{f}$, and $\tau_q=\tilde{\tau}\circ\phi$. Therefore, the decomposition $\adk{{g_q}}{3}{f_q}={\structfunctB}_q\ad{{g_q}}{f_q}\mod{\GGG}$ implies that $f_q=f^1_q\vec{z}+f^2_q\vec{y}+f^3_q\vec{w}$ satisfies 
    \begin{align}\label{eq:proof_thm_m1_nf_equation_fundamental}
        \frac{\partial^3f^i_q}{\partial w^3}={\structfunctB}(z,y)\frac{\partial f^i_q}{\partial w},
    \end{align}
    \noindent
    for $i=1,2$. Applying the feedback $u=f^3_q(z,y,w) +\tilde{u}$ we obtain the form $\Sigma_q$ {with $f_q=f_q^1\vec{z}+f_q^2\vec{y}$ and $g_q=\vec{w}$}. {The condition $\frac{\partial^2\textsf{f}_q}{\partial w^2}\wedge\frac{\partial \textsf{f}_q}{\partial w}\neq0$ follows from \cref{con:thm_nf_m1_c1} and feedback invariance of the latter.}
\qed\end{proof}

The following corollary shows that we can test on the structure functions of $\Sigma$ if the equivalent quadratic system $\Sigma_q$ will be of elliptic, hyperbolic, or parabolic type. 

\begin{corollary}\label{cor:thm_feedback_quadratization}
    Under conditions \cref{con:thm_nf_m1_c1} and \cref{con:thm_nf_m1_c2} of the previous theorem we have, locally around $\xi_0$,
    \begin{enumerate}[label=(\roman*),ref=\textit{(\roman*)}]
        \item $\Sigma$ is\st{ locally} feedback equivalent to $\Sigma_E$ if and only if $\chi(\xi_0)>0$,
        \item $\Sigma$ is\st{ locally} feedback equivalent to $\Sigma_H$ if and only if $\chi(\xi_0)<0$,
        \item $\Sigma$ is\st{ locally} feedback equivalent to $\Sigma_P$ if and only if $\chi \equiv 0$ in a neighbourhood of $\xi_0$,
    \end{enumerate}
    \noindent
    where $\chi=3\dL{g}{\rho}-2\rho^2-9\tau$.
\end{corollary}
Notice that $\Sigma$ is locally feedback equivalent to $\Sigma_P$ if and only if it satisfies \cref{con:thm_nf_m1_c1} and $\chi\equiv 0$, condition \cref{con:thm_nf_m1_c2} being satisfied automatically.
\begin{proof}\leavevmode
    From the necessity part of the proof of \cref{thm:feedback_quadratization_m_1} we know that for $\Sigma_q${, with structure functions $\rho=0$ and $\tau=\tau_q$,} we have ${\chi}=-9{\structfunctB}_q$ and we \st{saw}{observed} that under pure feedback transformations $(\alpha,\beta)$ we have $\tilde{\chi}=\beta^2\chi$, thus the sign of $\chi$ is invariant as well as the locus where it vanishes. Moreover, {by statement \cref{prop:three_normal_forms_quadratic_systems:2} of \cref{prop:three_normal_forms_quadratic_systems},} $\Sigma_q$ is elliptic if ${\tau_q} > 0$, equivalently ${\chi}<0$, $\Sigma_q$ is hyperbolic if ${\structfunctB_q}>0$, equivalently ${\chi}<0$, and $\Sigma_q$ is parabolic if ${\structfunctB_q}\equiv0$, equivalently ${\chi}\equiv 0$. Hence the necessity of the stated conditions is established. 
    
    Conversely, in the sufficiency part of the proof of \cref{thm:feedback_quadratization_m_1} we obtained {$\Sigma$ with} structure functions $(\tilde{\structfunctA},\tilde{\structfunctB})=(0,\tilde{\structfunctB}(z,y))$ {via a suitable} feedback transformations. Since $\tilde{\chi}=\beta^2\chi$, we have $-9\tilde{\structfunctB}=\beta^2\chi$ and thus we get $\sgn{\tilde{\structfunctB}} = -\sgn{\chi}$ and the conclusion follows by statement \cref{prop:three_normal_forms_quadratic_systems:2} of \cref{prop:three_normal_forms_quadratic_systems}.
\qed\end{proof}

\subsection{Normal form of quadratizable control-affine system}
{Any control-affine system $\Sigma$,} under the regularity assumption $g\wedge\ad{g}{f}(\xi_0)\neq0$, can be written (after applying a suitable feedback transformation) locally around $0\in\RR^3$ as
\begin{align*}
    \Sigma_h\,:\, \left\{\begin{array}{cl}
        \dot{{z}} &= h(z,y,w) \\
        \dot{{y}} &= w +\varepsilon\\
        \dot{{w}} &= u
    \end{array}\right. ,
\end{align*}
\noindent
with $h$ a smooth function {and $\varepsilon=0,1$. The parameter $\varepsilon$, which is invariant under local feedback transformations, is $\varepsilon=0$ if either $\left(f\wedge g\wedge\ad{g}{f}\right)(\xi_0)\neq0$ or $\left(f\wedge g\right)(\xi_0)=0$ and $\varepsilon=1$ otherwise, i.e. $\left(f\wedge g\right)(\xi_0)\neq0$ but $\left(f\wedge g\wedge \ad{g}{f}\right)(\xi_0)=0$}. By applying \cref{thm:feedback_quadratization_m_1}{, we will give in this subsection}\st{ we are able to give} a normal form of all smooth functions $h(z,y,w)$ that describe quadratizable systems, that is, control-affine systems feedback equivalent to $\Sigma_q$. In what follows, we assume to work locally around \st{$0\in\RR^3$}{$0\in\RR^3$} and all derivatives are taken with respect to $w${ and denoted by prime, double prime etc}. Whenever we apply $\ln(a)$, we assume that $a>0$ (if not, take the absolute value).
\begin{theorem}[Normal form of quadratizable control-affine systems]\label{thm:m1_normal_form_quadratizable_systems}
The following statements are equivalent{, locally around $0\in\RR^3$}:
\begin{enumerate}[label=(\roman*), ref=\textit{(\roman*)}]
    \item \label{prop:all-quadratization-submanifold-c1} $\Sigma_h$ is feedback equivalent to a quadratic system $\Sigma_q$;
    \item \label{prop:all-quadratization-submanifold-c2} The function $h$ satisfies $h''(0)\neq0$ and, in a neighbourhood, it holds
    \begin{align}
        \label{eq:all-quadratization-submanifold-full-equation-h}
        9h^{(5)}\left(h''\right)^2-45h^{(4)}h^{(3)}h''+40\left(h^{(3)}\right)^3=0,
    \end{align}
    \noindent
    {recall that the derivatives are taken with respect to $w$;}
    \item \label{prop:all-quadratization-submanifold-c3} The second derivative of $h$ is of the following form 
    \begin{align}\label{eq:all-quadratization-submanifold-form-h2}
        h''(x,w)=a(d w^2+e w+1)^{-3/2},
    \end{align}
    \noindent
    where $a=a(x)$, $d=d(x)$, and $e=e(x)$ are smooth functions satisfying $a(0)\neq0$;
    \item \label{prop:all-quadratization-submanifold-c4} The function $h$ is given by
    \begin{align}\label{eq:all-quadratization-submanifold-form-h}
        h(x,w) = 2a\left(\frac{w^2}{(\sqrt{d w^2+e w+1}+1)^2-d w^2}\right)+b w+c
    \end{align}
    \noindent
    where $a$, $b$, $c$, $d$, $e$ are any smooth functions of $x$ such that $a(0)\neq0$.
\end{enumerate}
\end{theorem}
\begin{proof}
    \cref{prop:all-quadratization-submanifold-c1}$\Rightarrow$\cref{prop:all-quadratization-submanifold-c2}. It is a straightforward application of the conditions of \cref{thm:feedback_quadratization_m_1} with the structure functions of $\Sigma_h$ given by $\structfunctA=\frac{h^{(3)}}{h''}$ and $\structfunctB=0$ {yielding $\chi=3\rho'-2\rho^2$}. By \cref{con:thm_nf_m1_c1}, we have $h''(0)\neq0$ and then condition \cref{con:thm_nf_m1_c2} reads
    \begin{align}\label{eq:prf-all-quadratization-submanifold-rho}
        \chi'-\frac{2}{3}\rho\chi=3\rho''-6\rho\rho'+\frac{4}{3}\rho^3=0,
    \end{align}
    \noindent
    which, by plugging $\rho=\frac{h^{(3)}}{h''}$ into the last equation, gives  \cref{eq:all-quadratization-submanifold-full-equation-h}. \\ 
    \cref{prop:all-quadratization-submanifold-c2}$\Rightarrow$\cref{prop:all-quadratization-submanifold-c3}. Assume that $h$ satisfies $h''(0)\neq0$ and \cref{eq:all-quadratization-submanifold-full-equation-h}. Set $\rho=\frac{h^{(3)}}{h''}$, then $\rho=\rho(x,w)$ fulfils $3\rho''-6\rho\rho'+\frac{4}{3}\rho^3=0$, namely, the second equation of \cref{eq:prf-all-quadratization-submanifold-rho}.
    By a change of variable, it is easy to obtain that the solutions of \cref{eq:prf-all-quadratization-submanifold-rho} are of the following form (see \cref{appendix:m1_solution_lineard_equation} for the proof)
    \begin{align}\label{eq:prf-all-quadratization-submanifold-rho-form}
        \rho(x,w)=-\frac{3}{2}\frac{2{d}(x)w+e(x)}{{d}(x)w^2+e(x)w+1}.
    \end{align}
    This form can be integrated{, using $\rho=\frac{h^{(3)}}{h''}=(\ln h'')'$,} into $h''(x,w)=a(x)\left(d(x)w^2+e(x)w+1\right)^{-3/2}$ with $a$, $d$, and $e$ any smooth functions such that $a(0)\neq0$.\\
    \cref{prop:all-quadratization-submanifold-c3}$\Rightarrow$\cref{prop:all-quadratization-submanifold-c4}. To show \cref{eq:all-quadratization-submanifold-form-h}, we integrate twice the second derivative of $h$ given by \cref{eq:all-quadratization-submanifold-form-h2}. Denote $p=p(x,w)=d(x)w^2+e(x)w+1$ and $\Delta =\Delta(x)=e(x)^2-4d(x)$. First, we obtain (see \cref{appendix:m1_smooth_derivative_h} for details)
    \begin{align*}
        h'(x,w) &= \frac{2aw\left(\sqrt{p}+1\right)}{\sqrt{p}(ew+2+2\sqrt{p})}+b,
    \end{align*}
    \noindent
    with $b$ an arbitrary smooth function of $x$. Integrate once more to get
    \begin{align*}
        h(x,w)&= \frac{2a}{\Delta\sqrt{p}}\left(ew\sqrt{p}-2p\right) +\frac{4a}{\Delta} +bw+ c =\frac{2a}{\Delta}\left(ew+2-2\sqrt{p}\right)+bw+c, \\
        & = \frac{2aw^2}{ew+2+2\sqrt{p}}+bw+c =\frac{2aw^2}{(\sqrt{p}+1)^2-dw^2}+bw+c.
    \end{align*}
    \cref{prop:all-quadratization-submanifold-c4}$\Rightarrow$\cref{prop:all-quadratization-submanifold-c1}. Given $\Sigma_h$ with $h$ defined by \cref{eq:all-quadratization-submanifold-form-h} we will construct a feedback transformation that brings the system into $\Sigma_q$. First, we introduce coordinates, centred at $0\in\RR^2$, $(\tilde{z},\tilde{y})=\phi(z,y)$, where $\tilde{y}=y$, such that $\phi_*\left(b\vec{z}+\vec{y}\right)=\vec{\tilde{y}}$. Those coordinates transform the system $\Sigma_h$ into 
    \begin{align*}
        \left\{\begin{array}{cl}
            \dot{\tilde{z}} &= 2\tilde{a} \frac{{w}^2}{(\sqrt{\tilde{p}}+1)^2-\tilde{d}w^2} + \tilde{c} \\
            \dot{\tilde{y}} &= w+\varepsilon\\
            \dot{{w}} &= u
        \end{array}\right. ,
    \end{align*}
    \noindent
    where $\tilde{a}$, $\tilde{c}$, $\tilde{d}$, and $\tilde{p}$ are new functions satisfying $\tilde{a}(0)\neq0$ and $d=\tilde{d}\circ\phi$ and $p =\tilde{p}\circ\phi$. Second, we set $\tilde{w}^2=\frac{w^2}{(\sqrt{\tilde{p}}+1)^2-\tilde{d}w^2}$ or, equivalently, $w=\tilde{w}\left(\tilde{e}\tilde{w}\pm2\sqrt{\tilde{d}\tilde{w}^2+1}\right)$, which brings the above system into {(after applying a suitable feedback along the last component)}
    \begin{align}\label{eq:prf:m1:pre_normal_form_of_conic_systems}
        \tilde{\Sigma}_h\,:\,\left\{\begin{array}{cl}
            \dot{\tilde{z}} &= 2\tilde{a} \tilde{w}^2 + \tilde{c}  \\
            \dot{\tilde{y}} &= \tilde{e}\tilde{w}^2\pm2\tilde{w}\sqrt{\tilde{d}\tilde{w}^2+1}+\varepsilon\\
            \dot{\tilde{w}} &= \tilde{u}
        \end{array}\right. .
    \end{align}
    \noindent
    The structure functions of $\tilde{\Sigma}_h$ are given by $\tilde{\rho}=\frac{-3\tilde{w}\bar{d}}{\bar{d}\tilde{w}^2+1}$ and $\tilde{\tau}=\frac{3\bar{d}}{\bar{d}\tilde{w}^2+1}$, then apply the feedback $\tilde{u}=\beta \bar{u}$, where $\beta(\tilde{x},\tilde{w})=\sqrt{\bar{d}\tilde{w}^2+1}$ (it is a solution of the equation $\frac{\partial\beta}{\partial \tilde{w}}=-\frac{\tilde{\rho}\beta}{3}$), to obtain a new vector field $\bar{g}=\beta\vec{\tilde{w}}$ and structure functions of $(\bar{f},\bar{g})$, where $\bar{f}$ is the drift of $\tilde{\Sigma}_h$, are given by $\bar{\structfunctA}=0$ and $\bar{\structfunctB}=4\bar{d}(x)$. To complete the form, it remains to find a new $\bar{w}=\psi(\tilde{x},\tilde{w})$ such that the diffeomorphism $(\bar{x},\bar{w})=\Psi(\tilde{x},\tilde{w})=(\tilde{x},\psi(\tilde{x},\tilde{w}))$ satisfies $\Psi_*\bar{g}=\Psi_*\beta\vec{\tilde{w}}=\vec{\bar{w}}$. In general, $\psi$ is given in terms of the following integral
    \begin{align*}
        \bar{w}=\psi(\tilde{x},\tilde{w})=\int_0^{\tilde{w}}\frac{1}{\sqrt{1+\tilde{d}\tilde{w}^2}}\,\diff\tilde{w}
    \end{align*}
    \noindent
    and, together with a suitable feedback along $\dot{\bar{w}}$, thus provides a quadratic system $\Sigma_q$.
\qed\end{proof}
\begin{remark}
    This theorem provides a normal form of submanifolds $\SM=\{\dot{z}=s(x,\dot{y})\}$ that are equivalent to a conic submanifold $\SM_q$, namely, they are represented by $\sm(x,\dot{y})=h(x,\dot{y})$ with $h$ as in \cref{eq:all-quadratization-submanifold-form-h}. Moreover, system \cref{eq:prf:m1:pre_normal_form_of_conic_systems} leads to a normal form for all conic submanifold $\SM_q$ (even for those that smoothly pass through $\Delta_2(x)=0${, i.e. from the elliptic to the hyperbolic submanifolds}) and thus completes {the characterisation of }\cref{lem:m1_classification_conic_submanifolds}. Indeed, for \cref{eq:prf:m1:pre_normal_form_of_conic_systems} we have (tildes have been removed for more readability)
    {
    \begin{align*}
        \Sm_q(x,\dot{x})=\left[e(\dot{z}-c)-2a(\dot{y}-\varepsilon)\right]^2-4d(\dot{z}-c)^2-8a(\dot{z}-c),
    \end{align*}
    \noindent
    }for which {$a\neq0$ and} we can compute $\Delta_1=-64a^4\neq0$ and $\Delta_2=-16a^2d$.
\end{remark}
In the last item of the proof we saw that the function $d$ plays an import role for the shape of the transformation, in the following corollary we show that this function is the key of the normal form of quadratizable control systems. 
\begin{corollary}\label{cor:thm:m1_normal_form_quadratizable_systems}
    Assume that $\Sigma_h$ is given by $h$ of the form \cref{eq:all-quadratization-submanifold-form-h}. Then $\Sigma$ is feedback equivalent to $\Sigma_P$, resp. $\Sigma_E$, resp. $\Sigma_H$ if and only if $d\equiv0$, resp. $d<0$, resp. $d>0$. Moreover, the normalizing feedback transformation is given by 
    \begin{align*}
    &\bar{w}=\frac{w}{1+\sqrt{ew+1}}\quad\textrm{for}\quad\Sigma_P,\\
    \textrm{resp.}\quad\sin^2(\sqrt{-d}\bar{w})=\frac{-dw^2}{ew+2+2\sqrt{p}}\quad\textrm{for}\quad&\Sigma_E,\quad \textrm{resp.}\quad\sinh^2(\sqrt{d}\bar{w})=\frac{dw^2}{ew+2+2\sqrt{p}}\quad\textrm{for}\quad\Sigma_H,
    \end{align*}
    \noindent
    where $p=p(x,w)=d(x)w^2+e(x)w+1$.
\end{corollary}
\begin{proof}
    First we show that $\Sigma_h$ is feedback equivalent to $\Sigma_P$, resp. $\Sigma_E$, resp. $\Sigma_H$, if and only if $d\equiv 0$, resp. $d<0$, and resp. $d>0$. From \cref{cor:thm_feedback_quadratization} we know that we have to compute the sign of $\chi$, which is given by $\chi = -9\frac{d}{p}$
    for $\Sigma_h$ (this can easily be deduced from the expression of $\rho$ given by \cref{eq:prf-all-quadratization-submanifold-rho-form}). Since $p(0)=1>0$ we have $\sgn{\chi}=-\sgn{d}$ and thus the conclusion follows.\\ 
    We now show how to explicitly transform $\Sigma_h$ into $\Sigma_P$, resp. $\Sigma_E$, resp. $\Sigma_H$. From the last part of the proof of the previous theorem we know that a suitable parametrisation $\bar{w}$ is given by the following two steps
    \begin{align*}
        \tilde{w}^2=\frac{w^2}{ew+2+2\sqrt{p}}, \quad\textrm{and}\quad\bar{w}=\int_0^{\tilde{w}} \frac{1}{\sqrt{1+d\tilde{w}^2}}\,d\tilde{w}.
    \end{align*}
    \noindent
    Assume $d\equiv0$, then the procedure reduces to the first step only and thus $\bar{w}^2=\tilde{w}^2=\frac{w^2}{ew+2+2\sqrt{ew+1}}=\left(\frac{w}{1+\sqrt{ew+1}}\right)^2$ and we choose $\tilde{w}=\frac{w}{1+\sqrt{ew+1}}$. Assume $d<0$, then the second step of the procedure leads to $\bar{w}=\frac{1}{\sqrt{-d}}\arcsin(\sqrt{-d}\tilde{w})$. Hence a reparametrisation is given by $\sin^2(\sqrt{-d}\bar{w})=\frac{-dw^2}{ew+2+2\sqrt{p}}$. Assume $d>0$, then from the second step of the procedure we have $\bar{w}=\frac{1}{\sqrt{d}}\arcsinh(\sqrt{d}\tilde{w})$. Hence a reparametrisation is given by $\sinh^2(\sqrt{d}\bar{w})=\frac{dw^2}{ew+2+2\sqrt{p}}$.
\qed\end{proof}

\begin{remark}[Interpretation of parametrising functions]
In the normal form \cref{eq:all-quadratization-submanifold-form-h}, there are $5$ parametrising functions\footnote{On \url{https://www.geogebra.org/m/tyb4ygpb} the reader can play with those parameters (the functions $a$, $b$, $c$, $d$, and $e$ become real numbers when fixing $x\in\XXX$).}. However, only $d=d(x)$ and $e=e(x)$ play a significant role in the shape of the submanifold $\SM_q$. Indeed, $a$ is a scaling of the submanifold, $c$ is the value of $h$ at $w=0$, and by an appropriate choice of coordinates {(as in \cref{eq:prf:m1:pre_normal_form_of_conic_systems})} we can always assume that $b\equiv0$. From the above corollary, the role of $d$ is clear: its sign around $x_0=0\in\RR^2$ determines the nature of the submanifold, that is, whether the submanifold is elliptic, hyperbolic, or parabolic. 

The role of the function $e$ is, however, more subtle. Clearly, $h$ is well defined whenever $p>0$ and, for a given $d$, the function $e$ determines the region in which $p>0$ (in particular, whether $h$ is defined globally with respect to $w$ or not). If $d\equiv0$, then $p>0$ holds everywhere ($h$ is defined globally) if and only if $e\equiv0$ that is, $h$ is explicitly given by $h=2aw^2+bw+c$. If $d<0$, then we have $p>0$ only between its roots and the parametrisation is never global. Finally, if $d>0$ then the parametrisation is global if and only if $\Delta < 0$ (where $\Delta$ is the discriminant of $p$), that is $|e| <  2\sqrt{d}$.

\end{remark}

\section{Classification of quadratic systems}\label{sec:classification_conic}
\textbf{\textit{{Notations.}}}
\begin{table}[H]
\centering
\begin{tabularx}{\textwidth}{lX}
%
%
%
{$\Xi_Q=\{\Xi_{\ttiny{EH}}, \Xi_P\}$, $\Xi_{\ttiny{EH}}=\{\Xi_E,\Xi_H\}$} & {First extension of elliptic, hyperbolic, parabolic submanifolds, seen as a control-nonlinear system on $\XXX$.} \\
{$\Xi_Q=(A,B,C)$} & {Triple of vector fields attached to a first extension of a quadratic submanifold\st{, satisfying $A\wedge B\neq0$}.} \\
{$(\alpha,\beta)$} & {Reparametrisation (feedback) acting on quadratic nonlinear systems, given by \cref{prop:reparam_m1_quadratic_systems}.} \\
{$(\mu_0,\mu_1)$, $\gamma=(\gamma_0,\gamma_1)$} & {Structure functions attached to $(A,B,C)$ by $\lb{A}{B}=\mu_0A+\mu_1B$ and $C=\gamma_0A+\gamma_1B$; see \cref{eq:definition_structure_function}.} \\
{$\antigamma$} & {Notation for $\antigamma=(\gamma_1,\mp\gamma_0)$, used in the elliptic and hyperbolic cases.}\\
{$\Gamma_{E}$, $\Gamma_H$, $\Gamma_P$} & {Functions given by $\Gamma_E=(\gamma_0)^2+(\gamma_1)^2$, $\Gamma_H=(\gamma_0)^2-(\gamma_1)^2$, and $\Gamma_P=\gamma_0+(\gamma_1)^2$.} \\
{$\textsf{g}_{\pm}$, $\kappa_{\pm}$} & {A (pseudo-)Riemanian metric defined by $\textsf{g}_{\pm}(A,A)=1$, $\textsf{g}_{\pm}(B,B)=\pm 1$, and $\textsf{g}_{\pm}(A,B)=0$, and $\kappa_{\pm}$ its Gaussian curvature; see \cref{eq:m1:definition_metric_AB}.}\\
\end{tabularx}%
\end{table}%
From \cref{thm:feedback_quadratization_m_1} we know how to characterise control-affine system equivalent to the quadratic form $\Sigma_q$ and, in particular, we know how to characterise the subclasses {of elliptic, hyperbolic, and parabolic systems}\st{$\Sigma_E$, $\Sigma_H$, and $\Sigma_P$} (see \cref{cor:thm_feedback_quadratization}). We are now interested in classifying, under feedback transformations, the systems inside those three subclasses {because, due to \cref{prop:equiv_of_equiv}, their classification provides a classification of the elliptic, hyperbolic, and parabolic submanifolds that they parametrise (see \cref{lem:classification_ehp_system_to_submanifolds} below).} To this end, we consider the quadratic nonlinear system 
\begin{align*}
    \Xi_Q\,:\,\dot{x}=\textsf{f}_Q(x,w),
\end{align*}
\noindent
where $x\in\XXX$ is the $2$-dimensional state, $w\in\RR$ plays the role of a control that enters in a nonlinear way and $\textsf{f}_Q$ is {a $w$-parameterised vector field on $\XXX$} given by either 
\begin{align*}
    \begin{array}{rll}
        \textsf{f}_E &= A(x)\cos(w) + B(x)\sin(w) + C(x),   &\quad\textrm{defining}\;\;\Xi_E,\quad \textrm{or} \\ 
        \textsf{f}_H &= A(x)\cosh(w)+ B(x)\sinh(w) + C(x),  &\quad \textrm{defining}\;\;\Xi_H,\quad  \textrm{or}\\
        \textsf{f}_P &= A(x)w^2  +B(x)w + C(x),             &\quad\textrm{defining}\;\;\Xi_P.
    \end{array}
\end{align*}
\noindent
In each of the three cases, $A$, $B$, and $C$ are smooth vector fields on $\XXX$ satisfying $(A\wedge B)(x_0)\neq 0$. {We call $\Xi_E$ an elliptic system, $\Xi_H$ a hyperbolic system, and $\Xi_P$ a parabolic system, because in each fiber $T_x\XXX$, the system $\Xi_E$, resp. $\Xi_H$, resp. $\Xi_H$ parametrises an ellipse, resp. a hyperbola, resp. a parabola.} A quadratic nonlinear system $\Xi_Q$ is then represented by the triple $(A,B,C)$ of three smooth vector fields satisfying $A\wedge B\neq 0$. We call the pair $(A,B)$ a \emph{Q-frame}, and if additionally $\lb{A}{B}=0$, then we call $(A,B)$ a \emph{commutative Q-frame}. We will denote by the index $EH$ object attached to either the elliptic or the hyperbolic case, as those two are treated in a similar manner.  \\
\noindent
{For quadratic submanifolds, elliptic and hyperbolic, of the form $\SM_{\ttiny{EH}}=\left\{a^2(\dot{z}-c_0)^2\pm b^2(\dot{y}-c_1)^2=1\right\}$ and parabolic of the form $\SM_P=\left\{a\dot{y}^2-\dot{z}+b\dot{y}+c=0\right\}$, we distinguished specific classes that are summarised in the following \cref{tab:classification_nomenclature}.}
\begin{table}[H]
    \centering
    \begin{tabular}{c|c||c|c}
        \multicolumn{2}{c||}{Elliptic and hyperbolic submanifolds classification} & \multicolumn{2}{c}{Parabolic submanifolds classification} \\
        \hline
        $a=b$ & conformally-flat & $a=1$ & weakly-flat \\
        $a=b=1$ & flat & $a=1$, $b=0$ & strongly-flat \\
        $a=b=1$, $(c_0,c_1)\in\RR^2$ & constant-form & $a=1$, $b=0$, $c\in\RR$ & constant-form \\
        $a=b=1$, $c_0=c_1=0$ & null-form & $a=1$, $b=0$, and $c=0$ & null-form
    \end{tabular}
    \caption{Classification of elliptic, hyperbolic, and parabolic submanifolds.}
    \label{tab:classification_nomenclature}
\end{table}
{Our goal is to characterise the above types of quadratic submanifolds and we now show in the next lemma that the classification of elliptic, hyperbolic, and parabolic submanifolds $\SM_Q$ presented in \cref{tab:classification_nomenclature} is reflected in properties of the control system $\Xi_Q=(A,B,C)$.
\begin{lemma}\label{lem:classification_ehp_system_to_submanifolds}
    Consider a quadratic submanifold $\SM_Q$ together with its regular parametrisation $\Xi_Q=(A,B,C)$. 
    \begin{enumerate}[label=(\roman*),ref=\textit{(\roman*)}]
        \item\label{lem:classification_ehp_system_to_submanifolds:1} $\SM_{\ttiny{EH}}$ is locally equivalent to a conformally-flat elliptic/hyperbolic submanifold if and only if $\Xi_{\ttiny{EH}}$ is locally feedback equivalent to $\Xi_{\ttiny{EH}}$, whose EH-frame $(A,B)$ is given by $A=r(x)\vec{z}$ and $B=r(x)\vec{y}$ for some nonvanishing function $r(x)$.  
        \item\label{lem:classification_ehp_system_to_submanifolds:2} $\SM_{\ttiny{EH}}$ is locally equivalent to a flat elliptic/hyperbolic submanifold if and only if $\Xi_{\ttiny{EH}}$ is locally feedback equivalent to $\Xi_{\ttiny{EH}}$, whose frame $(A,B)$ is commutative.  
        \item\label{lem:classification_ehp_system_to_submanifolds:3} $\SM_{\ttiny{EH}}$ is locally equivalent to a constant-form elliptic/hyperbolic submanifold if and only if $\Xi_{\ttiny{EH}}$ is locally feedback equivalent to $\Xi_{\ttiny{EH}}$, whose EH-frame $(A,B)$ is commutative and, additionally, $\lb{A}{C}=\lb{B}{C}=~0$. 
        \item\label{lem:classification_ehp_system_to_submanifolds:4} $\SM_{\ttiny{EH}}$ is locally equivalent to a null-form elliptic/hyperbolic submanifold if and only if $\Xi_{\ttiny{EH}}$ is locally feedback equivalent to $\Xi_{\ttiny{EH}}$, whose EH-frame $(A,B)$ is commutative and, additionally, $C=0$.
        \item\label{lem:classification_ehp_system_to_submanifolds:5} $\SM_P$ is locally equivalent to weakly-flat parabolic submanifold if and only if $\Xi_{P}$ is locally feedback equivalent to $\Xi_{P}$, whose P-frame $(A,B)$ is commutative. 
        \item\label{lem:classification_ehp_system_to_submanifolds:6} $\SM_P$ is locally equivalent to strongly-flat parabolic submanifold if and only if $\Xi_{P}$ is locally feedback equivalent to $\Xi_{P}$, whose P-frame $(A,B)$ is commutative and, additionally, $A\wedge C=0$.
        \item\label{lem:classification_ehp_system_to_submanifolds:7}$\SM_P$ is locally equivalent to constant-form parabolic submanifold if and only if $\Xi_{P}$ is locally feedback equivalent to $\Xi_{P}$, whose P-frame $(A,B)$ is commutative and, additionally, $\lb{A}{C}=\lb{B}{C}=0$.
        \item \label{lem:classification_ehp_system_to_submanifolds:8} $\SM_P$ is locally equivalent to null-form parabolic submanifold if and only if $\Xi_{P}$ is locally feedback equivalent to $\Xi_{P}$, whose P-frame $(A,B)$ is commutative and, additionally, $C=0$.
    \end{enumerate}
\end{lemma}
Recall that a general conic submanifold $\SM_q$ is equivalent to an elliptic $\SM_E$, resp. a hyperbolic $\SM_H$, resp. a parabolic $\SM_P$, submanifold if and only the determinant $\Delta_2$ satisfies $\Delta_2>0$, resp. $\Delta_2<0$, resp. $\Delta_2\equiv0$; see \cref{lem:m1_classification_conic_submanifolds}. Therefore the above lemma allows to check equivalence of $\SM_q$ to a submanifold of any of the subclasses listed in \cref{tab:classification_nomenclature}.
\begin{proof}
    It is a straightforward computation to check that for the submanifolds $\SM_Q$ of the indicated forms, their first extensions $\Xi_Q=(A,B,C)$ have the triple $(A,B,C)$ or the $Q$-frame $(A,B)$ satisfying the above stated conditions. Conversely, in local coordinates $x=(z,y)$ in which either $A=r(x)\vec{z}$ and $B=r(x)\vec{y}$, for \cref{lem:classification_ehp_system_to_submanifolds:1}, or $A=\vec{z}$ and $B=\vec{y}$, for \cref{lem:classification_ehp_system_to_submanifolds:2} to \cref{lem:classification_ehp_system_to_submanifolds:8}, the system $\Xi_Q$ is a regular parametrisation of $\SM_Q$ with the desired properties. Now, all items \cref{lem:classification_ehp_system_to_submanifolds:1} to \cref{lem:classification_ehp_system_to_submanifolds:8} follow from \cref{prop:equiv_of_equiv}. 
\qed\end{proof}
The above lemma asserts that to achieve the classification of elliptic, hyperbolic, and parabolic submanifolds presented in \cref{tab:classification_nomenclature}, it is crucial to classify, under feedback transformations, quadratic control systems $\Xi_Q=(A,B,C)$ with the properties presented in \cref{tab:classification_nomenclature_systems}, which will be the goal of the remaining part of this section. In particular, we will show that the characterisation of item \cref{lem:classification_ehp_system_to_submanifolds:1}, resp. item \cref{lem:classification_ehp_system_to_submanifolds:5}, is always satisfied by any $\Xi_{\ttiny{EH}}$ and thus by the corresponding $\SM_{\ttiny{EH}}$, resp. by any $\Xi_P$ and thus by the corresponding $\SM_P$, while the characterisations of the remaining classes of systems, and thus of the corresponding submanifolds, require non-trivial conditions. 
}

\begin{table}[H]
    \centering
    \begin{tabular}{c|l||c|l}
        \multicolumn{2}{c||}{Elliptic and hyperbolic classification} & \multicolumn{2}{c}{Parabolic classification} \\
        \hline
        conformally-flat & $A=r(x)\vec{z}$ and $B=r(x)\vec{y}$ & weakly-flat & $\lb{A}{B}=0$ \\
        flat & $\lb{A}{B}=0$ & strongly-flat & $\lb{A}{B}=0$ and $A\wedge C=0$ \\
        constant-form & $\lb{A}{B}=\lb{A}{C}=\lb{B}{C}=0$, & constant-form & $\lb{A}{B}=\lb{A}{C}=\lb{B}{C}=0$ \\
        null-form & $\lb{A}{B}=0$ and $C=0$ & null-form & $\lb{A}{B}=0$ and $C=0$
    \end{tabular}
    \caption{Reflection of the classification of elliptic, hyperbolic, and parabolic submanifolds in properties of the triple $\Xi_Q=(A,B,C)$.}
    \label{tab:classification_nomenclature_systems}
\end{table}

{Although the systems of the form $\Xi_Q$ are nonlinear with respect to the control $w$, the feedback transformations that preserve this class are not as general as possible. Indeed, feedback transformations which preserve the class of quadratic system $\Xi_Q$ are \emph{affine} (and even of Brockett type, in the case of elliptic and hyperbolic systems) with respect to the control $w$, as assured by the next proposition, which also shows how feedback acts on the triple $(A,B,C)$.}\st{The following proposition gives the class of admissible reparametrisations (pure feedback) for each quadratic class and shows how those transformations act on the triple $(A,B,C)$.}
\begin{proposition}[Reparametrisation of quadratic nonlinear systems]\label{prop:reparam_m1_quadratic_systems}
{Consider two quadratic systems $\Xi_Q$ and $\tilde{\Xi}_Q$ around $(x_0,w_0)$ and $(\tilde{x}_0,\tilde{w}_0)$, respectively.}
    \begin{enumerate}[label=(\roman*),ref=\textit{(\roman*)}]
        \item Two elliptic systems $\Xi_E$ and $\tilde{\Xi}_E$ are {locally} feedback equivalent if and only if there exists a local diffeomorphism $\tilde{x}=\phi(x)$ and a reparametrisation (\st{nonlinear }feedback) $w=\psi(x,\tilde{w})$, given by $\psi=\pm\tilde{w}+\alpha(x)$, and satisfying
        \begin{align}\label{eq:reparam_elliptic_abc_m1}
            \tilde{A}=\phi_*\left(A\cos\alpha+B\sin\alpha\right), \; \tilde{B}=\pm\phi_*\left(-A\sin\alpha+B\cos\alpha\right), \; \tilde{C}=\phi_*\left(C\right).
        \end{align}
        \item Two hyperbolic systems $\Xi_H$ and $\tilde{\Xi}_H$ are {locally} feedback equivalent if and only if there exists a local diffeomorphism $\tilde{x}=\phi(x)$ and a reparametrisation (\st{nonlinear }feedback) $w=\psi(x,\tilde{w})$, given by $\psi=\pm\tilde{w}+\alpha(x)$, and satisfying
        \begin{align}\label{eq:reparam_hyperbolic_abc_m1}
            \tilde{A}=\phi_*\left(A\cosh\alpha+B\sinh\alpha\right), \; \tilde{B}=\pm\phi_*\left(A\sinh\alpha+B\cosh\alpha\right), \; \tilde{C}=\phi_*\left(C\right). 
        \end{align}
        \item Two parabolic systems $\Xi_P$ and $\tilde{\Xi}_P$ are {locally} feedback equivalent if and only if there exists a local diffeomorphism $\tilde{x}=\phi(x)$ and an invertible reparametrisation (\st{nonlinear }feedback) $w=\psi(x,\tilde{w})$, given by $\psi=\alpha(x)+\beta(x)\tilde{w}$ and $\beta(\cdot)\neq0$, and satisfying
        \begin{align}\label{eq:reparam_parabolic_abc_m1}
            \tilde{A}=\phi_*\left(A\beta^2\right), \quad \tilde{B}=\phi_*\left(2A\alpha\beta+B\beta\right), \quad     \tilde{C}=\phi_*\left(C+A\alpha^2+B\alpha\right). 
        \end{align}
    \end{enumerate}
\end{proposition}
\begin{proof}\leavevmode 
    We show the necessity of each statement as the converse implications are immediate. 
    \begin{enumerate}[label=\textit{(\roman*)}]
        \item Assume that $\Xi_E$ and $\tilde{\Xi}_E$ are locally equivalent via a diffeomorphism $\tilde{x}=\phi(x)$ and a reparametrisation $w=\psi(x,\tilde{w})$. Then we have the following relation $\phi_*\textsf{f}_E(x,\psi(x,\tilde{w}))=\tilde{\textsf{f}}_E(\tilde{x},\tilde{w})$, which we differentiate $3$ times with respect to $\tilde{w}$ and using $\frac{\partial^3\tilde{\textsf{f}}_E}{\partial\tilde{w}^3}=-\frac{\partial\tilde{\textsf{f}}_E}{\partial\tilde{w}}$, we conclude the relation $\phi_*\frac{\partial^3}{\partial \tilde{w}^3}\textsf{f}_E=-\phi_*\frac{\partial }{\partial \tilde{w}}\textsf{f}_E$, which translates into
        \begin{align*}
            &A\left(-\psi'''\sin(\psi)+(\psi')^3\sin(\psi)-3\psi'\psi''\cos(\psi)\right) \\
            &+B\left(\psi'''\cos(\psi)-(\psi')^3\cos(\psi)-3\psi'\psi''\sin(\psi)\right)= A\psi'\sin(\psi)-B\psi'\cos(\psi),
        \end{align*}
        \noindent
        where the derivatives are taken with respect to $\tilde{w}$. Since the functions $\cos$ and $\sin$ are linearly independent, we obtain $\psi''=0$ and $(\psi')^2=1$. Thus $\psi(x,\tilde{w})=\pm\tilde{w}+\alpha(x)$. Applying this reparametrisation we obtain relations \cref{eq:reparam_elliptic_abc_m1}.
        \item Exactly the same reasoning, using $\textsf{f}_H$ and the fact $\phi_*\frac{\partial^3}{\partial \tilde{w}^3}\textsf{f}_H=\phi_*\frac{\partial }{\partial \tilde{w}}\textsf{f}_H$, implies that $\psi(x,\tilde{w})=\pm\tilde{w}+\alpha(x)$. Applying $\phi(x)$ and $w=\psi(x,\tilde{w})=\pm\tilde{w}+\alpha$ we obtain relations \cref{eq:reparam_hyperbolic_abc_m1}
        \item We repeat again the same reasoning to $\textsf{f}_P$ with the property $\phi_*\frac{\partial^3}{\partial \tilde{w}^3}\textsf{f}_P=0$. However, this time we obtain the conditions $\psi'''=0$ and $\psi\psi'''+3\psi'\psi''=0$ on the reparametrisation $\psi$, which implies $\psi''=0$, that is $\psi(x,\tilde{w})=\beta(x)\tilde{w}+\alpha(x)$, with $\beta$ satisfying $\beta(\cdot)\neq0$. Applying this reparametrisation together with a diffeomorphism $\phi$ yields relations \cref{eq:reparam_parabolic_abc_m1}.
    \end{enumerate}
\qed\end{proof}
\begin{remark}[{Local character of the results}]
    Initially, $\Xi_Q$ was considered {locally around a point $x_0$ and a control $w_0$}, however, {since $\Xi_Q$ is defined globally with respect to $w$ and, moreover,} by the last proposition, the transformations $w=\psi(x,\tilde{w})$ are \emph{global} with respect to $w$, so we will consider the systems $\Xi_Q$ {and their equivalence} locally in $x$ and globally with respect to $w$. All results below are stated assuming this structure.
\end{remark}
We will develop relations involving structure functions attached to any fixed triple $(A,B,C)$ in a unique way and thus change accordingly with diffeomorphisms $\tilde{x}=\phi(x)$. So we will act on $(A,B,C)$ by $(\alpha,\beta)$ only ($\beta$ is $\pm1$ in the elliptic and hyperbolic cases) and we will denote by $(\tilde{A},\tilde{B},\tilde{C})$ the result of that action (given by \cref{eq:reparam_elliptic_abc_m1}, or \cref{eq:reparam_hyperbolic_abc_m1}, or \cref{eq:reparam_parabolic_abc_m1}, with $\phi=\text{id}$), called a reparametrisation. \\
Observe that the reparametrisations of $\Xi_P$ depend on two smooth functions $\alpha$ and $\beta$ while those of $\Xi_E$ and $\Xi_H$ depend on one smooth function $\alpha$ only. Therefore we expect the classification of parabolic systems to be less rich (less parametrising functions) than the classification of elliptic and hyperbolic systems.\st{ For the elliptic (resp. hyperbolic) case, in order to avoid unnecessary computations, we assume that E-frames, resp. H-frames, $(A,B)$ and $(\tilde{A},\tilde{B})$ of two equivalent systems have the same orientation (we will come back to this simplification in \cref{prop:m1:charact_class_eh_flat_systems}). Thus we restrict reparametrisations \cref{eq:reparam_elliptic_abc_m1}, resp. \cref{eq:reparam_hyperbolic_abc_m1}, to those with $\beta=1$}. In the following subsections we will first classify elliptic and hyperbolic systems as the procedures are similar, and then we will classify parabolic systems under reparametrisation actions.
\subsection{Classification of elliptic and hyperbolic systems}\label{subsec:classification_eh}
In this subsection, we classify elliptic and hyperbolic systems under the action of reparametrisations. {Recall, that our aim is to classify elliptic and hyperbolic nonholonomic constraints (submanifolds of $T\XXX$) $\SM_E$ and $\SM_H$, which are parameterised by systems of the form $\Xi_E$ and $\Xi_H$, respectively. The classification of submanifolds given in \cref{tab:classification_nomenclature} is reflected in special properties of the vector fields $(A,B,C)$, attached to the control system $\Xi_E$ and $\Xi_H$, that we list in \cref{tab:classification_nomenclature_systems} above and summarise in \cref{lem:classification_ehp_system_to_submanifolds}.} Firstly, we will give a normal form for both types of systems $\Xi_E$ and $\Xi_H$ showing that they actually depend on three smooth functions{, that normal form corresponds to conformally-flat elliptic and hyperbolic submanifolds.} Secondly, we will further develop their classification, in particular we will give conditions for the existence of commutative frames {(corresponding to flat elliptic and hyperbolic submanifolds)} and a complete characterisation of forms without functional parameters{, corresponding to constant-form (and, in particular, null-form) elliptic and hyperbolic submanifolds}. 
\paragraph{Notations.}In order to simplify {and unify} notations, in the following formulae the upper sign always corresponds to the elliptic case and the lower sign to the hyperbolic case, e.g. we will use the symbol $\pm$ to design similar objects attached to the elliptic ($+$ case) and to the hyperbolic ($-$ case) systems and in the case of a $\mp$ symbol we have $-$ for elliptic systems and $+$ for hyperbolic ones. We denote $\Xi_{\ttiny{EH}}$ elliptic and hyperbolic systems, and an EH-frame stands for an E-frame or an H-frame. {To avoid unnecessary computations, we assume that EH-frames $(A,B)$ and $(\tilde{A},\tilde{B})$ of two equivalent systems have the same orientation (we will come back to this simplification in \cref{prop:m1:charact_class_eh_flat_systems}), therefore we restrict reparametrisations of the control $w$, given by \cref{eq:reparam_elliptic_abc_m1,eq:reparam_hyperbolic_abc_m1} to those with $\beta=1$}. Denoting by $\bar{R}_{\ttiny{EH}}(\alpha)$ the (trigonometric or hyperbolic) rotation matrix given by
\begin{align*}
    \bar{R}_{\ttiny{E}}(\alpha)=\begin{pmatrix}\cos(\alpha) & -\sin(\alpha) \\ \sin(\alpha) &\cos(\alpha) \end{pmatrix},\quad\textrm{and}\quad \bar{R}_{\ttiny{H}}(\alpha)=\begin{pmatrix}\cosh(\alpha) & -\sinh(\alpha) \\ -\sinh(\alpha) &\cosh(\alpha) \end{pmatrix},
\end{align*}
\noindent
respectively, we see from \cref{eq:reparam_elliptic_abc_m1} and \cref{eq:reparam_hyperbolic_abc_m1} that {EH-}frames are transformed by $(\tilde{A},\tilde{B})=(A,B)\bar{R}_{\ttiny{EH}}(\pm\alpha)$ under reparametrisations {of the form $w=\tilde{w}+\alpha$}. Introduce structure functions $(\mu_0,\mu_1)$ and $(\gamma_0,\gamma_1)$ uniquely defined by 
\begin{align}\label{eq:definition_structure_function}
\lb{A}{B}=\mu_0A+\mu_1B\quad\textrm{and}\quad C=\gamma_0A+\gamma_1B,  
\end{align}
\noindent
respectively. We denote $\gamma=(\gamma_0,\gamma_1)$, and $\antigamma=(\gamma_1,\mp\gamma_0)$, and set $\Gamma_{\ttiny{EH}}=(\gamma_0)^2\pm(\gamma_1)^2$. \\
We begin by a technical lemma showing how structure functions behave under reparametrisations {of the control $w$}. 
\begin{lemma}[Transformation of structure functions]\label{lem:m1:reparam_struct_funct_ellip_hyper}
    Consider an elliptic/hyperbolic system $\Xi_{\ttiny{EH}}$ with structure functions $(\mu_0,\mu_1,\gamma_0,\gamma_1)$. Then under the reparametrisation $w=\tilde{w}+\alpha(x)$ we have
    \begin{align}\label{eq:lem:m1:reparam_struct_funct_EH}
        (\tilde{\mu}_0,\tilde{\mu}_1) &=\left(\mu_0\mp\dL{A}{\alpha} , \mu_1-\dL{B}{\alpha} \right) \bar{R}_{\ttiny{EH}}(\alpha),\quad \textrm{and}\quad \quad\tilde{\gamma}=\gamma\,\bar{R}_{\ttiny{EH}}(\alpha).
    \end{align}
\end{lemma}
\begin{proof}
    Details of the computations can be found in \cref{apdx:m1:reparam_struct_funct_ellip_hyper}.
\qed\end{proof}
Clearly from \cref{eq:lem:m1:reparam_struct_funct_EH}, $\Gamma_{\ttiny{EH}}=(\gamma_0)^2\pm(\gamma_1)^2$ is invariant under reparametrisations{, i.e. $\tilde{\Gamma}_{\ttiny{EH}}=\Gamma_{\ttiny{EH}}$}.
\begin{proposition}[Conformal form of elliptic and hyperbolic systems]\label{prop:m1:equivalence_prenormal_elliptic_hyperbolic}
    \begin{enumerate}[label=(\roman*),ref=\textit{(\roman*)}]
        \item Any elliptic system $\Xi_E$, resp. hyperbolic system $\Xi_H$, always admits under a reparametrisation the following conformal form, {locally around $x_0$,}
        \begin{align*}
            \Xi_E^{c}\,:\,\dot{x}&=r(x)\begin{pmatrix}\cos(w) \\ \sin(w)\end{pmatrix}+\begin{pmatrix}c_0(x) \\ c_1(x) \end{pmatrix}, 
            \textrm{  resp.  }\Xi_H^{c}\,:\,\dot{x}=r(x)\begin{pmatrix}\cosh(w)\\ \sinh(w)\end{pmatrix}+\begin{pmatrix}c_0(x) \\ c_1(x) \end{pmatrix},
        \end{align*}
        \noindent
        with $r$ { a smooth function satisfying} $r>0$.
        \item Two conformal forms $\Xi_E^{c}$ and $\tilde{\Xi}_E^{c}$, resp. $\Xi_H^{c}$ and $\tilde{\Xi}_H^{c}$, are locally feedback equivalent if and only if there exists a local diffeomorphism $\tilde{x}=\phi(x)=(\phi_1(x),\phi_2(x))$, where $x=(z,y)$, satisfying 
        \begin{align}\label{eq:prop:m1:equivalence_prenormal_elliptic_hyperbolic}
            \frac{\partial \phi_1}{\partial z} = \frac{\partial \phi_2}{\partial y},\quad \frac{\partial \phi_1}{\partial y} = \mp \frac{\partial \phi_2}{\partial z},\quad \left(\frac{\partial \phi_1}{\partial z}\right)^2 \pm\left(\frac{\partial \phi_1}{\partial y}\right)^2=\left(\frac{\tilde{r}}{r}\right)^2,\; \textrm{and}\quad \phi_*C=\tilde{C}.
        \end{align}
    \end{enumerate}
\end{proposition}
{We call $\Xi^{c}_{\ttiny{EH}}$ a conformal-form because the systems of that class parametrise elliptic and hyperbolic submanifolds for which the quadratic term, interpreted as a (pseudo-)Riemannian metric, is conformally flat.}
%
\begin{proof}\leavevmode
    \begin{enumerate}[label=\textit{(\roman*)}]
        \item For the system $\Xi_{\ttiny{EH}}=(A,B,C)$, define a (pseudo-)Riemannian metric $\textsf{g}_{\pm}$ on $\XXX$ by $\textsf{g}_{\pm}(A,A)=1$, $\textsf{g}_{\pm}(B,B)=\pm1$, and $\textsf{g}_{\pm}(A,B)=0$. It is known that any non-degenerate metric on a manifold of dimension two is conformally flat (see \cite[pp 15-35]{bers1958RiemannsurfacesLectures} or \cite[Addendum 1 of chapter 9]{spivak1999ComprehensiveIntroductionDifferentiala} for the elliptic case and \cite[theorem 7.2]{schottenloher2008MathematicalIntroductionConformal} for the hyperbolic case). Therefore, there exists a local diffeomorphism $(\tilde{z},\tilde{y})=\tilde{x}=\phi(x)$ such that $\textsf{g}_{\pm}=\phi^*\tilde{\textsf{g}}_{\pm}$, where $\tilde{\textsf{g}}_{\pm}=\varrho(\tilde{x})\left(d\tilde{z}^2\pm d\tilde{y}^2\right)$, $\varrho>0$. The vector fields $\tilde{A}=\phi_*A$ and $\tilde{B}=\phi_*B$ satisfy $\tilde{\textsf{g}}_{\pm}(\tilde{A},\tilde{A})=1$, $\tilde{\textsf{g}}_{\pm}(\tilde{B},\tilde{B})=\pm1$, and $\tilde{\textsf{g}}_{\pm}(\tilde{A},\tilde{B})=0$ which implies that $(\tilde{A},\tilde{B})$ is a (pseudo-)orthonormal frame. Finally, using feedback $\alpha$ we can smoothly \emph{rotate} $(\tilde{A},\tilde{B})$ into $\left(r\vec{\tilde{z}},r\vec{\tilde{y}}\right)$ with $r=\frac{1}{\sqrt{\varrho}}$, which gives the desired form ${\Xi}_{\ttiny{EH}}^{{c}}=(\tilde{A},\tilde{B},\tilde{C})$, with $\tilde{C}=\phi_*C$.
        \item By relations \cref{eq:reparam_elliptic_abc_m1} and \cref{eq:reparam_hyperbolic_abc_m1}, the reparametrisations do not act on $C$ and thus the relation $\tilde{C}=\phi_*C$ is necessary for the equivalence of conformal forms. Consider two elliptic conformal systems $\Xi_E^{c}$ and $\tilde{\Xi}_E^{c}$, resp. two hyperbolic conformal systems $\Xi_H^{c}$ and $\tilde{\Xi}_H^{c}$, with frames $(A,B)$ and $(\tilde{A},\tilde{B})$ and related by a feedback $w=\tilde{w}+\alpha$ and a diffeomorphism $\phi$. Thus, using relation \cref{eq:reparam_elliptic_abc_m1}, resp. \cref{eq:reparam_hyperbolic_abc_m1}, we obtain
        \begin{align*}
            \frac{\partial \phi_1}{\partial z}&=\frac{\tilde{r}}{r}\cos(\alpha)=\frac{\partial \phi_2}{\partial y},\quad \frac{\partial \phi_1}{\partial y}=\frac{\tilde{r}}{r}\sin(\alpha)=-\frac{\partial \phi_2}{\partial z}, \\ 
            \textrm{resp. }\frac{\partial \phi_1}{\partial z}&=\frac{\tilde{r}}{r}\cosh(\alpha)=\frac{\partial \phi_2}{\partial y},\quad \frac{\partial \phi_1}{\partial y}=-\frac{\tilde{r}}{r}\sinh(\alpha)=\frac{\partial \phi_2}{\partial z},
        \end{align*}
        \noindent
        from which we deduce condition \cref{eq:prop:m1:equivalence_prenormal_elliptic_hyperbolic}. Conversely, applying the diffeomorphism $\phi$ given by \cref{eq:prop:m1:equivalence_prenormal_elliptic_hyperbolic}, together with the feedback $w=\tilde{w}+\alpha$, with $\alpha$ being a solution of 
            \begin{align*}
                \cos(\alpha)&=\frac{r}{\tilde{r}}\frac{\partial\phi_1}{\partial z},\quad \sin(\alpha)=-\frac{r}{\tilde{r}}\frac{\partial \phi_1}{\partial y},\\
                \textrm{resp. }\cosh(\alpha)&=\frac{r}{\tilde{r}}\frac{\partial\phi_1}{\partial z},\quad \sinh(\alpha)=-\frac{r}{\tilde{r}}\frac{\partial \phi_1}{\partial y}, 
            \end{align*}
            \noindent
            we transform $\Xi_E^{c}$ into $\tilde{\Xi}_E^{c}$, resp. $\Xi_H^{c}$ into $\tilde{\Xi}_H^{c}$.
    \end{enumerate}
\qed\end{proof}
\noindent
\begin{remark}
     In the above proof we used the metric $\textsf{g}_{\pm}$ on $\XXX$ defined by 
     \begin{align}\label{eq:m1:definition_metric_AB}
         \textsf{g}_{\pm}(A,A)=1,\quad \textsf{g}_{\pm}(B,B)=\pm1,\quad \textsf{g}_{\pm}(A,B)=0.
     \end{align}
     \noindent
     This object will play a special role in the interpretation of the conditions describing the existence of a commutative EH-frame.
\end{remark}
The above proposition shows that elliptic and hyperbolic systems  $\Xi_{\ttiny{EH}}$ are parametrized by three smooth functions of two variables (and not by 6 functions defining the triple $(A,B,C)$). Moreover, relation \cref{eq:prop:m1:equivalence_prenormal_elliptic_hyperbolic} shows that the group of diffeomorphisms, conjugating the EH-frames of two given conformal forms, is parametrized by any function $\phi_1$ of two variables satisfying the third equation of \cref{eq:prop:m1:equivalence_prenormal_elliptic_hyperbolic}, $\phi_2$ being given in terms of $\phi_1$ via the first and second equations. 
Now we will pass to the problem of commutative frames and the following proposition gives equivalent algebraic and geometric conditions for the existence of a commutative EH-frame.
\begin{proposition}[Existence of a commutative EH-frame] \label{prop:m1:EH:existence_commutative_frame}\leavevmode
    Consider an elliptic/hyperbolic system $\Xi_{\ttiny{EH}}=(A,B,C)$  with structure functions $(\mu_0,\mu_1)$ of the EH-frame $(A,B)$. The following statements are equivalent {locally around $x_0$}:
    \begin{enumerate}[label=(\roman*),ref=\textit{(\roman*)}]
        \item \label{prop:m1:EH:existence_commutative_frame:1} There exists a commutative EH-frame.
        \item \label{prop:m1:EH:existence_commutative_frame:2} The structure functions $(\mu_0,\mu_1)$ attached to the EH-frame $(A,B)$ satisfy
        \begin{align}\label{eq:prop_existence_EH_commutative_frame}
            -(\mu_0)^2 \mp (\mu_1)^2 \pm \dL{A}{\mu_1}-\dL{B}{\mu_0} = 0.
        \end{align}
        \item \label{prop:m1:EH:existence_commutative_frame:3} The Gaussian curvature $\kappa_{\pm}$ of the metric $\textsf{g}_{\pm}$ vanishes.
    \end{enumerate}
\end{proposition}
Notice that item \cref{prop:m1:EH:existence_commutative_frame:1} describes the following normal forms, 
\begin{align*}
    \Xi_E'\,:\,\left\{\begin{array}{cl}
        \dot{z} &= \cos(w)+c_0(x)  \\
        \dot{y} &= \sin(w) + c_1(x) 
    \end{array}\right.,\quad\textrm{and}\quad \Xi_H'\,:\,\left\{\begin{array}{cl}
        \dot{z} &= \cosh(w)+c_0(x)  \\
        \dot{y} &= \sinh(w) + c_1(x) 
    \end{array}\right.,
\end{align*}
\noindent
whose structure functions are $\mu_0=\mu_1=0$, $\gamma_0=c_0$, and $\gamma_1=c_1$. We call $\Xi_E'$ a \emph{flat elliptic system} and $\Xi_H'$ a \emph{flat hyperbolic system}.
\begin{proof}\leavevmode
    The equivalence between \cref{prop:m1:EH:existence_commutative_frame:2} and \cref{prop:m1:EH:existence_commutative_frame:3} is immediate since \cref{eq:prop_existence_EH_commutative_frame} is the Gaussian curvature $\kappa_{\pm}$ of $\textsf{g}_{\pm}$ (details of the computations are in \cref{apdx:m1:gaussian_curvature}). We show that \cref{prop:m1:EH:existence_commutative_frame:1} is equivalent to \cref{prop:m1:EH:existence_commutative_frame:2}. If the EH-frame $(A,B)$ is equivalent via $w=\tilde{w}+\alpha(x)$ to a commutative EH-frame $(\tilde{A},\tilde{B})$, then by \cref{eq:lem:m1:reparam_struct_funct_EH} we immediately have $\dL{A}{\alpha}=\pm\mu_0$ and $\dL{B}{\alpha}=\mu_1$; the integrability condition of this system of first order partial differential equations gives \cref{eq:prop_existence_EH_commutative_frame}. Conversely, consider the system $\Xi_{\ttiny{EH}}=(A,B,C)$ and construct $\alpha$ as a solution of the system $\dL{A}{\alpha}=\pm\mu_0$ and $\dL{B}{\alpha}=\mu_1$, whose solvability is guaranteed by the integrability condition given by \cref{eq:prop_existence_EH_commutative_frame}. Then by \cref{eq:lem:m1:reparam_struct_funct_EH} we see that the resulting EH-frame $(\tilde{A},\tilde{B})$, of the system $\tilde{\Xi}_{\ttiny{EH}}$ obtained by the reparametrisation $w=\tilde{w}+\alpha$, is commutative.
\qed\end{proof}
Notice that when proving \cref{prop:m1:EH:existence_commutative_frame} we have shown that the Gaussian curvature $\kappa_{\pm}$ of  the metric $\textsf{g}_{\pm}$ is given by the left hand side of \cref{eq:prop_existence_EH_commutative_frame}. Moreover relation \cref{eq:lem:m1:reparam_struct_funct_EH} implies that $\kappa_{\pm}$ is invariant under reparametrisations $w=\tilde{w}+\alpha$ and is therefore an equivariant of the feedback transformations of the system $\Xi_{\ttiny{EH}}$.
In the following proposition, we give first a classification of flat elliptic/hyperbolic systems, second we characterise those without functional parameters, i.e. constant-forms, and third we provide a canonical form for the latter. Recall that $\antigamma=(\gamma_1,\mp\gamma_0)$ and that for flat elliptic/hyperbolic systems $\Xi_{\ttiny{EH}}'$ we have $(\gamma_0,\gamma_1)=(c_0,c_1)$ so all statements of the proposition below are actually expressed in terms of structure functions. From now on, we will consider the group of feedback transformations consisting of $\tilde{x}=\phi(x)$ and $w=\pm\tilde{w}+\alpha(x)$. The additional transformation $w=-\tilde{w}+\alpha$ implies $(\tilde{A},\tilde{B})=(A,B)\bar{\bar{R}}_{\ttiny{EH}}(\pm\alpha)$, where $\bar{\bar{R}}_{\ttiny{E}}(\alpha)=\begin{pmatrix}\cos(\alpha)&\sin(\alpha) \\ \sin(\alpha)&-\cos(\alpha)\end{pmatrix}$ and $\bar{\bar{R}}_{\ttiny{H}}(\alpha)=\begin{pmatrix}\cosh(\alpha)&\sinh(\alpha) \\ -\sinh(\alpha)&-\cosh(\alpha)\end{pmatrix}$, and the corresponding structure functions change by, compare \cref{eq:lem:m1:reparam_struct_funct_EH},
\begin{align}
\label{eq:lem:m1:reparam_struct_funct_EH_prime}
    (\tilde{\mu}_0,\tilde{\mu}_1)=-(\mu_0\mp\dL{A}{\alpha},\mu_1-\dL{B}{\alpha})\bar{\bar{R}}_{\ttiny{EH}}(\alpha) \quad\textrm{and}\quad\tilde{\gamma}=\gamma \bar{\bar{R}}_{\ttiny{EH}}(\alpha).
\end{align}
\noindent
\begin{proposition}[Characterisation and classification of flat elliptic/hyperbolic systems]\label{prop:m1:charact_class_eh_flat_systems}\leavevmode
    \begin{enumerate}[label=(\roman*),ref=\textit{(\roman*)}]
        \item \label{prop:m1:charact_class_eh_flat_systems:1} Two flat elliptic systems $\Xi_E'$ and $\tilde{\Xi}_E'$, resp. two flat hyperbolic systems $\Xi_H'$ and $\tilde{\Xi}_H'$, are locally feedback equivalent around $x_0=0\in\RR^2$ if and only if there exists a constant $\alpha\in\RR$ satisfying
        \begin{align}\label{eq:classification_flat_EH_systems}
            R_{\ttiny{EH}}(\pm\alpha)^{-1}C(x)=\tilde{C}\left(R_{\ttiny{EH}}(\pm\alpha)^{-1}x\right),
        \end{align}
        \noindent
        where $R_{\ttiny{EH}}$ stands for either $\bar{R}_{\ttiny{EH}}$ or $\bar{\bar{R}}_{\ttiny{EH}}$.
        \item \label{prop:m1:charact_class_eh_flat_systems:2} An elliptic/hyperbolic system $\Xi_{\ttiny{EH}}$ is locally feedback equivalent to a constant-form, i.e. $\Xi_{\ttiny{EH}}'$ with $(c_0,c_1)\in\RR^2$, if and only if one of the equivalent conditions of \cref{prop:m1:EH:existence_commutative_frame} holds and, additionally,
        \begin{align}\label{eq:m1:prop:EH:normalization_flat_systems_cond}
            \dL{A}{\gamma}+\antigamma\mu_0=0\quad\textrm{and}\quad \dL{B}{\gamma}\pm \antigamma\mu_1=0.
        \end{align}
        \item \label{prop:m1:charact_class_eh_flat_systems:3}  A constant-form elliptic system is always feedback equivalent, locally around $x_0=0\in\RR^2$, to the canonical form
        \begin{align*}
            \Xi_E^{\Gamma_{\ttiny{E}}}\,:\,\left\{\begin{array}{cl}
                \dot{z} &= \cos(w) + \sqrt{\Gamma_{E}} \\
                \dot{y} &= \sin(w)
            \end{array}\right.,
        \end{align*}
        \noindent
        where $\Gamma_{E}=(c_0)^2+(c_1)^2\in\RR$ is an invariant. 
        \item \label{prop:m1:charact_class_eh_flat_systems:4} A constant-form hyperbolic system is always feedback equivalent, locally around $x_0=0\in\RR^2$, to one of the following canonical form
        \begin{align*}
            \begin{array}{ll}
                \Xi_H^{\Gamma_H,\varepsilon}\,:\,\left\{\begin{array}{cl}
                    \dot{z} &= \cosh(w) + \varepsilon\sqrt{\Gamma_H} \qquad \\
                    \dot{y} &= \sinh(w)  
                \end{array}\right.,\;&\textrm{or}\quad 
                \Xi_H^{-\Gamma_H}\,:\,\left\{\begin{array}{cl}
                    \dot{z} &= \cosh(w)  \\
                    \dot{y} &= \sinh(w) + \sqrt{-\Gamma_H}         
                \end{array}\right.,\\
                \textrm{or}\quad \Xi_H^{0,\varepsilon}\,:\,\left\{\begin{array}{cl}
                    \dot{z} &= \cosh(w) + \varepsilon \\
                    \dot{y} &= \sinh(w) + 1
                \end{array}\right.,\;&\textrm{or}\quad 
                \Xi_H^{0,0}\,:\,\left\{\begin{array}{cl}
                    \dot{z} &= \cosh(w)  \\
                    \dot{y} &= \sinh(w)        
                \end{array}\right.,
            \end{array}
        \end{align*}
        \noindent
        where $\Gamma_H=(c_0)^2-(c_1)^2\in\RR$ and satisfies $\Gamma_H>0$ for the first form, $\Gamma_H<0$ for the second form, and $\Gamma_H=0$ for the third and fourth ones, where $\varepsilon=\sgn{c_0}=\pm 1$. Moreover $(\Gamma_H,\varepsilon)$ is a complete invariant.
    \end{enumerate}
\end{proposition}
Observe that if in items \cref{prop:m1:charact_class_eh_flat_systems:1}, \cref{prop:m1:charact_class_eh_flat_systems:3}, and \cref{prop:m1:charact_class_eh_flat_systems:4} the considered systems are defined globally, then their feedback equivalence is also global and, in particular, the proposed canonical forms are also global.
\begin{remark}
     In item \cref{prop:m1:charact_class_eh_flat_systems:4}, notice that there are two orbits of the local action of feedback transformations group for $\Gamma_H>0$, corresponding to $\sgn{c_0}=\varepsilon=\pm1$, one orbit for $\Gamma_H<0$, and three orbits for $\Gamma_H=0$ corresponding, respectively, to $\sgn{c_0}=\varepsilon=\pm1$ or $({c_0},{c_1})=(0,0)$. The invariant $\varepsilon=\pm1$ corresponds to the parametrisation of one of two branches of the hyperbola $(\dot{z}-\sqrt{\Gamma_H})^2-\dot{y}^2=1$.
\end{remark}
\begin{proof}\leavevmode
    \begin{enumerate}[label=\textit{(\roman*)}]
        \item Consider, locally around $0\in\RR^2$, two equivalent flat elliptic/hyperbolic systems $\Xi_{\ttiny{EH}}'$ and $\tilde{\Xi}_{\ttiny{EH}}'$ given by structure functions $({\mu}_0,{\mu}_1, {\gamma}_0, {\gamma}_1)=(0,0,c_0,c_1)$ and $(\tilde{\mu}_0,\tilde{\mu}_1, \tilde{\gamma}_0, \tilde{\gamma}_1)=(0,0,\tilde{c}_0,\tilde{c}_1)$, respectively. Since they both have a commutative EH-frame, by \cref{eq:lem:m1:reparam_struct_funct_EH} (and \cref{eq:lem:m1:reparam_struct_funct_EH_prime}) they differ by a reparametrisation $w=\pm\tilde{w}+\alpha$ satisfying $\dL{A}{\alpha}=\dL{B}{\alpha}=0$ and thus $\alpha\in\RR$. Applying this reparametrisation together with a diffeomorphism $\phi$ satisfying $\phi_*=R_{\ttiny{EH}}(\pm\alpha)^{-1}$, that is $\tilde{x}=\phi(x)=R_{\ttiny{EH}}(\pm\alpha)^{-1}x$, transforms $\Xi_{\ttiny{EH}}'$ into $\tilde{\Xi}_{\ttiny{EH}}'$ if and only if
        \begin{align*}
            \begin{pmatrix} \tilde{c}_0(\tilde{x})\\ \tilde{c}_1(\tilde{x})\end{pmatrix}=R_{\ttiny{EH}}(\pm\alpha)^{-1}\begin{pmatrix} {c}_0(x)\\ {c}_1(x)\end{pmatrix},
        \end{align*}
        \noindent
        which is \cref{eq:m1:prop:EH:normalization_flat_systems_cond}.
        \item Assume that $\Xi_{\ttiny{EH}}$, given by structure functions $(\mu_0,\mu_1,\gamma_0,\gamma_1)$, is equivalent via $\tilde{x}=\phi(x)$ and $w=\tilde{w}+\alpha$ to $\Xi_{\ttiny{EH}}'$ with structure functions $(\tilde{\mu}_0,\tilde{\mu}_1, \tilde{\gamma}_0, \tilde{\gamma}_1) = (0,0,c_0,c_1)$, where $(c_0,c_1)\in\RR^2$. Necessity of one (and thus any) of the conditions of \cref{prop:m1:EH:existence_commutative_frame} is clear, and by \cref{eq:lem:m1:reparam_struct_funct_EH} and \cref{eq:lem:m1:reparam_struct_funct_EH_prime} we have first, $\dL{A}{\alpha}=\pm\mu_0$ and $\dL{B}{\alpha}=\mu_1$ and second, $\gamma R_{\ttiny{EH}}(\alpha)=\tilde{\gamma}=(c_0,c_1)$. By differentiating this last relation along $A$ and $B$ we obtain
        \begin{align*}
            0&=\dL{A}{\gamma}\,R_{\ttiny{EH}}(\alpha)+\gamma\dL{A}{R_{\ttiny{EH}}(\alpha)}  = \dL{A}{\gamma}R_{\ttiny{EH}}(\alpha)+\gamma\left(\pm\dL{A}{\alpha}\begin{pmatrix}0&\mp1\\1&0 \end{pmatrix}R_{\ttiny{EH}}(\alpha)\right), \\
            0&=\dL{A}{\gamma}+\antigamma\mu_0,\quad\textrm{and} \\ 
            0&= \dL{B}{\gamma}\,R_{\ttiny{EH}}(\alpha)+\gamma\dL{B}{R_{\ttiny{EH}}(\alpha)}  = \dL{B}{\gamma}R_{\ttiny{EH}}(\alpha)+\gamma\left(\pm\dL{B}{\alpha}\begin{pmatrix}0&\mp1\\1&0 \end{pmatrix}R_{\ttiny{EH}}(\alpha)\right), \\
            0&= \dL{B}{\gamma}\pm\antigamma\mu_1,
        \end{align*}
        \noindent
        thus proving that \cref{eq:m1:prop:EH:normalization_flat_systems_cond} holds. Conversely, assume that \cref{eq:prop_existence_EH_commutative_frame} and \cref{eq:m1:prop:EH:normalization_flat_systems_cond} hold for $\Xi_{\ttiny{EH}}$. By \cref{prop:m1:EH:existence_commutative_frame}, $\Xi_{\ttiny{EH}}$ is equivalent to $\Xi_{\ttiny{EH}}'$ with a commutative EH-frame $(A,B)$, and applying \cref{eq:m1:prop:EH:normalization_flat_systems_cond} to the latter we get $\dL{A}{\gamma}=\dL{B}{\gamma}=0$ and therefore, we have $(c_0,c_1)\in\RR^2$. 
        \item Consider a flat elliptic system $\Xi_E'$ with $(c_0,c_1)\in\RR^2$, then relation \cref{eq:classification_flat_EH_systems} reads
        \begin{align}\tag{\ref*{eq:classification_flat_EH_systems}'}\label{eq:classification_flat_EH_systems_E}
            \begin{pmatrix} \tilde{c}_0 \\ \tilde{c}_1 \end{pmatrix} &=\begin{pmatrix} \cos(\alpha) & \sin(\alpha) \\ \mp\sin(\alpha) & \pm\cos(\alpha) \end{pmatrix} \begin{pmatrix} {c}_0 \\ {c}_1 \end{pmatrix} .
        \end{align}
        Take $\alpha$ as a solution of $-\sin(\alpha)c_0+\cos(\alpha)c_1=0$, then we have $\tilde{c}_1=0$ and $\tilde{c}_0=\pm\sqrt{\Gamma_E}$, with $\Gamma_E=(c_0)^2+(c_1)^2$. If necessary, apply $\alpha=\pi$ to send $(\tilde{c}_0,\tilde{c}_1)=(-\sqrt{\Gamma_E},0)$ into $(\sqrt{\Gamma_E},0)$. The proof that $\Xi_E^{\Gamma_E}$ is equivalent to $\tilde{\Xi}_E^{\Gamma_E}$ if and only if $\Gamma_E=\tilde{\Gamma}_E$ is immediate from \cref{eq:classification_flat_EH_systems_E}.
        \item Consider a flat hyperbolic system $\Xi_H'$ with $(c_0,c_1)\in\RR^2$ and denote $\Gamma_H=(c_0)^2-(c_1)^2$, then relation \cref{eq:classification_flat_EH_systems} reads 
        \begin{align}\tag{\ref*{eq:classification_flat_EH_systems}''}\label{eq:classification_flat_EH_systems_H}
            \begin{pmatrix} \tilde{c}_0 \\ \tilde{c}_1 \end{pmatrix} &=\begin{pmatrix} \cosh(\alpha) & -\sinh(\alpha) \\ \mp\sinh(\alpha) & \pm\cosh(\alpha) \end{pmatrix} \begin{pmatrix} {c}_0 \\ {c}_1 \end{pmatrix} .
        \end{align}
        \noindent
        We consider four cases. First, assume that $\Gamma_{H} > 0$, that is $c_0\neq0$ and $-1<\frac{c_1}{c_0}<1$, and take $\alpha$ as the solution of $\tanh(\alpha)=\frac{c_1}{c_0}$ yielding $\tilde{c}_1=0$ and $\tilde{c}_0=\sgn{c_0}\sqrt{\Gamma_H}$, which gives the canonical form $\Xi_H^{\Gamma_H,\varepsilon}$. Second, assume that $\Gamma_{H} < 0$, that is $c_1\neq0$ and $-1<\frac{c_0}{c_1}<1$, and take $\alpha$ as the solution of $\tanh(\alpha)=\frac{c_0}{c_1}$ yielding $\tilde{c}_0=0$ and $\tilde{c}_1=\sgn{c_1}\sqrt{-\Gamma_H}$. If $\sgn{c_1}=-1$, then by applying \cref{eq:classification_flat_EH_systems_H} with $\alpha=0$ and the bottom sign, we can always normalize $\sgn{c_1}$ yielding the canonical form $\Xi_H^{-\Gamma_H}$. Third, assume that $\Gamma_H=0$ and $c_0=0$ thus $c_1=0$ and therefore we immediately have the canonical form $\Xi_H^{0,0}$. Fourth, and finally, assume that $\Gamma_H=0$ and $c_0\neq0$, thus $c_1=\varepsilon c_0$ with $\varepsilon=\pm1$. If necessary, apply \cref{eq:classification_flat_EH_systems_H} with $\alpha=0$ and the bottom sign to obtain $c_1>0$. Take $\alpha=\varepsilon\ln{c_1}$ and apply \cref{eq:classification_flat_EH_systems_H} with the upper sign to obtain $\tilde{c}_1=1$ and $\tilde{c}_0=\varepsilon$. To show that $(\Gamma_H,\varepsilon)$ is a complete invariant is trivial by applying \cref{eq:classification_flat_EH_systems_H} to the canonical forms $\Xi_H^{\Gamma_H,\varepsilon}$, $\Xi_H^{-\Gamma_H}$, $\Xi_H^{0,\varepsilon}$, and $\Xi_H^{0,0}$.
    \end{enumerate}
\qed\end{proof}
We summarise the results of this subsection, we started from a general {elliptic, resp. hyperbolic, system $\Xi_E$, resp. $\Xi_H$, which parametrises a elliptic, resp. a hyperbolic, submanifold $\SM_E$, resp. $\SM_H$.}\st{elliptic/hyperbolic system $\Xi_{\ttiny{EH}}$, which depends} {We show that reparametrisations (pure feedback transformations) acting on the class of elliptic/hyperbolic systems $\Xi_{\ttiny{EH}}$ is affine (actually it is of the Brockett type $w=\pm\tilde{w}+\alpha$) with respect to the control and thus depend on one function $\alpha$ only. By transforming into the the conformal form, we show that elliptic/hyperbolic systems $\Xi_{\ttiny{EH}}$ are given by} three {arbitrary} smooth functions{.}\st{, then with the} {Next, we show that the vanishing of the} Gaussian curvature {of a (pseudo-)Riemannian metric}, associated with the EH-frame {of $\Xi_{\ttiny{EH}}$}, \st{being zero we reduced the systems to}{characterises} the flat elliptic/hyperbolic systems $\Xi_{\ttiny{EH}}'$\st{ depending}{, which depend} on two {arbitrary} smooth functions only. \st{And f}{F}inally, we gave conditions characterising the\st{ flat} {constant-form} systems{, i.e. ellitpic/hyperbolic systems} without functional parameters. In the elliptic case, equivalent {constant-form} systems correspond to the circles $\Gamma_E=(c_0)^2+(c_1)^2=\textrm{const.}$, and their canonical forms are parametrized by a closed half-line of real constants. On the other hand, in the hyperbolic case the structure is richer because equivalent systems correspond to connected branches of the hyperbolas $\Gamma_H=(c_0)^2-(c_1)^2=\textrm{const.}$; two connected components for $\Gamma_H>0$, one for $\Gamma_H<0$, and three for $\Gamma_H=0$. Thus canonical forms of hyperbolic systems are parametrized by a real line of constants (the value of $\Gamma_H$) and by a discrete invariant $\varepsilon=\pm1$ (if $\Gamma_H>0$ or $\Gamma_H=0$ and $c_0\neq0$). 
\st{We also obtained a similar characterisation of flat and strongly flat elliptic/hyperbolic submanifolds; the derived classification of those strongly flat submanifolds is similar to the one of their parametrisations. The only notable difference shows up in some cases of hyperbolic submanifolds, where we can have two non-equivalent parametrisations of the two branches of an hyperbola whereas for the hyperbolic submanifold there is no distinction between the branches.}{Our characterisation of elliptic/hyperbolic systems gives an equivalent classification for elliptic/hyperbolic submanifolds that is summarised in \cref{lem:classification_ehp_system_to_submanifolds} and explicitly given in \cite{schmoderer2018Studycontrolsystems}.}

\subsection{Classification of parabolic systems}\label{subsec:classification_p}%
We now turn to the classification of parabolic systems $\Xi_P=(A,B,C)$, of the form 
\begin{align*}
    \Xi_P\,:\,\dot{x}&=A(x)w^2+B(x)w+C(x),
\end{align*}
\noindent 
which is expected to be different from that of elliptic and hyperbolic systems because the allowed reparametrisations {of the control $w$} depend on $2$ smooth functions $(\alpha,\beta)$, see \cref{prop:reparam_m1_quadratic_systems}. {Our aim is to get a classification of parabolic submanifolds $\SM_P$, as the one presented in \cref{tab:classification_nomenclature}, via the properties of the triple $(A,B,C)$ presented in \cref{tab:classification_nomenclature_systems}, see the beginning of \cref{sec:classification_conic} for both tables. In particular, existence of a commutative P-frame $(A,B)$ corresponds to weakly-flat parabolic submanifolds; existence of a commutative P-frame, which additionally satisfies $A\wedge C=0$, corresponds to strongly-flat submanifolds; finally, existence of a commutative P-frame, which additionally satisfies $C$ constant (in the coordinates where $(A,B)$ is rectified), corresponds to constant-form parabolic submanifolds; in particular, the case of $C=0$ is called null-form. }As in the elliptic/hyperbolic cases, we introduce the structure functions $(\mu_0,\mu_1,\gamma_0,\gamma_1)$ uniquely defined for any triple $(A,B,C)$ by
\begin{align*}
 \lb{A}{B}=\mu_0A+\mu_1B\quad \textrm{and} \quad C=\gamma_0A+\gamma_1B.    
\end{align*}
\noindent
By a direct computation, we obtain that under reparametrisations of the form {$w=\beta\tilde{w}+\alpha$ the structure functions $(\mu_0,\mu_1,\gamma_0,\gamma_1)$ of $(A,B,C)$ are transformed intro the structure functions $(\tilde{\mu}_0,\tilde{\mu}_1,\tilde{\gamma}_0,\tilde{\gamma}_1)$ of $(\tilde{A},\tilde{B},\tilde{C})$ via the following relations:}
\begin{align}\label{eq:action_alpha_beta_on_gamma_nf2d}
    &\tilde{\gamma}_0 = \frac{1}{\beta^2}\left(\gamma_0-2\alpha\gamma_1-\alpha^2\right), \quad 
        \tilde{\gamma}_1 = \frac{1}{\beta}\left(\gamma_1+\alpha\right), \\
    \label{eq:reparam_m1_action_on_ab_bracket}
    &\left.\begin{array}{rl}
        \tilde{\mu}_0 &= \beta\mu_0-2\alpha\dL{A}{\beta}+2\beta\dL{A}{\alpha}-2\dL{B}{\beta}-2\alpha\left(\dL{A}{\beta}+\beta\mu_1\right),\\    
        \tilde{\mu}_1 &= \beta^2\mu_1+\beta\dL{A}{\beta}.
    \end{array}\right.
\end{align}
\noindent
There are two main questions that we will answer. First, when does a commutative P-frame $(\tilde{A},\tilde{B})$ exist, i.e. $\tilde{\mu}_0=\tilde{\mu}_1=0$? Second, provided that a commutative P-frame $(A,B)$  has been normalised, how can we additionally simplify $C$? Contrary to the elliptic and hyperbolic cases the answer to the first question is always positive without any additional assumption, as assured by the next result.
\begin{proposition}[Existence of a commutative P-frame]\label{prop:m1_exits_commutative_frame}
    \begin{enumerate}[label={(\roman*)}, ref=\textit{(\roman*)}]
        \item \label{itm:prop_m1_exits_ab_commutative} For any P-frame $(A,B)$ there exists{, locally around $x_0$,} a reparametrisation $(\alpha,\beta)$ such that $(\tilde{A},\tilde{B})$ is a commutative P-frame.
        \item If $(A,B)$ is a commutative P-frame, then $(\tilde{A},\tilde{B})$ is also a commutative P-frame if and only if the reparametrisation $(\alpha,\beta)$ satisfies
        \begin{align}\label{eq:m1_commutative_alpha_beta}
            \dL{A}{\beta}=0\  \textrm{  and  }\  \frac{1}{\beta}\dL{B}{\beta}=\dL{A}{\alpha}.
        \end{align}
    \end{enumerate}
\end{proposition}
\begin{proof}\leavevmode
    \begin{enumerate}[label=\textit{(\roman*)}]
        \item Consider a P-frame $(A,B)$ whose structure functions are $(\mu_0,\mu_1)$. Apply a reparametrisation $(\alpha,\beta)$, $\beta\neq0$, given by a solution of the following {system of} equations (we solve the first equation for $\beta$, and plug in $\beta$ into the second equation to solve it for $\alpha$)
        \begin{align*}
            \left\{\begin{array}{rl}
                \dL{A}{\beta} &= -\beta \mu_1  \\
                \dL{A}{\alpha}+\alpha\mu_1 &= \frac{1}{2\beta}\left(2\dL{B}{\beta}-\beta\mu_0\right)
            \end{array}\right. .
        \end{align*}
        \noindent
        Then formula \cref{eq:reparam_m1_action_on_ab_bracket} implies $\tilde{\mu}_0=\tilde{\mu}_1=0${, i.e. $\lb{\tilde{A}}{\tilde{B}}=0$}. Notice that to ensure $\beta\neq0$, we may actually solve $\dL{A}{\ln\beta}=-\mu_1$.
        \item Using relation \cref{eq:reparam_m1_action_on_ab_bracket} with $\mu_i=\tilde{\mu_i}=0$, for $i=0,1$, we see that all reparametrisations $(\alpha,\beta)$ have to satisfy relation \cref{eq:m1_commutative_alpha_beta}. Conversely, if $(A,B)$ is a commutative P-frame ($\mu_0=\mu_1=0$) and $(\alpha,\beta)$ is any solution of \cref{eq:m1_commutative_alpha_beta}, with $\beta\neq0$, then by \cref{eq:reparam_m1_action_on_ab_bracket} we obtain that $\tilde{\mu}_0=\tilde{\mu}_1=0${, i.e. $(\tilde{A},\tilde{B})$ is also a commutative P-frame}. 
    \end{enumerate}
\qed\end{proof}
\noindent
Immediately, item \cref{itm:prop_m1_exits_ab_commutative} of \cref{prop:m1_exits_commutative_frame} gives the following prenormal forms of {parabolic} systems $\Xi_P$. {Recall that systems and equivalence are considered locally with respect to the state $x$ and globally with respect to the control $w\in\RR$.}
\begin{corollary}[Prenormal forms of $\Xi_P$]\label{cor:m1_prenormal_form}
    The parabolic system $\Xi_P$ is always feedback equivalent to the following prenormal forms{, locally around $x_0$}: 
    \begin{align*}
    \Xi_P'\,:\,\left\{ \begin{array}{rl}
        \dot{z} &= w^2 + c_0(x)  \\
        \dot{y} &= w + c_1(x) 
    \end{array} \right., \qquad\qquad 
    \Xi_P''\,:\,\left\{ \begin{array}{rl}
        \dot{z} &= w^2+ b(x)w + \Gamma_P(x) \\
        \dot{y} &= w 
    \end{array} \right.,
    \end{align*}
    \noindent
    whose structure functions are $(\mu_0',\mu_1',\gamma_0',\gamma_1')=(0,0,c_0,c_1)$, and $(\mu_0'',\mu_1'',\gamma_0'',\gamma_1'')=\left(\frac{\partial b}{\partial z},0,\Gamma_P,0\right)$, respectively.
\end{corollary}
{A parabolic system of the form $\Xi_P'$ is called weakly-flat, this terminology is justified by generalisations of the previous result in higher dimensions that we will present in a future paper.}
\begin{remark}
     Since any parabolic system can be brought into $\Xi_P'$ (and into $\Xi_P''$), it follows that all parabolic systems are locally parametrized by two functions of two variables ($c_0$ and $c_1$, or, equivalently, $b$ and $\Gamma_P$). This is in contrast with elliptic/hyperbolic systems $\Xi_{\ttiny{EH}}$ parametrized by three functions of two variables (compare \cref{prop:m1:equivalence_prenormal_elliptic_hyperbolic}).
\end{remark}
\begin{proof}
    Apply to $\Xi_P$ a reparametrisation $(\alpha,\beta)$ transforming its P-frame into a commutative P-frame $(\tilde{A},\tilde{B})$ and let $\phi$ be a local diffeomorphism introducing coordinates $x=(z,y)$ such that $\phi_*\tilde{A}=\vec{z}$ and $\phi_*\tilde{B}=\vec{y}$. In this system of coordinates, $\Xi_P$ takes the form $\Xi_P'$. Then apply to $\Xi_P'$ the reparametrisation $\tilde{w}=w+c_1(x)$ to obtain the form $\Xi_P''$ with $b=-2c_1$ and $\Gamma_P=c_0+(c_1)^2$. The computation of the structure functions is straightforward.
\qed\end{proof}
Notice that the normal forms $\Xi_P'$ and $\Xi_P''$ are related by the reparametrisation $\tilde{w}=w+c_1(x)$. The function $\Gamma_P$ {(appearing in $\Xi_P''$)} will be of special importance in the remaining part of this section, and in any P-frame $(A,B)$, commutative or not, we define it by setting 
\begin{align*}
  \Gamma_P=\gamma_0+(\gamma_1)^2.
\end{align*}
\noindent
{Recall that $(\gamma_0,\gamma_1)$ are defined by $C=\gamma_0A+\gamma_1B$}. Clearly, diffeomorphisms act on $\Gamma_P$ by conjugation and reparametrisations $(\alpha,\beta)$ act by $\beta^2\tilde{\Gamma}_P=\Gamma_P$ (as it can be computed from formula \cref{eq:action_alpha_beta_on_gamma_nf2d}).\\

The remaining part of this section shows how to additionally normalize $\Xi_P'$ and $\Xi_P''$ while preserving the commutativity of the P-frame $(A,B)$. Although for a parabolic system $\Xi_P$, there always exists a commutative $P$-frame $(A,B)$, its explicit construction can be difficult or even impossible, as it requires to solve a system of first order PDEs. For this reason, we will state our results for a general, not necessarily commutative, P-frame $(A,B)$.
\begin{theorem}[Normalisation of parabolic systems] \label{thm:m1_normalisation_quadratic_system_gamma_1_0}
    Let $\Xi_P=(A,B,C)$ be a parabolic system with structure functions $(\mu_0,\mu_1, \gamma_0,\gamma_1)$. Then the following statements hold{, locally around $x_0$}:
    \begin{enumerate}[label={(\roman*)}, ref=\textit{(\roman*)}]
        \item \label{thm:m1_normalisation_quadratic_system_gamma_1_0:1} $\Xi_P$ is {strongly-flat,}\st{locally} {i.e.} feedback equivalent to $\Xi_P'$ with $c_1\equiv 0${, we will denote that form $\Xi_P'''$,} if and only if
        \begin{align}\label{eq:thm_m1_normalize_quadratic_submanifold_gamma1_0}
            \dLk{A}{2}{\gamma_1}+\gamma_1\left(\dL{A}{\mu_1}-(\mu_1)^2\right)&= \frac{\mu_0\mu_1}{2}+\frac{1}{2}\dL{A}{\mu_0}+\dL{B}{\mu_1}.
        \end{align}
        \item \label{thm:m1_normalisation_quadratic_system_gamma_1_0:2} $\Xi_P$ is {a constant-form, i.e.}\st{locally} feedback equivalent to $\Xi_P'$ with $(c_0,c_1)\in\RR^2$, satisfying $c_0+(c_1)^2\in\RR^*$, if and only if $\Gamma_P \neq 0$ and 
        \begin{align}\label{eq:thm_m1_normalize_quadratic_submanifold_g1_0_g0_pm}
            \left.\begin{array}{rl}
                \dL{A}{\Gamma_P}+2\mu_1\Gamma_P &=0,  \\
                \dL{B}{\Gamma_P}+2\Gamma_P\dL{A}{\gamma_1} &= \Gamma_P\mu_0-2\Gamma_P\gamma_1\mu_1 .
            \end{array}\right.
        \end{align}
        \noindent
        Moreover, in this case we can always normalize $c_0=\pm1$ and $c_1=0$ and we will denote that form $\Xi_P^{\pm}$.
        \item \label{thm:m1_normalisation_quadratic_system_gamma_1_0:3}  $\Xi_P$ is {a null-form, i.e.}\st{locally} feedback equivalent to $\Xi_P'$ with $c_0\equiv c_1\equiv 0$, which we will denote $\Xi_P^0${,} if and only if \cref{eq:thm_m1_normalize_quadratic_submanifold_gamma1_0} holds and, additionally, $\Gamma_P\equiv0$.
    \end{enumerate}
\end{theorem}

Notice that items \cref{thm:m1_normalisation_quadratic_system_gamma_1_0:1}, \cref{thm:m1_normalisation_quadratic_system_gamma_1_0:2}, and \cref{thm:m1_normalisation_quadratic_system_gamma_1_0:3} of the above theorem characterise, respectively, the following normal forms {around $x_0$ (recall that we work globally with respect to $w$)}
\begin{align*}
    \Xi_P'''\,:\,\left\{\begin{array}{cl}
        \dot{z} &= w^2+c_0(x)  \\
        \dot{y} &= w 
    \end{array}\right.,\quad\quad \Xi_P^{\pm}\,:\,\left\{\begin{array}{cl}
        \dot{z} &= w^2\pm 1 \\
        \dot{y} &= w 
    \end{array}\right.,\quad\textrm{and}\quad \Xi_P^{0}\,:\,\left\{\begin{array}{cl}
        \dot{z} &= w^2\\
        \dot{y} &= w 
    \end{array}\right..
\end{align*}
{A parabolic system of the form $\Xi_P'''$ is called strongly-flat (the terminology is justified by generalisations of the previous result to higher dimensions).} Equivalent statements can be formulated to obtain special structures of the second prenormal form $\Xi_P''$, e.g. the normal form $\Xi_P'''$ corresponds to $\Xi_P''$ with $b\equiv0$ and thus describes the intersection of the prenormal forms $\Xi_P'$ and $\Xi_P''$. The {constant-}forms $\Xi_P^{\pm}$ and $\Xi_P^0$ are clearly special cases of $\Xi_P'''$ where $c_0$ is constant. {Observe that contrary to the elliptic/hyperbolic case, where canonical forms of constant-form systems are parametrised by real-continuous parameters, in the parabolic case, there are only three distinct canonical forms of constant-form systems.} The difference between the normal form $\Xi_P^{\pm}$ and $\Xi_P^0$ lies in the existence or not of an equilibrium point: the control $w=0$ defines an equilibrium of $\Xi_P^0$ while there are no equilibria for $\Xi_P^{\pm}$. 
\begin{proof}\leavevmode
    \begin{enumerate}[label=\textit{(\roman*)}]
        \item \textbf{Sufficiency}: Consider a parabolic system $\Xi_P=(A,B,C)${, with structure functions $(\mu_0,\mu_1,\gamma_0,\gamma_1)$,} and assume that relation \cref{eq:thm_m1_normalize_quadratic_submanifold_gamma1_0} holds. {Introduce the reparametrisation $\alpha=-\gamma_1$ and $\beta$ given as a solution of the following system of equations
        \begin{align}\label{prf:eq:thm_m1_normalize_quadratic_submanifold_gamma1_0:system_beta}
            \left\{\begin{array}{rl}
                \frac{1}{\beta}\dL{A}{\beta} &= -\mu_1 \\ 
                \frac{1}{\beta}\dL{B}{\beta} &= \frac{1}{2}\left(\mu_0-2\gamma_1\mu_1-2\dL{A}{\gamma_1}\right)
            \end{array}\right..
        \end{align}
        \noindent
        This system, rewritten for $\ln(\beta)$, admits solutions since the integrability condition
        \begin{align}\tag{\ref*{prf:eq:thm_m1_normalize_quadratic_submanifold_gamma1_0:system_beta}'}\label{prf:eq:thm_m1_normalize_quadratic_submanifold_gamma1_0:system_beta:prime}
            \dL{A}{\dL{B}{\ln\beta}}-\dL{B}{\dL{A}{\ln\beta}} = \dL{\lb{A}{B}}{\ln \beta} &= \mu_0\dL{A}{\ln\beta} + \mu_1 \dL{B}{\ln\beta}
        \end{align}
        \noindent
        takes the form 
        \begin{align}\tag{\ref*{prf:eq:thm_m1_normalize_quadratic_submanifold_gamma1_0:system_beta}''}\label{prf:eq:thm_m1_normalize_quadratic_submanifold_gamma1_0:system_beta:prime:prime}
            \frac{1}{2}\dL{A}{\mu_0}-\gamma_1\dL{A}{\mu_1}-\dLk{A}{2}{\gamma_1}-\mu_1\dL{A}{\gamma_1}+\dL{B}{\mu_1}  &= -\frac{1}{2}\mu_0\mu_1-\gamma_1(\mu_1)^2-\mu_1\dL{A}{\gamma_1}
        \end{align}
        and is guaranteed by condition \cref{eq:thm_m1_normalize_quadratic_submanifold_gamma1_0}. Consider the system $\tilde{\Xi}_P=(\tilde{A},\tilde{B},\tilde{C})$ obtained by the above defined reparametrisation $(\alpha,\beta)$ and using relations \cref{eq:action_alpha_beta_on_gamma_nf2d,eq:reparam_m1_action_on_ab_bracket} we deduce that the structure functions $(\tilde{\mu}_0,\tilde{\mu}_1,\tilde{\gamma}_0,\tilde{\gamma}_1)$ of $\tilde{\Xi}_P$ satisfy $\tilde{\mu}_0=\tilde{\mu}_1=0$ and $\tilde{\gamma}_1=0$. By choosing local coordinates $(z,y)$ such that $(\tilde{A},\tilde{B})=\left(\vec{z},\vec{y}\right)$, we obtain the system $\Xi_P'$ with $c_1=\tilde{\gamma}_1=0${, i.e. a strongly-flat parabolic system $\Xi_P'''$}.
        }
        %
        \addtocounter{enumi}{-1}
        \item \textbf{Necessity}: Assume that $\Xi_P=(A,B,C)$ is feedback equivalent to $\tilde{\Xi}_P$ of the form $\Xi_P'$ with $c_1\equiv 0$, via $\phi$ and $(\alpha,\beta)$. Then for $\tilde{\Xi}_P$ we have $\tilde{\mu}_0=\tilde{\mu}_1=\tilde{\gamma}_1=0$. First, by \cref{eq:action_alpha_beta_on_gamma_nf2d} we obtain $\alpha = -\gamma_1$ and by  \cref{eq:reparam_m1_action_on_ab_bracket} we obtain {that $\beta$ satisfies} the relations given by system \cref{prf:eq:thm_m1_normalize_quadratic_submanifold_gamma1_0:system_beta} above. 
        Therefore, by computing the integrability condition \cref{prf:eq:thm_m1_normalize_quadratic_submanifold_gamma1_0:system_beta:prime}, equivalently given by \cref{prf:eq:thm_m1_normalize_quadratic_submanifold_gamma1_0:system_beta:prime:prime}, we conclude relation \cref{eq:thm_m1_normalize_quadratic_submanifold_gamma1_0}.
        \item \textbf{Sufficiency}: Consider a parabolic system $\Xi_P=(A,B,C)${, with structure functions $(\mu_0,\mu_1,\gamma_0,\gamma_1)$,} and assume that $\Gamma_P=\gamma_0+(\gamma_1)^2\neq0$ and \cref{eq:thm_m1_normalize_quadratic_submanifold_g1_0_g0_pm} holds. {We follow the same reasoning as in the sufficiency of statement \cref{thm:m1_normalisation_quadratic_system_gamma_1_0:1}. Introduce a reparametrisation $(\alpha,\beta)$, where $\alpha=-\gamma_1$ and $\beta$ is a solution of the system of equations \cref{prf:eq:thm_m1_normalize_quadratic_submanifold_gamma1_0:system_beta}. To assure the existence of $\beta$, we have to fulfil the integrability condition of \cref{prf:eq:thm_m1_normalize_quadratic_submanifold_gamma1_0:system_beta}, which is \cref{prf:eq:thm_m1_normalize_quadratic_submanifold_gamma1_0:system_beta:prime}, equivalently \cref{prf:eq:thm_m1_normalize_quadratic_submanifold_gamma1_0:system_beta:prime:prime}. To this end, we differentiate the second condition of \cref{eq:thm_m1_normalize_quadratic_submanifold_g1_0_g0_pm} along $A$ and use $\dL{A}{\dL{B}{\cdot}}=\dL{B}{\dL{A}{\cdot}}+\mu_0\dL{A}{\cdot}+\mu_1\dL{B}{\cdot}$ to conclude \cref{prf:eq:thm_m1_normalize_quadratic_submanifold_gamma1_0:system_beta:prime:prime}. Consider the system $\tilde{\Xi}_P=(\tilde{A},\tilde{B},\tilde{C})$ obtained by the above defined reparametrisation $(\alpha,\beta)$ and using relations \cref{eq:action_alpha_beta_on_gamma_nf2d,eq:reparam_m1_action_on_ab_bracket}, we deduce that the structure functions $(\tilde{\mu}_0,\tilde{\mu}_1,\tilde{\gamma}_0,\tilde{\gamma}_1)$ of $\tilde{\Xi}_P=(\tilde{A},\tilde{B},\tilde{C})$ satisfy $\tilde{\mu}_0=\tilde{\mu}_1=\tilde{\gamma}_1=0$. Therefore $\tilde{\Gamma}_P=\tilde{\gamma}_0$, for which condition \cref{eq:thm_m1_normalize_quadratic_submanifold_g1_0_g0_pm} yields $\dL{\tilde{A}}{\tilde{\Gamma}_P}=\dL{\tilde{B}}{\tilde{\Gamma}_P}=0$, implying that $\tilde{\Gamma}_P$ is constant (we still have $\tilde{\Gamma}_P\neq0$ since $\beta^2\tilde{\Gamma}_P=\Gamma_P$). Introduce coordinates $(\tilde{z},\tilde{y})$ such that $\tilde{A}=\vec{\tilde{z}}$ and $\tilde{B}=\vec{\tilde{y}}$, in which the system takes the form (recall that $\tilde{\gamma}_1=0$)
        \begin{align*}
            \left\{\begin{array}{cl}
                \dot{\tilde{z}} &= \tilde{w}^2 + c_0 \\
                \dot{\tilde{y}} &= \tilde{w} 
            \end{array}\right.,
        \end{align*}
        \noindent
        with $c_0\in\RR^*$. Finally, defining new coordinates $(z,y)$ by $z=\frac{\tilde{z}}{|c_0|}$ and $y=\frac{\tilde{y}}{\sqrt{|c_0|}}$, and reparametrising by $w=\frac{\tilde{w}}{\sqrt{|c_0|}}$, yields the normal form $\Xi_P^{\pm}$.} 

        \addtocounter{enumi}{-1}
        \item\textbf{Necessity}: Assume that $\Xi_P$, whose structure functions are  $(\mu_0,\mu_1,\gamma_0,\gamma_1)$ and $\Gamma_P=\gamma_0+\gamma_1^2$, is feedback equivalent, via $\phi$ and $(\alpha,\beta)$, to $\tilde{\Xi}_P$ of the form $\Xi_P'$ with $(c_0,c_1)\in\RR^2$ satisfying $c_0+c_1^2\neq 0$. For $\tilde{\Xi}_P$ we have $\tilde{\mu}_0=\tilde{\mu}_1=0$ and $\tilde{\gamma}_0=c_0$, $\tilde{\gamma}_1=c_1$, hence $\tilde{\Gamma}_P=c_0+c_1^2\neq 0$ implying $\Gamma_P\neq0$ since $\beta^2\tilde{\Gamma}_P=\Gamma_P$. By {relation} \cref{eq:reparam_m1_action_on_ab_bracket},  we obtain that $\beta$ satisfies the relations of system \cref{prf:eq:thm_m1_normalize_quadratic_submanifold_gamma1_0:system_beta}. 
        Differentiating $\Gamma_P=\beta^2\tilde{\Gamma}_P$ along $A$ we deduce
        \begin{align*}
            \dL{A}{\Gamma_P} &= \dL{A}{\beta^2\tilde{\Gamma}_P}= 2\tilde{\Gamma}_P\beta\dL{A}{\beta} = -2\tilde{\Gamma}_P\beta^2\mu_1 = -2\mu_1\Gamma_P,
        \end{align*}
        \noindent
        giving the first relation of \cref{eq:thm_m1_normalize_quadratic_submanifold_g1_0_g0_pm}. A similar computation, by taking the derivative of $\Gamma_P=\beta^2\tilde{\Gamma}_P$ along $B$, implies the second relation of \cref{eq:thm_m1_normalize_quadratic_submanifold_g1_0_g0_pm}. 
        \item The proof of that statement is a special case of the proof of item \cref{thm:m1_normalisation_quadratic_system_gamma_1_0:1} with the additional condition $\Gamma_P\equiv 0$.\\
        \textbf{Sufficiency}: 
        Use the proof of the sufficiency of item \cref{thm:m1_normalisation_quadratic_system_gamma_1_0:1} to bring the system $\Xi_P$ into $\Xi_P'''$. For the latter form we have $\Gamma_P=c_0(x)$, hence $c_0(x)\equiv 0$ (due to $\beta^2\tilde{\Gamma}_P=\Gamma_P$ and assumption $\Gamma_P\equiv0$) and we obtain the normal form $\Xi_P^0$.\\
        \textbf{Necessity}: Assume that $\Xi_P$, whose structure functions are $(\mu_0,\mu_1,\gamma_0,\gamma_1)$ and  $\Gamma_P=\gamma_0+\gamma_1^2$, is feedback equivalent, via $\phi$ and $(\alpha,\beta)$, to $\tilde{\Xi}_P$ of the form $\Xi_P^0$ (which is, actually, $\Xi_P'$ with $c_0\equiv c_1\equiv0$). For that system we have $\tilde{\mu}_0=\tilde{\mu}_1=0$ and $\tilde{\Gamma}_P\equiv 0$ and since $\Gamma_P$ is transformed under $(\alpha,\beta)$ by $\beta^2\tilde{\Gamma}_P=\Gamma_P$, we obtain the necessity of $\Gamma_P\equiv 0$. The necessity of \cref{eq:thm_m1_normalize_quadratic_submanifold_gamma1_0} is deduced from the necessity part of statement \cref{thm:m1_normalisation_quadratic_system_gamma_1_0:1}.
    \end{enumerate}
\qed\end{proof}
Observe that item \cref{thm:m1_normalisation_quadratic_system_gamma_1_0:2} {of the above theorem} does not explicitly require condition \cref{eq:thm_m1_normalize_quadratic_submanifold_gamma1_0}, while the normal form $\Xi_P^{\pm}$ satisfies $c_1\equiv 0$ and hence that condition has to be hidden in \cref{eq:thm_m1_normalize_quadratic_submanifold_g1_0_g0_pm}. Indeed, this can be observed by differentiating $\Gamma_P$ along $\lb{A}{B}$ and using the constraint \cref{eq:thm_m1_normalize_quadratic_submanifold_g1_0_g0_pm}, which after a short computation gives condition \cref{eq:thm_m1_normalize_quadratic_submanifold_gamma1_0}.
\begin{remark}[Interpretation of the conditions]
     We now give a tangible interpretation of our conditions. To this end, consider the system $\Xi_P'$ for which we have $\mu_0=\mu_1=0$, $\gamma_0=c_0(x)$, $\gamma_1=c_1(x)$ and thus $\Gamma_P(x)=c_0(x)+(c_1(x))^2$. First, condition \cref{eq:thm_m1_normalize_quadratic_submanifold_gamma1_0} implies $\frac{\partial^2 c_1}{\partial z^2}=0$, that is, $c_1$ is affine with respect to $z$, namely, $c_1(x)=c_1^0(y)z+c_1^1(y)$, and thus, $\Gamma_P(x)$ is given by $c_0(x)+(c_1^0(y)z+c_1^1(y))^2$. This means that if a {weakly-flat} system $\Xi_P'$ is feedback equivalent to {a strongly-flat system} $\Xi_P'''$, then it is parametrized by 3 smooth functions, two of them being functions of $y$ only, and it has the following form
    \begin{align*}
        \left\{\begin{array}{cl}
            \dot{z} &= w^2 + \Gamma_P(x) - \left(c_1^0(y) z + c_1^1(y)\right)^2 \\
            \dot{y} &=w  + c_1^0(y) z + c_1^1(y)
        \end{array}\right..
    \end{align*}
    \noindent
    By additionally applying the first equation of \cref{eq:thm_m1_normalize_quadratic_submanifold_g1_0_g0_pm}, we obtain $\frac{\partial c_0}{\partial z}+2c_1^0(y)^2z+2c_1^0(y)c_1^1(y)=0$ and thus $c_0(x)$ is a polynomial of degree 2 in $z$, related to $c_1(x)$ by $c_0(x) = -(c_1(x))^2 + c_2(y)$, for an arbitrary smooth function $c_2(y)$. We now have $\Gamma_P =\Gamma_P(y)= c_2(y)$ and we use the second equation of \cref{eq:thm_m1_normalize_quadratic_submanifold_g1_0_g0_pm}. Thus, we get $\Gamma_P(y)=G\exp\left(-2\int c_1^1(y)\,dy\right)$, with $G\in\RR$. To summarise, any system $\Xi_P'$ satisfying \cref{eq:thm_m1_normalize_quadratic_submanifold_gamma1_0} and \cref{eq:thm_m1_normalize_quadratic_submanifold_g1_0_g0_pm}{, i.e. equivalent to a constant-form parabolic system,} is parametrized by two arbitrary smooth functions of $y$ and a constant $G\in\RR$ and is expressed by the form 
    \begin{align*}
        \left\{\begin{array}{cl}
            \dot{z} &= w^2 + \Gamma_P(y)-\left(c_1^0(y) z + c_1^1(y)\right)^2 \\
            \dot{y} &=w  + c_1^0(y) z + c_1^1(y)
        \end{array}\right.,
    \end{align*}
    \noindent 
    where $\Gamma_P(y)=G\exp\left(-2\int c_1^1(y)\,dy\right)$. Finally, $\Xi_P'$ (satisfying \cref{eq:thm_m1_normalize_quadratic_submanifold_gamma1_0} and \cref{eq:thm_m1_normalize_quadratic_submanifold_g1_0_g0_pm}) is feedback equivalent to $\Xi_P^0$ if and only if $G=0$, and is feedback equivalent to $\Xi_P^+$, respectively to $\Xi_P^-$, if $G>0$, resp. $G<0$. The distinction between the three normal forms comes from the sign of $\Gamma_P$, which is thus a discrete invariant (and that sign is dictated by the value of the constant $G$).
\end{remark}
\section{Conclusions} 
In this paper, we studied {the characterisation and the classification problem of three dimensional submanifolds of the tangent bundle of a smooth surface. We showed that the equivalence of submanifolds is reflected in the equivalence, under feedback transformations, of their first and second extensions treated as control-nonlinear and control-affine systems, respectively. We gave a complete characterisation of non-degenerate quadratic submanifolds and proposed a classification of the regular ones, namely the classes of elliptic, hyperbolic, and parabolic, submanifolds. To achieve our characterisation results, we introduced the novel class of quadratic control-affine systems and we gave a characterisation of that class in terms of relations between structure functions. Using our characterisation, we identified the subclasses of elliptic, hyperbolic, and parabolic control-affine systems (second extensions of elliptic, hyperbolic, and parabolic submanifolds). Moreover, we constructed a normal-form for all quadratisable systems and thus gave a normal form of quadratic submanifolds, that may smoothly pass from the elliptic to the hyperbolic classes.} 
Finally, working within the class of control-nonlinear systems subject to a {regular quadratic}\st{conic} nonholonomic constraint, we exhibited several normal forms of elliptic, hyperbolic, and parabolic systems. In particular we highlighted a connection between the Gaussian curvature of a well-defined metric and the existence of a commutative frame for elliptic and hyperbolic systems. Normal forms include systems without functional parameters. {As a consequence of our classification of quadratic systems, we obtained a classification of elliptic, hyperbolic, and parabolic submanifolds.}

In future works we plan to extend our results to higher dimensional {quadratic nonholonomic}\st{conic} constraints, {in particular, we will generalise our results for parabolic systems to the case of control-nonlinear systems (with the state-space being an $n\geq3$ dimensional manifold) subject to paraboloid nonholomic constraints}\st{as well as we will study systems on $3$-dimensional manifolds subject to a parabolic-elliptic or parabolic-hyperbolic constraints}. Another interesting problem is to characterise control-nonlinear systems subject to any algebraic nonholonomic constraint. For instance, in order to generalise our results for parabolic systems, one can study polynomial systems, that is, systems subject to the nonholonomic constraint $\dot{z}-\sum_{i=0}^da_i(x)\dot{y}^i=0$.

\begin{acknowledgements}
The authors would like to thank anonymous reviewers whose suggestions and remarks have substantially improved the presentation of the paper.
\end{acknowledgements}


\bibliographystyle{references/spmpsci}
\bibliography{references/bibliography.bib}

\appendix

\appendix
\section{Simultaneous rectification of two rank \texorpdfstring{$1$}{1} distributions in \texorpdfstring{$\RR^2$}{R²}}\label{appendix:m1:rectification_2_distrib_rk_1}
We show that locally on $\RR^2$ we can simultaneously rectify two involutive complementary distributions. To this end, consider $\AAA=\distrib{A}$ and $\BBB=\distrib{B}$, two distributions of rank $1$ satisfying $\AAA(x_0)\oplus\BBB(x_0)=T_{x_0}\RR^2$. Choose two functions $\phi$ and $\psi$ satisfying $d\phi(x_0)\neq0$, $d\psi(x_0)\neq0$, $d\phi\perp\BBB$, and $d\psi\perp\AAA$. Obviously, $\left(d\phi\wedge d\psi\right)(x_0)\neq0$ and set $(z,y)=(\phi,\psi)$. So in the $(z,y)$-coordinates, we have $\AAA=\distrib{\vec{z}}$ and $\BBB=\distrib{\vec{y}}$.

\section{Resolution of equation \texorpdfstring{\cref{eq:proof_thm_m1_nf_equation_fundamental}}{\ref*{eq:proof_thm_m1_nf_equation_fundamental}}} \label{appendix:m1_solution_equation_quadratizable}
Consider $f$ a smooth function of $(z,y,w)$, we will abbreviate $(z,y)$ in $x$ to shorten the notation. We show how to construct all {smooth} solutions around $(x_0,w_0)$ of the equation
\begin{align*}
    \frac{\partial^3f}{\partial w}&= \structfunctB(x)\frac{\partial f}{\partial w}.
\end{align*}
\noindent
First, we will find an expression for $\frac{\partial f}{\partial w}$ that will be then integrated to obtain the desired form. Consider the following linear system of PDEs: 
\begin{align*}
    \begin{pmatrix}\frac{\partial f_1}{\partial w} \\ \frac{\partial f_2}{\partial w}\end{pmatrix} &= \begin{pmatrix}0&1\\ \structfunctB(x)&0\end{pmatrix} \begin{pmatrix}f_1\\ f_2\end{pmatrix} 
\end{align*}
\noindent
given for the functions $f_1=\frac{\partial f}{\partial w}$ and $f_2=\frac{\partial^2 f}{\partial w^2}$. Solutions of this system {(interpreted as a system of ODEs with $x$ being a parameter)} are expressed by the exponential of the matrix $\begin{pmatrix}0&w\\ \structfunctB(x)w&0\end{pmatrix}$ given by the formula
\begin{align*}
    \exp\begin{pmatrix}0&w\\ \structfunctB(x)w&0\end{pmatrix}&= \sum_{k=0}^{+\infty}\frac{w^{2k+1}\structfunctB(x)^k}{(2k+1)!}\begin{pmatrix}0&1\\ \structfunctB(x)&0\end{pmatrix}+\sum_{k=0}^{+\infty}\frac{w^{2k}\structfunctB(x)^k}{(2k)!}\begin{pmatrix}1&0\\ 0&1\end{pmatrix}
\end{align*}
\noindent
obtained by expressing the power series of the exponential and by regrouping the terms of odd and even degree. Denote ${b}(x)=f_1(x,w_0)$ and  ${a}(x)=f_2(x,w_0)$, thus we obtain
\begin{align*}
    \frac{\partial f}{\partial w}=f_1(x,w) &= {a}(x)\sum_{k=0}^{+\infty}\frac{(w-w_0)^{2k+1}}{(2k+1)!}\structfunctB(x)^k + {b}(x) \sum_{k=0}^{+\infty}\frac{(w-w_0)^{2k}}{(2k)!}\structfunctB(x)^k.
\end{align*}
\noindent
Integration of this expression yields
\begin{align*}
    f(x,w) &= {a}(x)\sum_{k=0}^{+\infty}\frac{(w-w_0)^{2k+2}}{(2k+2)!}\structfunctB(x)^k + {b}(x) \sum_{k=0}^{+\infty}\frac{(w-w_0)^{2k+1}}{(2k+1)!}\structfunctB(x)^k +c(x).
\end{align*}

\section{Detailed computation of the proof of  \texorpdfstring{\cref{thm:m1_normal_form_quadratizable_systems}}{Theorem \ref*{thm:m1_normal_form_quadratizable_systems}}}
We extend the computation made in the proof of \cref{thm:m1_normal_form_quadratizable_systems}.
\subsection{Resolution of equation \texorpdfstring{\cref{eq:prf-all-quadratization-submanifold-rho}}{\ref*{eq:prf-all-quadratization-submanifold-rho}}}\label{appendix:m1_solution_lineard_equation}
We show how to, locally around {$0\in\RR^3$}\st{$\xi_0=(x_0,w_0)\in\RR^3$}, solve the equation $\rho''-2\rho\rho'+\frac{4}{9}\rho^3=0$, where $\rho=\rho(x,w)$ and the derivatives are taken with respect to $w$. Take the new unknown function $R(x,w)=\exp\left(-\frac{2}{3}\int_{0}^w\rho(x,t)\,\diff t\right)$ which satisfies $R(x,0)\neq0$ and 
\begin{align*}
    R'&=-\frac{2}{3}\rho R,\quad R'' =-\frac{2}{3}R\left(\rho'-\frac{2}{3}\rho^2\right),\quad 
    R'''= -\frac{2}{3}R\left(\rho''-2\rho\rho'+\frac{4}{9}\rho^3\right)=0.
\end{align*}
\noindent
Thus, $R(x,w)=a(x)w^2+b(x)w+c(x)$ yielding $\rho =-\frac{3}{2}\frac{R'}{R}=-\frac{3}{2}\frac{2aw+b}{aw^2+bw+c}$. Since $R(x,0)=c(x_0)\neq0$, thus taking $d(x)=\frac{a}{c}$ and $e(x)=\frac{b}{c}$ we obtain 
\begin{align*}
    \rho(x,w)=-\frac{3}{2}\frac{2d(x)w+e(x)}{d(x)w^2+e(x)w+1}.
\end{align*}

\subsection{Smooth form of \texorpdfstring{$h'$}{h'}}\label{appendix:m1_smooth_derivative_h}

We integrate $h''(x,w)=a(x)(d(x)w^2+e(x)w+1)^{-3/2}$ to obtain a smooth expression of $h'(x,w)$ around {$0\in\RR^3$}\st{$(x_0,0)\in\RR^3$}. Recall that we denote $p=p(x,w)=dw^2+ew+1$ and $\Delta=\Delta(x)=e^2-4d(x)$. We have
\begin{align*}
    h'(x,w)&=a(x) \int \frac{1}{p(x,w)^{3/2}}\,dw = \frac{-2a(2dw+e)}{\Delta\sqrt{p}}+\bar{b}(x).
\end{align*}
\noindent
Since $h'(x,0)=b(x)$ is smooth, we have $\frac{-2ae}{\Delta}+\bar{b}=b(x)$, where $b(x)$ is a smooth function around $0$. We can then derive a smooth closed form expression of $h'(x,w)$:
\begin{align*}
    h'(x,w)&=\frac{-2a(2dw+e)}{\Delta\sqrt{p}}+\frac{2ae}{\Delta}+b = \frac{-2a}{\Delta\sqrt{p}}\left(2dw+e-e\sqrt{p}\right)+b, \\
    &= \frac{-2a}{\Delta\sqrt{p}(ew + 2 +2\sqrt{p})}\left(2dw+e-e\sqrt{p}\right)(ew+2+2\sqrt{p}) +b, \\
    &= \frac{-2a}{\Delta\sqrt{p}(ew+2+2\sqrt{p})}\left(2dew^2+4dw+4dw\sqrt{p}+e^2w+2e+2e\sqrt{p}-e^2w\sqrt{p}-2e\sqrt{p}-2ep\right)+b, \\
    &= \frac{-2a}{\Delta\sqrt{p}(ew+2+2\sqrt{p})}\left(w\sqrt{p}(4d-e^2)+4dw+2dew^2+2e+e^2w-2edw^2-2e^2w-2e\right)+b, \\
    &= \frac{-2a}{\Delta\sqrt{p}(ew+2+2\sqrt{p})}\left(w\sqrt{p}(4d-e^2)+w(4d-e^2)\right)+b, \\
    &=\frac{2aw}{\sqrt{p}(ew+2+2\sqrt{p})}\left(\sqrt{p}+1\right)+b.
\end{align*}

\section{Details of the computations of \texorpdfstring{\cref{lem:m1:reparam_struct_funct_ellip_hyper}}{Lemma \ref*{lem:m1:reparam_struct_funct_ellip_hyper}}}\label{apdx:m1:reparam_struct_funct_ellip_hyper}
Denoting $c_{\ttiny{E}}(x)=\cos(x)$, $c_{\ttiny{H}}(w)=\cosh(w)$, $s_{\ttiny{E}}(w)=\sin(w)$, and $s_{\ttiny{H}}(w)=\sinh(w)$ and starting from the system 
\begin{align*}
    \dot{x}=(A,B)\begin{pmatrix} c_{\ttiny{EH}}(w)\\s_{\ttiny{EH}}(w)\end{pmatrix}+(A,B)\begin{pmatrix} \gamma_0\\\gamma_1\end{pmatrix},
\end{align*}
\noindent
we apply a reparametrisation $w=\tilde{w}+\alpha(x)$: 
\begin{align*}
    \dot{x}&=(A,B)\bar{R}_{\ttiny{EH}}(\pm\alpha)\begin{pmatrix} c_{\ttiny{EH}}(\tilde{w})\\s_{\ttiny{EH}}(\tilde{w})\end{pmatrix}+(A,B)\begin{pmatrix} \gamma_0\\\gamma_1\end{pmatrix}, \\ 
    &=(\tilde{A},\tilde{B})\begin{pmatrix} c_{\ttiny{EH}}(\tilde{w})\\s_{\ttiny{EH}}(\tilde{w})\end{pmatrix}+(\tilde{A},\tilde{B})\bar{R}_{\ttiny{EH}}(\pm\alpha)^{-1}\begin{pmatrix} \gamma_0\\\gamma_1\end{pmatrix}.
\end{align*}
\noindent
This yields $\left(\tilde{\gamma_0}, \tilde{\gamma_1}\right) = \left( \gamma_0, \gamma_1\right)\bar{R}_{\ttiny{EH}}(\pm\alpha)^{-T}$, and by the definition of $\bar{R}_{\ttiny{EH}}(\alpha)$ we have, 
\begin{align*}
    \bar{R}_{\ttiny{E}}(\alpha)^{-T}=\bar{R}_{\ttiny{E}}(-\alpha)^T=\bar{R}_{\ttiny{E}}(\alpha),\quad \textrm{and} \quad \bar{R}_{\ttiny{H}}(-\alpha)^{-T}=\bar{R}_{\ttiny{H}}(-\alpha)^{-1}=\bar{R}_{\ttiny{H}}(\alpha).
\end{align*}
\noindent
Next, computing separately in the elliptic and hyperbolic cases, we have
\begin{align*}
    \lb{\tilde{A}}{\tilde{B}}&= (\mu_0\mp\dL{A}{\alpha})A+(\mu_1-\dL{B}{\alpha})B, \\
    (A,B)\bar{R}_{{\ttiny{E}}}(\pm\alpha)\begin{pmatrix}\tilde{\mu}_0 \\ \tilde{\mu}_1 \end{pmatrix} &=(A,B) \begin{pmatrix}\mu_0\mp\dL{A}{\alpha} \\ \mu_1-\dL{B}{\alpha} \end{pmatrix}.
\end{align*}
\noindent
Thus
\begin{align*}
    \begin{pmatrix}\tilde{\mu}_0 , \tilde{\mu}_1 \end{pmatrix}=\begin{pmatrix}\mu_0\mp\dL{A}{\alpha} , \mu_1-\dL{B}{\alpha} \end{pmatrix}\bar{R}_{{\ttiny{E}}}(\pm\alpha)^{-T}
\end{align*}
\noindent
and relation \cref{eq:lem:m1:reparam_struct_funct_EH} follows.

\section{Gaussian curvature for metric given in terms of vector fields}\label{apdx:m1:gaussian_curvature}
Consider a $2$-dimensional manifold $\XXX$ and two smooth vector fields $A$, and $B$ satisfying $A\wedge B\neq0$. Construct the (pseudo)-Riemanian metric $\textsf{g}_{\pm}$ defined by 
\begin{align*}
    \textsf{g}_{\pm}(A,A)=1,\quad\textsf{g}_{\pm}(B,B)=\pm1,\quad\textrm{and}\quad \textsf{g}_{\pm}(A,B)=0.
\end{align*}
\noindent
We will give a formula for the Gaussian curvature of $\textsf{g}_{\pm}$ in terms of the structure functions $(\mu_0,\mu_1)$ uniquely defined by $\lb{A}{B}=\mu_0A+\mu_1B$. We will use the following formula for the covariant derivative 
\begin{align*}
    \nabla_{E_i}E_j=\frac{1}{2}\sum_k\left(\textsf{g}_{\pm}(\lb{E_i}{E_j},E_k)-\textsf{g}_{\pm}(\lb{E_i}{E_k},E_j)-\textsf{g}_{\pm}(\lb{E_j}{E_k},E_i)\right)E_k
\end{align*}
\noindent
for $E_i,E_j,E_k\in\{A,B\}$, and the following formula for the Gaussian curvature of a $2$-dimensional manifold
\begin{align*}
    \kappa_{\pm}=\frac{\textsf{g}_{\pm}\left( (\nabla_B\nabla_A-\nabla_A\nabla_B+\nabla_{\lb{A}{B}})A, B\right)}{\det(\textsf{g}_{\pm})}.
\end{align*}
Computing, we have
\begin{align*}
    \nabla_AA&= -\mu_0\, B,\quad \nabla_AB=\mu_0\,A,\quad  \nabla_BA=\mp\mu_1\,B,\quad \nabla_BB=\pm\mu_1\,A.
\end{align*}
\noindent
Then we can deduce
\begin{align*}
    \nabla_B\nabla_AA&= -\dL{B}{\mu_0}B\mp\mu_0\mu_1A,\quad \nabla_A\nabla_BA= \mp\dL{A}{\mu_1}B\mp\mu_0\mu_1A,\\ 
    \nabla_{\lb{A}{B}}A&=\mu_0\nabla_AA+\mu_1\nabla_BA = -(\mu_0)^2B\mp(\mu_1)^2B.
\end{align*}
Thus
\begin{align*}
    \kappa_{\pm}&=\pm\textsf{g}_{\pm}\left( \left(-\dL{B}{\mu_0}\pm\dL{A}{\mu_1}-(\mu_0)^2\mp(\mu_1)^2\right)B,B\right), \\
    &=-\dL{B}{\mu_0}\pm\dL{A}{\mu_1}-(\mu_0)^2\mp(\mu_1)^2.
\end{align*}        
\noindent
{Therefore, we have obtained the expression of relation \cref{eq:prop_existence_EH_commutative_frame}}.

\end{document}